\documentclass[oneside,11pt]{amsart}
\usepackage{amssymb,latexsym,amsmath,amsthm,enumitem,mathrsfs}
\usepackage[margin=1in]{geometry}
\usepackage{hyperref}
\usepackage{fancyhdr}
\usepackage{color}
\usepackage{stmaryrd}
\usepackage{mathrsfs}
\usepackage{lineno}
\usepackage{diagbox}
\usepackage{array}
\usepackage{tikz}
\usetikzlibrary{arrows}
\usetikzlibrary{positioning}
\usepackage{float}
\usepackage{pgfplots}
\usepackage{soul}

\usepackage[normalem]{ulem}
\newcolumntype{C}[1]{>{\centering\let\newline\\\arraybackslash\hspace{0pt}}m{#1}}
\numberwithin{equation}{section}

\makeatletter
\def\@tocline#1#2#3#4#5#6#7{\relax
  \ifnum #1>\c@tocdepth % then omit
  \else
    \par \addpenalty\@secpenalty\addvspace{#2}%
    \begingroup \hyphenpenalty\@M
    \@ifempty{#4}{%
      \@tempdima\csname r@tocindent\number#1\endcsname\relax
    }{%
      \@tempdima#4\relax
    }%
    \parindent\z@ \leftskip#3\relax \advance\leftskip\@tempdima\relax
    \rightskip\@pnumwidth plus4em \parfillskip-\@pnumwidth
    #5\leavevmode\hskip-\@tempdima
      \ifcase #1
       \or\or \hskip 1em \or \hskip 2em \else \hskip 3em \fi%
      #6\nobreak\relax
    \hfill\hbox to\@pnumwidth{\@tocpagenum{#7}}\par% <---- \dotfill -> \hfill
    \nobreak
    \endgroup
  \fi}
\makeatother

\newtheorem{letterthm}{Theorem}

\newtheorem{letterconjecture}[letterthm]{Conjecture}

\newtheorem{thm}{Theorem}[section]
\newtheorem*{thm*}{Theorem}
\newtheorem{prop}[thm]{Proposition}
\newtheorem{lemma}[thm]{Lemma}
\newtheorem{property}[thm]{Property}
\newtheorem{cor}[thm]{Corollary}
\newtheorem{conjecture}[thm]{Conjecture}
\newtheorem{condition}[thm]{Condition}
\newtheorem{question}[thm]{Question}

\newtheorem{dfn}[thm]{Definition}
\newtheorem*{dfn*}{Definition}

\theoremstyle{remark}
\newtheorem{remark}[thm]{Remark}

\newcommand{\ep}{\varepsilon}

\newcommand{\con}{\equiv}

\newcommand{\ndiv}{\nmid}
\newcommand{\modd}[1]{\; ( \text{mod} \; #1)}
\newcommand{\bstack}[2]{#1 \atop #2}
\newcommand{\tstack}[3]{#1\atop {#2 \atop #3}}

\newcommand{\maps}{\rightarrow}
\newcommand{\intersect}{\cap}

\newcommand{\Union}{\bigcup}

\newcommand{\al}{\alpha}
\newcommand{\be}{\beta}

\newcommand{\ga}{\gamma}
\newcommand{\del}{\delta}
\newcommand{\Del}{\Delta}
\newcommand{\ka}{\kappa}
\newcommand{\om}{\omega}

\newcommand{\sig}{\sigma}
\newcommand{\lam}{\lambda}

\newcommand{\Dbf}{\mathbf{D}}

\newcommand{\tbf}{{\bf t}}

\newcommand{\beq}{\begin{equation}}
\newcommand{\eeq}{\end{equation}}

\usepackage{amsmath, amssymb, amscd, verbatim, xspace,amsthm}
\usepackage{latexsym, epsfig, color}
\usepackage[OT2,T1]{fontenc}
\DeclareSymbolFont{cyrletters}{OT2}{wncyr}{m}{n}
\DeclareMathSymbol{\Sha}{\mathalpha}{cyrletters}{"58}

\newcommand{\ra}{\rightarrow}

\newcommand\Hom{\operatorname{Hom}}

\newcommand\Aut{\operatorname{Aut}}

\newcommand\isom{\simeq}
\newcommand\sub{\subset}

\renewcommand\O{\mathcal{O}}

\newcommand\bq{\begin{equation}}
\newcommand\eq{\end{equation}}

\newcommand{\Q}{\mathbb{Q}}
\newcommand{\R}{\mathbb{R}}

\newcommand{\C}{\mathbb{C}}

\newcommand{\Z}{\mathbb{Z}}

\newcommand{\Cbf}{\mathbf{C}}

\newcommand{\pfrak}{\mathfrak{p}}

\newcommand{\qfrak}{\mathfrak{q}}
\newcommand{\Cl}{\mathrm{Cl}}
\newcommand{\RI}{\mathcal{O}}

\newcommand{\Ocal}{\mathcal{O}}

\newcommand{\Cscr}{\mathscr{C}}

\newcommand{\Fscr}{\mathscr{F}}
\newcommand{\Iscr}{\mathscr{I}}
\newcommand{\Lscr}{\mathscr{L}}

\newcommand{\Gal}{\mathrm{Gal}}
\newcommand{\Disc}{\mathrm{Disc}\,}
\newcommand{\Li}{\mathrm{Li}}
\newcommand{\Cond}{\mathrm{Cond}}
\newcommand{\GL}{\mathrm{GL}}

\newcommand{\reg}{\mathrm{reg}}

\newcommand{\Nm}{\mathrm{Nm}}

\newcommand{\Ker}{\mathrm{Ker}}

\begin{document}

\title[Pierce, Turnage-Butterbaugh, Wood]{An effective Chebotarev density theorem \\
for families of number fields, \\
with an application to $\ell$-torsion in class groups}

\author[Pierce]{Lillian B. Pierce}
\address{Department of Mathematics, Duke University, 120 Science Drive, Durham, NC 27708 USA\\
and Institute for Advanced Study, 1 Einstein Drive, Princeton, NJ 08540 USA}
\email{pierce@math.duke.edu}

\author[Turnage-Butterbaugh]{Caroline L. Turnage-Butterbaugh}
\address{Department of Mathematics and Statistics, Carleton College, 1 North College Street, Northfield, MN 55057 USA}
\email{cturnageb@carleton.edu}

\author[Wood]{Melanie Matchett Wood}
\address{Department of Mathematics, 
University of California, Berkeley,
970 Evans Hall \# 3840,
Berkeley, CA 94720-3840 USA} 
\email{mmwood@berkeley.edu}

\keywords{Chebotarev, automorphic L-functions, counting number fields, class groups}
\subjclass[2010]{11R29, 11R42, 11R45, 11N75}

%11R29, %Class numbers, class groups, discriminants
%11R42 % zeta and L-functions of number fields
%11R45 % density theorems for algebraic number theory, for field counting
%11N75   	Applications of automorphic functions and forms to multiplicative problems 

\begin{abstract}

We prove a new effective Chebotarev density theorem for Galois extensions $L/\Q$ that allows one to count small primes (even as small as an arbitrarily small power of the discriminant of $L$); this theorem holds for the Galois closures of ``almost all'' number fields that lie in an appropriate family of field extensions. Previously, applying Chebotarev in such small ranges required assuming the Generalized Riemann Hypothesis. 
The error term in this new Chebotarev density theorem also avoids the effect of an exceptional zero of the Dedekind zeta function of $L$, without assuming GRH.
We give many different ``appropriate families,'' including families of arbitrarily large degree. To do this, we first prove a new effective Chebotarev density theorem that requires a zero-free region of the Dedekind zeta function. Then we prove that almost all number fields in our families yield such a zero-free region. The   innovation that allows us to achieve this is a
 delicate new method for controlling zeroes of certain families of \emph{non-cuspidal} $L$-functions. This builds on, and greatly generalizes the applicability of, work of Kowalski and Michel on the average density of zeroes of a family of \emph{cuspidal} $L$-functions.
A surprising feature of this new method, which we expect will have independent interest, is that we control the number of zeroes in the family of $L$-functions by bounding the number of certain associated fields with fixed discriminant. As an application of the new Chebotarev density theorem, we prove the first nontrivial upper bounds for $\ell$-torsion in class groups, for all integers $\ell \geq 1$, applicable to infinite families of fields of arbitrarily large degree. 
 \end{abstract}

\maketitle 
\setcounter{tocdepth}{0}

\section{Overview}
 In this paper, we give unconditional effective Chebotarev density theorems for almost all number fields in certain families of fields, of a strength that previously required the assumption of GRH.  We achieve this by a new method to control zeroes of non-cuspidal $L$-functions in families, and we give applications including the first non-trivial bounds on $\ell$-torsion for all $\ell \geq 1$ in class groups in infinite families of fields of arbitrarily large degree.  Our method requires only crude bounds on the number of fields in our families,  allowing us to treat families of arbitrarily high degree and more general families than in 
  \cite{EPW16}, which gives $\ell$-torsion bounds as a result of very precise counting of the families.

\subsection{Historical introduction}\label{sec_intro}
For any fixed number field $k$ and Galois extension $L/k$ of number fields, consider the counting function of prime ideals of bounded norm in $\Ocal_k$ and specified splitting type in $L$, defined by
\beq\label{pi_dfn}
 \pi_{\mathscr{C}}(x,L/k) := \# \{ \pfrak \subseteq \Ocal_k: \pfrak  \; \text{unramified in $L$}, \left[ \frac{L/k}{\pfrak} \right] = \Cscr, \Nm_{k/\Q} \pfrak \leq x \},
 \eeq
in which $\left[ \frac{L/k}{\pfrak} \right]$ is the Artin symbol and $\Cscr$ is any fixed conjugacy class in $\Gal(L/k)$. 
A central goal is to prove an asymptotic for $\pi_{\mathscr{C}}(x,L/k)$ that is valid for $x$ as small as possible (relative to the absolute discriminant of the number field $L$), which is a regime in which many of the most interesting applications arise.
  The celebrated Chebotarev density theorem \cite{Che26} provides the main term in the asymptotic,
  \beq\label{pi_C_asymp0}
    \pi_{\mathscr{C}}(x,L/k)  \sim \frac{|\Cscr |}{|G|} \mathrm{Li} (x),
    \eeq
 as $x \maps \infty$, where $\Gal(L/k) = G$ and $\mathrm{Li}(x) = \int_2^x dt/\log t.$ When $L=k=\Q$, this is the familiar Prime Number Theorem for $\pi(x)$; when $L=k$, this is the Prime Ideal Theorem, counting prime ideals $\pfrak \subset \Ocal_k$ with $\Nm_{k/\Q} \pfrak \leq x$; when $k=\Q$ and $L=\Q(e^{2\pi i /q})$, this provides Dirichlet's theorem, counting rational primes $p \con a \modd{q}$ with $p \leq x$, for any $(a,q)=1$.
  
An effective Chebotarev theorem, conditional on GRH, was proved by Lagarias and Odlyzko (with an improvement by Serre). Given any field extension $F/\Q$ we let $n_F=[F:\Q]$ and set $D_F = |\Disc F/\Q|$.

\begin{letterthm}[Conditional on GRH, {\cite[Theorem 1.1]{LagOdl75}, \cite[Th\'{e}or\`{e}me 4]{Ser82}}]\label{thm_LO}
There exists an effectively computable absolute constant $C_0>0$ such that for any Galois extension $L/k$ of number fields, if GRH holds for the Dedekind zeta function $\zeta_L$ and  $G:=\Gal(L/k)$, then for any fixed conjugacy class $\mathscr{C} \subseteq G$ and every $x \geq 2$,
\[ \left| \pi_{\mathscr{C}} (x,L/k) - \frac{|\mathscr{C}|}{|G|} \Li(x)  \right| \leq C_0 \frac{|\mathscr{C}|}{|G|} x^{1/2} \log (D_L x^{n_L}).\]
\end{letterthm}
Lagarias and Odlyzko also proved an unconditional result: 
\begin{letterthm}[{\cite[Corollary 1.3]{LagOdl75}}]\label{thm_LagOdl_uncond}
There exist effectively computable absolute constants $C_1,C_2>0$ such that the following holds.
Let $L/k$ be a Galois extension of number fields with $G:=\Gal(L/k)$. If $n_L>1$ then $\zeta_L(s)$ has at most one zero $s=\sig +it$ in the region
\beq\label{LO_region}
 \sig \ge 1-(4\log D_L)^{-1}, \qquad |t| \le (4\log D_L)^{-1}.
\eeq
This exceptional zero, denoted $\be_0$ if it exists, is real and simple. For all  
$x \geq \exp(10n_L(\log D_L)^2),
$
\beq\label{LO_asymp}
  \left| \pi_{\mathscr{C}}(x,L/k) - \frac{|\mathscr{C}|}{|G|} \Li(x) \right| \leq \frac{|\mathscr{C}|}{|G|} \Li (x^{\be_0}) + C_1 x \exp(-C_2 n_L^{-1/2} (\log x)^{1/2} ),
  \eeq
with the understanding that the $\be_0$ term is present only if $\be_0$ exists.
\end{letterthm}

Theorem \ref{thm_LO} holds for all $x \geq 2$.
Theorem \ref{thm_LagOdl_uncond} requires at least that $x \geq D_L^{10n_L}$, a power of the discriminant that is too large for many applications.  Consequently, citations of the Lagarias-Odlyzko work often use Theorem \ref{thm_LO} and are hence conditional on GRH.
Recent unconditional work that considers lower or upper bounds for $\pi_\Cscr(x,L/k)$ instead of asymptotics also leads to thresholds for $x$ that are too large for certain applications. For example,  \cite{ThoZam16a}, \cite[Eqn. 1.6]{ThoZam16b} prove lower bounds for $\pi_{\Cscr}(x,L/k)$ that require $x$ to be as large as a relatively large power of $D_L$; upper bounds for $\pi_\Cscr(x,L/k)$  in the classic work \cite[Thm. 1.4]{LMO79}  require $x \geq C \exp \{( \log D_L)( \log\log D_L) (\log \log \log D_L)\}$, for some constant $C$, with  improvements e.g. in \cite{ThoZam16b, Deb16x}.

\subsection{New results I: effective Chebotarev theorems}\label{sec_new_results_Cheb}
We prove a new effective Chebotarev theorem that includes two breakthroughs: we remove the term corresponding to the exceptional zero in (\ref{LO_asymp}),
and simultaneously we obtain an  asymptotic with an effective error term, which in particular holds   for $x$ as small as $D_L^\del$ for any small fixed $\del>0$ (for $D_L$ sufficiently large). Both aspects are critical to applications such as our new bound for $\ell$-torsion in class groups.
It is unlikely that we could accomplish these goals for \emph{all} fields  without proving something significant toward GRH; instead, we prove that within appropriate families of fields, ``almost all'' of the fields satisfy such an effective Chebotarev theorem.

We first state an inexplicit, general version of our result for a ``family'' $\Fscr(G)$ of fields (precise quantitative statements appear in Theorems  \ref{thm_cheb_main_quant_pt1}, \ref{thm_cheb_main_quant_pt2}, \ref{thm_cheb_main_quant_pt3}, \ref{thm_cheb_main_quant_pt2'}, \ref{thm_cheb_main_quant_pt4}, and Corollary \ref{cor_primes_families}  of  \S \ref{sec_quant_thms}).  
By a family $\Fscr(G)$ we mean a set of degree $n$ extensions $K/\Q$ with corresponding Galois closures $\tilde{K}/\Q$ having $\Gal(\tilde{K}/\Q)\simeq G$ for a fixed transitive subgroup $G \subseteq S_n$. We use $\Fscr(G;X)$ to denote those fields $K \in \Fscr(G)$ with $D_K \leq X$. We also use Vinogradov's notation: $A \ll B$ denotes that there exists a constant $C$ such that $|A| \leq C B,$ and $A\ll_\kappa B$ denotes that  $C$ may depend on $\kappa$.
\begin{thm}\label{thm_gen}
Fix an appropriate family $\Fscr(G)$ (described explicitly in \S \ref{sec_appropriate}), and constants $A \geq 2$ and $\ep>0$.  Then there exist constants $0<\tau < \be$  and $\ka_1,\ka_2,\ka_3>0$, such that for all $X \geq 1$ we have
$|\Fscr(G;X)| \gg X^\be,$ and   aside from at most 
$\ll_{\Fscr,A,\ep}X^{\tau + \ep}$ possible exceptions, each field $K \in \Fscr(G;X)$ 
 has the property that for every conjugacy class $\mathscr{C} \subseteq G$,
\beq\label{final_pi_quant_thm_gen}
 \left| \pi_{\mathscr{C}} (x,\tilde{K}/\Q) - \frac{|\mathscr{C}|}{|G|} \mathrm{Li}(x) \right| \leq  \frac{|\mathscr{C}|}{|G|}\frac{x}{(\log x)^A},
 \eeq
for all
 \beq\label{x_lower_bound_00_quant_thm_gen}
 x \geq \ka_1 \exp \{ \ka_2 (\log \log (D_{\tilde{K}}^{\ka_3}))^{5/3}(\log\log \log (D_{\tilde{K}}^{2}))^{1/3}\}  .
 \eeq
\end{thm}

In comparison to Theorem \ref{thm_LagOdl_uncond}, for each field to which this result applies, this theorem removes the effect of the possible exceptional zero on the error term, and holds for $x$ as small as an arbitrarily small power of $D_{\tilde{K}}$ (and hence of $D_K$), capabilities critical for many applications.

\subsubsection{The appropriate families of fields}\label{sec_appropriate}
 In general, we construct a set (or ``family'') of fields as follows. 
For a number field $k$, we  let
\[ Z_n(k, G;X)  = \{ K/k : K\sub\bar{\Q},  \deg{K/k} = n, \Gal(\tilde{K}/k) \simeq G, \Nm_{k/\Q} \Disc K/k \leq X \},\]
 where $\tilde{K}$ is the Galois closure of $K$ over $k$,  the Galois group is considered as a permutation group on the $n$ embeddings of $K$ in $\overline{\Q}$, and the isomorphism with $G$ is one of permutation groups.
We let $Z_n(k, G)=Z_n(k,G;\infty)$.
For our main results we will work over $\Q$, and study families of the form $Z_n^{\Iscr}(\Q,G;X)$, defined to be the subset of those fields $K\in Z_n(\Q,G;X)$ such that for each rational prime $p$ that is tamely ramified in $K$ (i.e. those $p$ not dividing any of the exponents of their factorization into prime ideals in the ring of integers of $K$), the  inertia group in $\Gal(\tilde{K}/k)$  of every prime ideal $\wp$ of $\tilde{K}$ dividing $p$ 
is  generated by an element of $\Iscr$, where $\Iscr$ specifies one or more conjugacy classes in $G$.  
The use of ramification restrictions will play a large role in our method of proof.

The most general families we treat are degree $n$ extensions with square-free discriminant, which are a positive proportion of all degree $n$ fields for $n \leq 5$, and conjecturally so for $n \geq 6$. (These families are recorded in entries \eqref{E:s3}, \eqref{E:s4}, \eqref{E:s5} in the lists below; square-free discriminant corresponds to $\Iscr$ being transpositions, as explained in  \S \ref{sec_sym_gp}.)  We give further examples to show the range of the method.
We prove, unconditionally, that Theorem \ref{thm_gen} applies to the following families $Z_n^\Iscr(\Q,G)$ of fields:
\begin{enumerate}
\item $G$ a cyclic group of order $n \geq 2$, with $\Iscr$ comprised of all generators of $G$ (equivalently every rational prime that is tamely ramified in $K$ is totally ramified).
\item  $n=p$ an odd prime, $G = D_p$ the order $2p$ dihedral  group of symmetries of a regular $p$-gon, $\Iscr$ being the conjugacy class of order $2$ elements. \label{E:Dp}
\item  $n=3,$ $G \simeq S_3$, $\Iscr$ is transpositions. \label{E:s3}
\item  $n=4,$ $G \simeq  S_4$, $\Iscr$ is transpositions. \label{E:s4}
\item  $n=4$, $G \simeq A_4$, $\Iscr$ the two conjugacy classes in $A_4$ of order $3$ elements.
\end{enumerate}
This is the content of Theorem  \ref{thm_cheb_main_quant_pt1}. 
Note for family \eqref{E:Dp} that asymptotic counting of the fields is essentially equivalent to knowing the exact average size of $p$-torsion of class groups of quadratic fields (and thus is open and very difficult).  Our method does not require this counting.

We furthermore prove, conditional on the strong Artin conjecture and (in some cases certain hypotheses for counting number fields), that Theorem \ref{thm_gen} applies to the following families of fields:
\begin{enumerate}
\setcounter{enumi}{5}
\item  $n \geq 5,$ $G \simeq  S_n$, $\Iscr$ is transpositions (Theorems  \ref{thm_cheb_main_quant_pt2} and \ref{thm_cheb_main_quant_pt3}).  \label{E:s5}
\item $n \geq 5$, $G \simeq A_n$, no ramification restriction (Theorem \ref{thm_cheb_main_quant_pt2'}).
\item $G \subseteq S_n$ a transitive simple group, no ramification restriction (Theorem \ref{thm_cheb_main_quant_pt4}).
\end{enumerate}
In addition, in Corollary \ref{cor_primes_families} we record quantitative results for counting certain types of primes.

\subsubsection{The proof strategy}
To describe our strategy to prove Theorem \ref{thm_gen} we define the notion of a $\del$-exceptional field:
\begin{property}[$\del$-exceptional field]\label{prop_del_exceptional}
For a fixed $0 < \del < 1/2$, a number field $K$ is $\del$-\emph{exceptional}  precisely when the Dedekind zeta function of the Galois closure $\tilde{K}$ of $K$ over $\Q$ has the property that $\zeta_{\tilde{K}}(s)/\zeta(s)$ has a zero in the region 
$[1-\del,1] \times [-(\log D_{\tilde{K}})^{2/\del}, (\log D_{\tilde{K}})^{2/\del}].$
\end{property}
(Under GRH, no field is $\del$-exceptional for any $0 \leq \del <1/2$.) 
Our first step toward Theorem \ref{thm_gen} is to prove the following (for a quantitative version over any fixed number field $k$, see Theorem \ref{thm_assume_zerofree1}):
\begin{thm}[Effective Chebotarev for non-$\delta$-exceptional fields]\label{thm_del_ex}
 For every integer $n \geq 1$, and every transitive group $G \subseteq S_n$, for every $A \geq 2$ and every $0 < \del \leq 1/(2A)$, there exist real numbers $D_0, \kappa_1,\kappa_2,\kappa_3$  (depending on $\del,n,A$)
  such that the following holds: for any extension $K/\Q$ with $\Gal(\tilde{K}/\Q) \simeq G$ such that
 $D_{\tilde{K}} \geq D_0$ and $K/\Q$ is not $\delta$-exceptional,
we have that for any conjugacy class $\mathscr{C} \subseteq G$, (\ref{final_pi_quant_thm_gen}) holds for all  $x$ satisfying (\ref{x_lower_bound_00_quant_thm_gen}). 
\end{thm}
Theorem \ref{thm_gen} relies on the following crucial step: we prove that within appropriate families, for sufficiently small $\del$,  almost all fields are not $\del$-exceptional.\footnote{See also  \S \ref{sec_app_Cheb_quad} on an unconditional approach to rule out exceptional zeroes in the standard zero-free region for $\zeta_{L}(s)$, and thus remove the $\be_0$ term in (\ref{LO_asymp}), for extensions with no quadratic subfields.  However, that approach does not rule out the extensions being $\del$-exceptional, and in particular, does not lead to an effective Chebotarev theorem that can count primes small enough  for our purposes.
} 
 We achieve this by developing a new method for controlling zeroes of certain families of non-cuspidal $L$-functions. 
 Previously, work of Kowalski and Michel \cite{KowMic02} provided density results for zeroes within appropriate families of cuspidal $L$-functions.  But we require zero-free regions for Dedekind zeta functions of Galois fields, and these correspond (in some cases conjecturally) to  automorphic $L$-functions that are \emph{not cuspidal}.
 This restriction of \cite{KowMic02} to the cuspidal case has been a significant barrier in many previous applications (such as an effective prime ideal theorem in \cite{ChoKim14a}, or \cite{ChoKim13b}; see Remark \ref{remark_Cho_Kim}).
We expect that our new approach to proving density results for zeroes in a family of non-cuspidal $L$-functions will have many further applications.

Precisely, let $G$ be a fixed transitive subgroup of $S_n$ and let $\rho_0, \rho_1, \ldots, \rho_s$ denote the irreducible representations of $G$, with $\rho_0$ being the trivial representation. Then for each  $K \in Z_n(\Q,G;X)$,  we may write $\zeta_{\tilde{K}}(s)$ as a product of Artin $L$-functions
\beq\label{zeta_prod_intro}
 \zeta_{\tilde{K}}(s) = \zeta(s) \prod_{i=1}^s L(s,\rho_i, \tilde{K}/\Q)^{\mathrm{dim}\rho_i}. \eeq
In particular,   consider a set $\Fscr(X)$ of fields $K \in Z_n(\Q,G;X)$ with distinct Galois closures $\tilde{K}$ over $\Q$, and denote the set of Galois closures by $\tilde{\Fscr}(X)$.  For each field $\tilde{K} \in \tilde{\Fscr}(X)$ and each representation $\rho_j$, there is an associated cuspidal automorphic representation $\pi_{\tilde{K},j}$ of $\GL(m_j)/\Q$ (in some cases conditional on the Strong Artin Conjecture), and then $L(s,\pi_{\tilde{K},j}) = L(s,\rho_j,{\tilde{K}}/\Q)$. For each $1 \leq j \leq s$, we let $\Lscr_j(X)$ denote the set of cuspidal automorphic representations $\pi_{\tilde{K},j}$ of $\GL(m_j)/\Q$ associated   to the fields $\tilde{K} \in \tilde{\mathscr{F}}(X)$ and the  representation $\rho_j$.
We   show using  \cite{KowMic02} that for each $j$, $\Lscr_j(X)$ has the property that aside from at most a possible small ``bad'' exceptional subset, each representation $\pi \in \Lscr_j(X)$ is such that its associated $L$-function $L(s,\pi)$ is zero-free in an appropriate region. (Of course, if GRH is true, there are no such exceptional $L$-functions, but we are working without GRH.)  In order to deduce that amongst the Dedekind zeta functions $\zeta_{\tilde{K}}(s)$ for $\tilde{K} \in \tilde{\Fscr}(X)$, almost all of them also possess this zero-free region, we need to build up the products as in  (\ref{zeta_prod_intro}), and  we need to understand the following question: given a representation $\pi \in \Lscr_j(X)$ (i.e. possibly a ``bad'' exceptional representation), how many fields $\tilde{K} \in \tilde{\Fscr}(X)$ can have the property that $L(s,\rho_j,\tilde{K}/\Q)=L(s,\pi)$? 

 This is  subtle, and relies on delicate properties of the families considered. 
At its heart the question is: for a fixed irreducible representation $\rho_i$ of $G$, for how many fields $K_1,K_2  \in \Fscr(X)  \subseteq Z_n(\Q,G;X)$ can we have $L(s,\rho_i, \tilde{K}_1/\Q)  = L(s,\rho_i, \tilde{K}_2/\Q)$?
We  transform this into a question of counting how many fields $K_1,K_2 \in \Fscr(X)$ have  fixed fields $\tilde{K}_1^{H} = \tilde{K}_2^H$, where $H = \Ker(\rho_i)$ (see Proposition \ref{prop_KluNic_new}). A challenge then appears:  for certain groups $G$, is it possible that such collisions can occur amongst a positive proportion of $K \in Z_n(\Q,G;X)$? (If so, a positive proportion of these fields could have $\zeta_{\tilde{K}}(s)$ containing a factor that is not zero-free in the desired region.)

For certain $G$, the answer is yes (see \S \ref{sec_rationale_ram}). 
In contrast, we show that for the groups $G$ and the corresponding families of fields we construct in our main theorems, the answer is no. 
Precisely, we  define each family $\Fscr(X) \subseteq Z_n(\Q,G;X)$ according to carefully chosen ramification restrictions on tamely ramified primes, and within these carefully constructed families we can transform the problem of counting fields that share a certain fixed field into a problem of counting number fields of degree $n$ with \emph{fixed} discriminant. This method of constructing families of fields so that we can  control the zeroes of associated $L$-functions by counting number fields is a key innovation of this paper.

Within our chosen families, by counting fields of fixed discriminant, we ultimately  show that such collisions of the fixed fields must be relatively rare. We can then prove that aside from at most a possible ``small'' exceptional subset of $\Fscr(X)$, each field has the property that its Dedekind zeta function  is  zero-free in an appropriate region.

In general, our approach can be seen as a new strategy that vastly generalizes the applicability of the result of Kowalski and Michel  to families  of automorphic $L$-functions corresponding not just to cuspidal automorphic representations but also to isobaric automorphic representations. We expect this new method will be relevant to other problems of interest.

\subsection{New results II: counting number fields}\label{sec_overview_counting}

Our new effective Chebotarev theorem for families of fields relies on quantitative counts for number fields in two ways. First, we must bound from above the number of fields in the family that have a fixed discriminant; second we must bound from below the  number of fields in the family with bounded discriminant.
In general, such questions lie in the arena of Malle's conjecture \cite{Mal02} and the Malle-Bhargava principle \cite[Section 10]{MWoo16}, and many questions remain open.

\begin{dfn}
Within a certain family $Z_n^{\Iscr}(\Q,G)$, 
we say a subset $E$ has \emph{density zero} if  for some $\gamma>0$ and some  $c_1>0$, for all $X \geq 1$,
\[|Z_n^{\Iscr}(\Q,G;X)|/ |E \intersect Z_n^{\Iscr}(\Q,G;X)|\geq c_1 X^{\gamma}.\]
  \end{dfn}
Each of our main results takes the form of an effective Chebotarev density theorem that holds for each field within a family of fields, except for fields belonging to a possible subfamily of density zero. In all cases, proving  an upper bound for $|E \intersect Z_n^{\Iscr}(\Q,G;X)|$  is a significant part of our new work; in many cases, proving a lower bound for $|Z_n^{\Iscr}(\Q,G;X)|$ is also a significant part of our new work.

For certain of the families of fields we consider, we prove the first recorded lower bounds. For example, we prove the following general result, from which we deduce the first lower bound in the literature for $|Z_n(\Q,A_n;X)|$ that grows like a power of $X$ (Theorem \ref{thm_counting_An_fields}).

\begin{thm}\label{thm_count_G_fields}
Fix an integer $n \geq 2$ and a transitive subgroup $G \sub S_n$.
Suppose $f(X,T_1,\dots,T_j) \in \Q[X,T_1,\ldots, T_j]$ is a regular polynomial of total degree $d$ in the $T_i$ and of degree $n$ in $X$ with transitive Galois group $G\sub S_n$ over $\Q(T_1,\dots,T_j)$.  
Then, for every $X \geq 1$ and every $\ep>0$,
\[  |Z_n (\Q,G;X) | \gg_{f,\ep}X^{\frac{1-|G|^{-1}}{d(2n-2)} - \ep}.\]
\end{thm}
Note that  a recent paper of D\`ebes \cite{Deb17} proves an analogous result for counting the degree $|G|$ Galois extensions in $\Z_{|G|}(\Q,G;X)$ rather than  the degree $n$ extensions we consider in Theorem \ref{thm_count_G_fields} (or equivalently, only in the case that $G$ is simply transitive).

In a different direction, as mentioned above, at a key step of extending the Kowalski-Michel zero density theorem to our setting (related to bounding $|E \intersect Z_n^{\Iscr}(\Q,G;X)|$ from above), we  require an upper bound for  how many fields have any given \emph{fixed} discriminant. To make things precise, 
 we define the following property (always defining extensions within $\overline{\Q}$):

\begin{property}[${\bf D}_n(G,\varpi)$]\label{property_D}
Let $n \geq 2$ be fixed and let $G$ be a fixed transitive subgroup of $S_n$. 
We say that property ${\bf D}_n(G,\varpi)$ holds if  for every fixed integer $D>1$ and for every $\ep>0$ there exist at most $\ll_{n,G,\ep} D^{\varpi+\ep}$ fields $K/\Q$ of degree $n$ and $\Gal(\tilde{K}/\Q) \simeq G$ such that $D_K=D$. Moreover, we say that property ${\bf D}_n(\varpi)$ holds if  for every fixed integer $D>1$ and for every $\ep>0$ there exist at most $\ll_{n,\ep} D^{\varpi+\ep}$ fields $K/\Q$ of degree $n$ such that $D_K=D$.
\end{property} 

For appropriate families, we can control $|E \intersect Z_n^{\Iscr}(\Q,G;X)|$ if we can prove Property $\Dbf_n(G,\varpi)$ for a sufficiently small $\varpi$. In particular we prove new results for $\Dbf_4(A_4,\varpi)$, $\Dbf_5(\varpi)$, and $\Dbf_p^\Iscr(D_p,\varpi)$, in the latter case assuming a certain ramification restriction.

The way Property $\Dbf_n(\varpi)$ arises in our work on families of automorphic $L$-functions appears to be completely new. But it is actually the subject of a well-known conjecture which  occupies a rather central role in number theory. Specifically, Duke \cite[\S 3]{Duk98} and Ellenberg and Venkatesh \cite[Conjecture 1.3]{EllVen05} conjecture:
 \begin{conjecture}[Discriminant Multiplicity Conjecture]\label{conj_DMC}
For each $n \geq 2$,  $\Dbf_n(0)$ holds.
\end{conjecture}

Of course,  $\Dbf_2(0)$ holds; for $n \geq 3$, much less is known, and results toward Conjecture \ref{conj_DMC} would have strong implications.
First, the ``pointwise'' counts encapsulated in Property $\Dbf_n(\varpi)$ relate to ``average'' counts for the number of extensions of degree $n$ with bounded discriminant. In one direction, this is trivial: Property $\Dbf_n(\varpi)$ immediately implies there are at most $\ll_{n,\ep} X^{1+\varpi + \ep}$ degree $n$ extensions of $\Q$ with discriminant at most $X$.  It may be surprising that there is also an implication in the other direction; this has been proved by Ellenberg and Venkatesh \cite[Prop. 4.8]{EllVen05}.

Second, questions about $\Dbf_n(\varpi)$ are directly connected to questions about $\ell$-torsion in class groups, for primes $\ell$. 
As just one example (see Duke \cite{Duk95}), quartic fields of fixed discriminant $-q$ ($q$ prime) can be explicitly classified by odd octahedral Galois representations of conductor $q$, and the number of such fields can be expressed as in \cite{Hei71} as an appropriate average of the number of 2-torsion elements in the class groups of cubic number fields of discriminant $-q$.  More generally, as noted in  \cite[p. 164]{EllVen05}, if Conjecture \ref{conj_DMC} holds (for all $n$), then it implies the main pointwise conjecture, Conjecture \ref{conj_class}, for upper bounds for $\ell$-torsion in class groups (for all $n,\ell$). The way we employ property $\Dbf_n(\varpi)$ in the present work is in some sense more efficient, since to study $\ell$-torsion (for all $\ell \geq 1$) in class groups of degree $n_0$ fields we only require information about $\Dbf_{n}(\varpi)$ for $n=n_0$, not for all $n$.

\subsection{New results III: applications}
We expect that the new effective Chebotarev theorems for families of fields will have many applications, and we exhibit two. First, we prove nontrivial bounds for $\ell$-torsion, for all integers $\ell \geq 1$, in class groups of ``almost all'' fields in each of the families to which our Chebotarev theorems apply (Theorem \ref{thm_class_main}). In many cases, these are the first ever nontrivial bounds for $\ell$-torsion, and in particular the first that apply to families of fields of arbitrarily large degree. As a second (related) application, we prove a result on the density of number fields with small generators, spurred by a question of Ruppert (Theorem \ref{thm_VaaWid_app}). Further applications will be described in later work.

\subsection{Organization of the paper}
In Part I, we state and prove the results we require for counting number fields, both with bounded discriminant and with fixed discriminant. In Part II, we turn to the Chebotarev theorems:
in \S \ref{sec_quant_thms} we state quantitative versions of all the effective Chebotarev theorems; in \S \ref{sec_ChebotarevProof} we  prove the quantitative version of Theorem \ref{thm_del_ex}, and in \S \ref{sec_zero_density} and \S \ref{sec_app_Dedekind} we prove  the quantitative versions of Theorem \ref{thm_gen}.
In Part III, we treat the two applications mentioned above.

\tableofcontents

\addcontentsline{toc}{part}{Part  I: Counting Number Fields}
\begin{center}
{\bf Part I: Counting Number Fields}\end{center}

%%%%%%%%%%%%%%%%%%%%%%%%%%%%%%%%%%%%
\section{Counting families of fields}\label{sec_counting_num_fields_proofs}
%%%%%%%%%%%%%%%%%%%%%%%%%%%%%%%%%%%%
As described in Section~\ref{sec_overview_counting}, we  require results counting number fields, and we prove those in this section.
Our principal concern is families of the form $Z_n^{\Iscr}(\Q,G;X)$, defined to be the subset of those fields $K\in Z_n(\Q,G;X)$ such that for each rational prime $p$ that is tamely ramified in $K$, an inertia group for $p$ is generated by an element of $\Iscr$.
 We  require an upper bound for  $|Z_n^{\Iscr}(\Q,G;X)|$, which can be an overestimate, a lower bound for $|Z_n^{\Iscr}(\Q,G;X)|$, which we aim to make as sharp as currently feasible, and  upper bounds on the number of fields in $Z_n^{\Iscr}(\Q,G)$ of discriminant  $D$.

\subsection{Cyclic fields}\label{sec_count_cyclic}
   The strategy for counting cyclic extensions goes back to Cohn \cite{Coh54}; see \cite{Mak85, Wri89,MWoo10,FLN15x} for results  counting abelian extensions of arbitrary degree.
   Let $G$ be cyclic of order $n \geq 2$ and let $g$ denote the smallest prime divisor of $n$.  Then we have (see, e.g. \cite{Wri89}) that 
\beq\label{cyclic_asymp}
 |Z_n(\Q,G;X)| \sim cX^{\frac{1}{n-n/g}}
 \eeq
for a certain constant $c=c(n)>0$.
 We require the following refinement:

\begin{prop}[Cyclic groups]\label{prop_counting_cyclic_fields}
Let $n \geq 2$ be fixed and let $G$ be a cyclic group of order $n$. Let $Z_n^{\Iscr}(\Q,G;X)$ count those fields $K\in Z_n(\Q,G;X)$ such that every rational prime that ramifies tamely in $K$ is totally ramified in $K$, that is, the inertia group is generated by an element that is of full order in $G$. Then there exists a constant $c_n>0$ such that
\beq\label{cyclic_full}
 |Z^{\Iscr}_n(\Q,G;X)| \sim c_n X^{\frac{1}{n-1}}.
\eeq
Furthermore, Property $\Dbf_n(G,0)$ holds.
\end{prop}

\begin{remark}\label{remark_cyclic_prime2}
If $|G|=n$ is prime then $1/(n-n/g) = 1/(n-1)$.  
However, when $|G|=n$ is not prime then $Z_n^\Iscr(\Q,G;X)$ is itself of density zero in $Z_n(\Q,G;X)$, by comparison of (\ref{cyclic_asymp}) and (\ref{cyclic_full}). 

\end{remark}

\begin{proof}
Let $a_1=1$ and for $m\geq 2$,  let $a_m$ be $|\Aut(G)|$ times the number of fields counted by $Z^{\Iscr}_n(\Q,G;X)$
with absolute discriminant $m$. We define a Dirichlet series $A(s):=\sum_{m\geq 1} a_m m^{-s}$, and by class field theory and now standard arguments we have
\begin{equation}\label{E:D}
A(s)  =P(s) \prod_{p \con 1 \modd{n}} ( 1 + \phi(n) p^{-(n-1)s}),
\end{equation}
where $P(s)$ is a product over $p|n$ of  polynomials in $p^{-s}$.  Briefly, by class field theory we are counting certain homomorphisms from the id\`ele class group to $G$, by \cite[Lemma 4.2]{MWoo10} we can replace  the id\`ele class group with a product of $p$-adic units, and then we can easily count the local homomorphisms
 (see, e.g. \cite[Section 4]{MWoo10}, for a similar analysis in a more difficult case).
 
When $a_m$ is non-zero, we have $a_m \leq C_n \phi(n)^{\om(m)}$ where $\om(m)$ is the number of distinct prime divisors of $m$ and $C_n$ is a constant depending only on $n$. (In particular, $C_n$ can be bounded above by the sum of the absolute values of all coefficients of the polynomial factors in the finite product $P(s)$.) Thus $a_m  \ll_{n,\ep} m^\ep$ for any $\ep>0$, proving Property $\Dbf_n(G;0)$.

For comparison to $A(s)$ we consider the product $B(s)$ over all Dirichlet characters defined modulo $n$, given for $\Re(s)>1$ by
$$
B(s)= \prod_{\chi} L(s,\chi) = \prod_{\chi} \prod_p ( 1 - \chi(p)p^{-s})^{-1} 
$$
which has a pole of order 1 at $s=1$ and otherwise may be analytically continued as a holomorphic function. 
Writing the Euler product as $ \prod_p \mu_p(s)^{-1}$, note that $\mu_p(s) = 1  -  \sum_{\chi} \chi(p) p^{-s} + O(p^{-2s})$; by orthogonality of characters, the coefficient $\sum_{\chi} \chi(p) = \phi(n)$ if $p \con 1 \modd{n}$ and zero otherwise. 
We can then check that $A(s)/B((n-1)s)$ is holomorphic in $\Re(s) > (2(n-1))^{-1}$.
Thus $A(s)$ has a  meromorphic continuation in $\Re(s) > (2(n-1))^{-1}$ with only a simple pole at $s=(n-1)^{-1}$; moreover $A(s)$ inherits a standard convexity estimate from $B(s)$ (see e.g. \cite[Lemma 5.2, Thm. 5.23]{IK}).  So, by the main term in a standard Tauberian theorem (see for example \cite[Thm. A.1]{CLT01} and \cite[Section 6.4]{Nar2000}), we have
$$
|Z^{\Iscr}_n(\Q,G;X)|=c_nX^{1/(n-1)}+o(X^{1/(n-1)}),
$$
for a certain constant $c_n$.
\end{proof}

\subsection{Dihedral groups $D_p$}

For $p$ an odd prime, let $D_p$ be the order $2p$ group of symmetries of the vertices of a regular $p$-gon. Kl\"{u}ners  \cite[Theorem 3.5]{Klu06} obtained the lower bound 
$
 |Z_p(\Q,D_p;X)| \gg X^{2/(p-1)}
$
predicted by Malle's conjecture \cite{Mal02}. Kl\"{u}ners also showed that Malle's conjectured upper bound $X^{2/(p-1)+\ep}$ follows from a special case of the Cohen-Lenstra heuristics  \cite[Thm. 2.5]{Klu06}, as well as proving \cite[Theorem 2.7]{Klu06} an  unconditional
 upper bound 
$
 |Z_p(\Q,D_p;X)| \ll_\ep X^{3/(p-1)+\ep}.
$
 This has recently been improved by  Cohen and Thorne \cite[Thm 1.1]{CohTho16x}, based on nontrivial bounds of \cite{EPW16} for averages of $\ell$-torsion over quadratic fields, to 
  \beq\label{Thorne_upper}
 |Z_p(\Q,D_p;X)|  \ll_\ep X^{\frac{3}{p-1} - \frac{1}{p(p-1)}+\ep}.
  \eeq

We require a lower bound that includes a ramification restriction. We let Property $\Dbf_p^\Iscr(D_p,\varpi)$ 
be the analog of Property $\Dbf_p(D_p,\varpi)$ in which we only count $D_p$-fields and with the ramification restriction $\Iscr$ for all tamely ramified primes.

\begin{prop}[Dihedral group $D_p$ of order $2p$]
\label{prop_counting_Dp_fields}
For $p$ an odd prime, let $D_p$ act on the $p$ vertices of the regular $p$-gon in the usual way, and let $Z_p^{\mathrm{\Iscr}}(\Q,D_p;X)$ count those fields $K\in Z_p(\Q,D_p;X)$ with the following ramification restriction $\Iscr$: every rational prime that ramifies tamely in $K$ has inertia group generated by an element in the conjugacy class $ [(2 \; p) (3  \; p-1)\cdots(\frac{p+1}{2} \; \frac{p+3}{2})]$ of reflections. Then 
$
 |Z^{\mathrm{\Iscr}}_p(\Q,D_p;X)| \gg_p X^{\frac{2}{p-1}}.
$

Further, $\Dbf_p^\Iscr(D_p,1/(p-1))$ holds.  More generally, if we know that for all quadratic fields $L$ we have $|\Cl_L[p]|=O_{p}(D_L^{b})$ for a certain exponent $b>0$, then $\Dbf_p^\Iscr(D_p,2b/(p-1))$ holds.

\end{prop}
Note: here we use the notation $\Cl_L[p]$ to denote the $p$-torsion subgroup of the class group $\Cl_L$ of the field $L/\Q$; see e.g. (\ref{torsion_subgroup_dfn}) for the definition.

\subsubsection{Proof of the upper bound}
Next we count degree $p$ $D_p$-fields with a fixed discriminant. 
We may trivially state that  $\Dbf_p(D_p,\varpi)$ holds with $\varpi= 3/(p-1) -1/p(p-1)$, by applying (\ref{Thorne_upper}). We improve on this by only counting fields with a fixed discriminant and using our additional ramification restriction. 

Let $K \in Z_p^\Iscr(\Q,D_p)$ be a degree $p$ $D_p$-field with absolute discriminant $D$.  Let $\tilde{K}$ be the Galois closure of $K$ and $L$ be the quadratic field inside $\tilde{K}$, so $\tilde{K}/L$ is a cyclic $p$ extension.  Our ramification restriction implies that $\tilde{K}/L$ is unramified except perhaps at primes dividing $2p$.  We have, by our ramification restriction, that $|\Disc K|=2^a p^b Q^{(p-1)/2}$, where $Q$ is square-free and relatively prime to $2p$.  Then $|\Disc L| = 2^{a'}p^{b'} Q$ for some $a',b'$ that are bounded in terms of $p$.  Thus, given $D$, there are a constant (in terms of $p$) possible quadratic fields $L$, and for each of them we will count the possible cyclic $p$ extensions $\tilde{K}/L$ that could arise.

Let $J_L$ be the id\`ele class group of $L$.  
For a finite place $v$ of $L$, let $\O_{v}$ be the elements of non-negative valuation in the completion $L_{v}$, and for an infinite place $v$ let $\O_v=L_v$.
From the exact sequence $\prod_v \O_v^* \ra J_L \ra \Cl_L\ra 1$ \cite[Ch. VI Prop. 1.3]{Neu} (where the product is over all places of $L$), and the left-exactness of $\Hom_{cts}( -,C_p)$, we have an exact sequence 
$$
1 \ra \Hom (\Cl_L,C_p) \ra \Hom_{cts} (J_L,C_p) \ra \Hom_{cts}( \prod_v \mathcal{O}_v^*,C_p),
$$
where we can take the product just above over finite places $v$ of $L$, since there are no continuous homormorphisms from $\R^*$ or $\C^*$ into $C_p$ for $p$ odd.  
Our desired $C_p$-extensions of $L$ correspond via class field theory to elements of $\Hom_{cts} (J_L,C_p)$ that for each $v \nmid 2p$ map $\mathcal{O}_v^*$ to the identity, since they are unramified at such $v$.  Thus the number of possible images in $\Hom_{cts} (\prod_v \mathcal{O}_v^*,C_p)$ for our desired elements of $\Hom_{cts} (J_L,C_p)$
is $|\Hom (\prod_{v\mid 2p} \mathcal{O}_v^*,C_p)|.$  The number of $v\mid 2p$ is at most $4$ since $L$ is quadratic.  Since $L_v$ is either $\Q_p$ or $\Q_2$ or quadratic over $\Q_p$ or $\Q_2$ (and there are only finitely many possibilities for the latter), the number of homomorphisms from $\mathcal{O}_v^*$ to $C_p$ for $v\mid 2p$ is bounded in terms of $p$.   Also, $|\Hom (\Cl_L,C_p)|=|\Cl_L[p]|$.  Note $\Disc L=O_p(|\Disc K|^{2/(p-1)})$. 
  So if we assume $|\Cl_L[p]|=O_{p}(|\Disc L|^{b})$, then the number of possible $\tilde{K}$, and thus the number of possible $K$, is $O_{p}(D^{2b/(p-1)})$.

\subsubsection{Proof of the lower bound}
Given a quadratic field $L$, if $\Cl_L[p]$ is non-trivial, 
class field theory gives an unramified cyclic degree $p$ extension $L'/L$.
The group $\Gal(L/\Q)=\langle \sigma \rangle$ acts on $\Cl_L$ by inversion (since for an ideal $\mathfrak{a}$ of $L$, we have that $\mathfrak{a}\sigma({\mathfrak{a}})$ is principal).
It follows that $L'/\Q$ is a degree $2p$ $D_p$-extension, with all inertia trivial or in a subgroup generated by a reflection. 

Now given an imaginary quadratic field $L$ with units $\pm 1$ such that $\Cl_L[p]$ is trivial and $p$ splits completely in $L$, we
 we will show by other means that we still can obtain a degree $2p$ $D_p$-extension $L'/\Q$ containing $L$, and  with our required ramification condition. 
We will construct a surjection $\phi$ from $J_L$ to the cyclic group $C_p$ of order $p$.  
Let $v_1,v_2$ be the two places of $L$ above $p$.    We let $\phi_{v_1}: \O_{v_1}^* \ra C_p$ be any surjection.  We let $\phi_{v_2} :  \O_{v_2}^* \ra C_p$ be defined by
$\phi_{v_2}(u) =\phi_{v_1}(\sigma(u))^{-1}$.  At every other place $v\ne v_1, v_2$, we let $\phi_v : \O_{v}^* \ra C_p$ be trivial.  
Then at each place $v$, we pick an element $\alpha_v\in L$ that has valuation $1$ at $v$ and valuation divisible by $p$ at all other places (which we can do since $\Cl_L[p]$ is trivial).   We extend $\phi_v$ to $\phi_v: L_v^* \ra C_p$ by letting $\phi_v(\alpha_v)=\prod_{w\ne v} \phi_w(\alpha_v)^{-1}$.  The $\phi_v$ combine to give a map $\phi : \prod_v L_v^* \ra C_p$, that is trivial on the diagonal embeddings of $p$th powers, the $\alpha_v$, and units.  These elements generate $L^*$ (since  $\Cl_L[p]$ is trivial), and so $\phi$ descends to a map $\phi: J_L \ra C_p$.  We can check that it follows from our definitions that $\phi(\sigma(x))=\phi(x)^{-1}$.   We recall from class field theory that the Artin map for $L$ is equivariant for the usual action of $\Gal(L/\Q)$
on $J_L$ and the action of $\Gal(L/\Q)$ on $\Gal(L^{ab}/L)$ given by conjugation by a lift in $\Gal(L^{ab}/\Q)$ \cite[Thm 11.5 (i)]{Tate67}.
So since $\ker \phi$ is $\Gal(L/\Q)$ invariant, it follows from Galois theory that the degree $p$ cyclic extension $L'$ of $L$ corresponding to $\phi$ (from class field theory) is actually Galois over $\Q$.  Since $p$ is odd, we have that $\Gal(L'/\Q)$ is a semi-direct product
$\Gal(L'/L)\rtimes \Gal(K/\Q)$, and the action of $\Gal(K/\Q)$ on the index $p$ subgroup $\Gal(L'/L)$ given above shows that $\Gal(L'/\Q)\isom D_p.$
 Since $L'/L$ has no tame ramification by choice of the $\phi_v|_{\O_{v}^*}$, all tame ramification of $L'/\Q$ has inertia in the subgroup of a reflection.

So for all but finitely many imaginary quadratic fields $L$ in which $p$ splits completely, we have constructed a degree $2p$ $D_p$-extension $L'/\Q$ containing $L$ with our required ramification condition, which in particular contains a degree $p$ $D_p$-extension $K$.  At primes $\ell\nmid 2p$ of $\Q$, the exponent of $\ell$ in $\Disc K$ is $(p-1)/2$ if $\ell$ is ramified in $L$ and $0$ otherwise.  So we have that $\Disc K$ is within a constant (depending on $p$) factor of $(\Disc L)^{(p-1)/2}$.  Since we have $\gg_p X$ of these quadratic fields $L$, we conclude we have $\gg_p X^{2/(p-1)}$ fields counted by $Z^{\mathrm{\Iscr}}_p(\Q,D_p;X)$.

 \subsection{Symmetric groups $S_n$}\label{sec_sym_gp}
Our work on $S_n$-fields requires understanding the size of $Z_n^\Iscr(\Q,S_n;X)$ with $\Iscr =[(1 \; 2)]$; this is equivalent to requiring that the tamely-ramified part of $D_K$ is square-free. This is a consequence of a standard fact (see Lemma \ref{lemma_power_p}) that $p$ is tamely ramified in $K$ with inertia group generated by a transposition if and only if $p\| D_K$. 
 We record for $n=3,4,5$, that by  work of Bhargava \cite[Theorem 1.3]{Bha14x}, 
\begin{equation}\label{E:Sncount}
 |Z_n^\Iscr(\Q,S_n;X)| \sim c_n X.
\end{equation}
By the asymptotic counts of $S_3$-fields due to  Davenport and Heilbronn \cite{DavHei71} and $S_4$-fields and $S_5$-fields due to Bhargava \cite{Bha05, Bha10a}, the fields in $Z_n^\Iscr(\Q,S_n)$ are a positive proportion of all $S_n$ fields for $n=3,4,5$. 
Moreover, it is conjectured by Malle \cite{Mal02} and Bhargava \cite{Bha07, Bha14x} that  asymptotics of order $X$
hold for $Z_n^\Iscr(\Q,S_n;X)$ and $Z_n(\Q,S_n;X)$
when $n\geq 6$.

For symmetric groups $S_n$ with $n \geq 6$, the best proven results are much weaker.
 For $n>2$, we have an upper bound of Ellenberg and Venkatesh \cite{EllVen06} on all degree $n$ number fields $Z_n(\Q)$,
 \beq\label{EllVen_upper}
 |Z_n(\Q; X)|\ll(\alpha_nX)^{\exp(C\sqrt{\log n})},
 \eeq
  where $\al_n$ is a constant depending only on $n$ and $C$ is an absolute constant.   
    The best known lower bound for $S_n$-fields is $ |Z_n(\Q, S_n; X)|\gg_n X^{1/2 + 1/n}$  by Bhargava, Shankar and Wang \cite[Thm. 1.3]{BSW16x}, and importantly for us, all of the fields they construct to deduce this new lower bound have square-free discriminant. 
 As a consequence,  for all $n \geq 6$ and  $\Iscr = [(1 \; 2)]$,
 \beq\label{Z_Sn_lower}
| Z_n^\Iscr(\Q,S_n;X) | \gg_{n} X^{1/2 +1/n}.
 \eeq

 We also require upper bounds on $S_n$-fields of a fixed discriminant, and we state the best known results here.
Ellenberg and Venkatesh \cite[p. 1]{EllVen07} prove Property $\Dbf_3(S_3,1/3)$. 
 Kl\"{u}ners \cite{Klu06b} proves Property $\mathbf{D}_4(S_4,1/2)$.
 From Bhargava's count for quintic fields, we may trivially deduce that $\Dbf_5(S_5,1)$ holds.
For our work, knowing $\Dbf_5(S_5,\varpi)$ for any $\varpi < 1$ would suffice, so we make the following simple observation:
\begin{prop}\label{prop_counting_S5_fields}
Property $\Dbf_5(\varpi)$ holds for $\varpi = 199/200$.
\end{prop}
This follows immediately from the power-saving count for quintic $S_5$-fields proved by Shankar and Tsimerman \cite{ShaTsi14} (see also the power-saving count for all quintic fields in \cite[Thm. 2.4]{EPW16}). Indeed, letting $Z_5(\Q;X)$ denote all quintic fields with $D_K \leq X$, we have a constant $c_{5a}>0$ such that
\[ |Z_5(\Q;X)| = c_{5a} X + O_\ep(X^{199/200 +\ep}) \]
for every $\ep>0$, so that upon differencing this for $X=D$ and $X=D-1$, Proposition \ref{prop_counting_S5_fields} follows.

\subsection{The alternating group $A_4$}\label{sec_A4_proof}

 For $A_4$, it is known by Baily \cite{Bai80} that the lower bound conjectured by Malle \cite{Mal02} holds,
$ |Z_4(\Q, A_4; X)|  \gg X^{1/2},
$
and by Wong  \cite{Won05} that a weaker  upper bound holds,
\beq\label{Wong_upper}
|Z_4(\Q, A_4; X)|  \ll X^{5/6+\ep}.
\eeq

We require a lower bound that includes a ramification restriction and an upper bound for fields of fixed discriminant.
\begin{prop}[Alternating group $A_4$]\label{prop_counting_A4_fields}
Let $Z_4^{\Iscr}(\Q,A_4;X)$ count those fields  $K\in Z_4(\Q,A_4;X)$ such that every rational prime that ramifies tamely in $K$ has inertia group generated by an element in either of the conjugacy classes 
$ \{(1\; 2\; 3), (1\; 3\; 4), (1\; 4\; 2), (2\; 4\; 3)\}$ or $\{(1\; 3\; 2), (1\; 4\; 3), (1\; 2\; 4), (2\; 3\; 4)\}.$
 Then 
$
|Z^{\Iscr}_4(\Q,A_4;X)| \gg X^{1/2}.
$
Moreover, $\Dbf_4(A_4,\varpi)$ holds for $\varpi=0.2784...$.
\end{prop}
 Property $\Dbf_4(A_4,\varpi)$  was previously known for $\varpi= 3/4$ due to Wong \cite[Thm. 6]{Won99b}, but this is not small enough for our purposes.

\subsubsection{The upper bound}
To show that $\Dbf_4(A_4,\varpi)$ holds with $\varpi=0.2784...$, we will apply Baily's connection \cite{Bai80} of $A_4$ fields to certain quadratic ray class characters of cyclic cubic fields, in combination with the bound on the $2$-torsion in class groups of cubic fields due to Bhargava,  Shankar,  Taniguchi,  Thorne, Tsimerman, and  Zhao \cite{BSTTTZ17}.
We thank Manjul Bhargava for suggesting this approach.

Let $K_4$ be a quartic $A_4$-field of discriminant $D$.  Let $K_3$ be the fixed  field of the subgroup of $A_4$ generated by $\{ (1 \; 2)(3 \; 4), (2 \; 3)(1 \;4)\}$, and note that $K_3$ is cyclic cubic.  We can check using Lemma~\ref{lemma_power_p}
that tame rational primes with inertia type in the conjugacy class of $(1 \; 2)(3 \; 4)$ appear squared in the discriminant of  $K_4$ and do not appear in the discriminant of $K_3$.  Similarly, 
tame rational primes with inertia type in the conjugacy class of $(1 \; 2 \; 3)$ appear squared in the discriminants of both $K_4$ and $K_3$.  
So $\Disc K_3 \mid 2^a3^b \Disc K_4$, for some absolute  positive integers $a,b$.

Let $K_6$ be one of the (conjugate) sextic subfields of the Galois closure of $K_4$.  Note that $K_4$ and $K_6$ have the same Galois closure, and so to count $K_4$ we may equivalently (up to a fixed constant) count the associated $K_6$.  By \cite[Lemmas 13 and 15]{Bai80} we have that $K_6=K_3(b^{1/2})$, where $b\in \O_{K_3}\setminus \{\Z \cup \O_{K_3}^2\}$ and $N_{K_3/\Q}(b)$ is a square rational integer.  We have that $N_{K_3/\Q}(\Disc(K_6/K_3))=\Disc K_4/\Disc K_3$ (see \cite[Lemma 11]{Bai80}).

Now, we sum over each divisor $d$ of $2^a3^bD$ the number of quartic $A_4$-fields $K_4$ of discriminant $D$ with $\Disc K_3=d$. 
 There are $O(2^{\omega(d)})$ cyclic cubic fields of discriminant $d$ \cite{Coh54}.  Given a fixed cyclic cubic field $K_3$ of discriminant $d$, for an upper bound, it suffices to bound the number of sextic fields of the form $K_3(b^{1/2})$, where $b\in \O_{K_3}\setminus \{\Z \cup \O_{K_3}^2\}$ and $N_{K_3/\Q}(b)$ is a square rational integer.  We do this following the argument in \cite[Lemma 10]{Bai80}. Such a sextic field corresponds to a quadratic ray class character of conductor $\mathfrak{d}$ with finite part
$\mathfrak{d^*}=\Disc(K_6/K_3)$, and  such a character is a product of a character on $(\O_{K_3}/\mathfrak{d^*})^\times$, a character on the class group of $K_3$, and a character on signature (see \cite[(4)]{Bai80}).  
  Baily \cite[Lemma 8]{Bai80} describes the possible forms of $\mathfrak{d}$, and in the proof of 
\cite[Lemma 9]{Bai80} gives a generating function for all the primitive quadratic characters on  $(\O_{K_3}/\mathfrak{d^*})^\times$.
From this it follows there are $O(3^{\omega(D/d)})$ choices of $\mathfrak{d^*}$ with characters on  $(\O_{K_3}/\mathfrak{d^*})^\times$ such that we will have $\Disc(K_6/K_3)=D/d$. Let $h_2(K_3)$ denote the size of the $2$-torsion subgroup of the class group of $K_3$. There are at most $h_2(K_3)$ class group characters, and $h_2(K_3)=O_{\ep}(d^{0.2784\dots+\ep})$ by \cite[Equation (4)]{BSTTTZ17}.  There are at most $8$ characters of signature, and so in conclusion, there are at most
$$
O_{\ep}\left( \sum_{d|D} 2^{\omega(d)}  3^{\omega(D/d)} d^{0.2784\ldots+\ep} \right)=O_{\ep}(D^{0.2784\ldots+\ep})
$$
quartic $A_4$-fields of discriminant $D$.

\subsubsection{The lower bound}
The lower bound on the number of quartic $A_4$-fields with our required ramification condition follows from the proof of
\cite[Theorem 3]{Bai80}.  As stated in line 2 of the proof of \cite[Lemma 16]{Bai80}, the degree $6$ fields $K_6$ constructed by Baily are unramified over the relevant degree $3$ cyclic field $K_3$, except perhaps at primes of $K_3$ dividing $2$.  These fields $K_6$ have Galois closure $K_{12}$ of degree 12 with Galois group $A_4$.  The fact that $K_6/K_3$ is unramified except at primes of $K_3$ dividing $2$ means that for each odd rational prime $p$
the inertia groups of $p$ in $\Gal(K_{12}/\Q)$ must be trivial or generated by a three-cycle.
The same holds for $p=2$ if the primes of $K_3$ that divide $2$ are unramified in $K_6/K_3$.  If a prime dividing $2$ is ramified in $K_6/K_3$, it is wildly ramified, and thus $2$ in wildly ramified in $K_{12}$.  

\subsection{The alternating groups $A_n$ and proof of Theorem \ref{thm_count_G_fields}}

The fact that for $n \geq 5$, $A_n$ is a simple group  will make a later part of our argument much simpler, but on the other hand we require a lower bound for the number of degree $n$ $A_n$-extensions of $\Q$ with bounded discriminant, which was not previously in the literature. We prove:
\begin{thm}[Alternating groups $A_n$, $n \geq 3$]\label{thm_counting_An_fields}
For each integer $n \geq 3$, there exists a real number $\be_n >0$ such that for all $X \geq 1$, for every $\ep>0$,
$|Z_n(\Q,A_n;X)| \gg_{n,\ep} X^{\be_n - \ep}.$
In fact we may take 
$ \be_n = (1-\frac{2}{n!})/(4n-4).
$
\end{thm}

We first observe that Theorem \ref{thm_count_G_fields} implies Theorem \ref{thm_counting_An_fields} when we specialize $G$ to $A_n$.  For each $n\geq 3$, Hilbert \cite{Hil92} gave polynomials $f(x,t)\in\Q[x,t]$ that have Galois group $A_n$ over $\Q(t)$ and are degree $n$ in $x$ and degree $2$ in $t$. (Hilbert in turn credits Hurwitz with the examples: see \cite[p. 125]{Hil92} for $n$ even and \cite[p. 126]{Hil92} for $n$ odd; see also \cite[Section 10.3]{Ser97}.)  Moreover, these same polynomials (by the same argument) have Galois group $A_n$ over $E(t)$, for any number field $E$, and thus their splitting fields do not contain a non-trivial finite extension of $\Q$ (i.e. they are regular).
Thus Theorem \ref{thm_count_G_fields} with $|G| = |A_n| = n!/2$, $j=1$, $m=n$ and $d=2$ verifies Theorem \ref{thm_counting_An_fields}.

We now prove Theorem \ref{thm_count_G_fields}; we thank Akshay Venkatesh and Manjul Bhargava for suggesting the approach we use, and for a number of helpful discussions.
The method of proof, in imprecise terms,  is as follows. Suppose that $f(x,\tbf)$ has Galois group $G$ over $\Q(\tbf)$, resulting in, say, $y$ different fields with Galois group $G$ as $\tbf$ varies over all integral tuples with coordinates at most $T$ in absolute value. Then by showing that $f(x,\tbf)f(x,\tbf')$ typically has Galois group $G \times G$ and very rarely has Galois group $G$ (which occurs when the fields provided by $f(x,\tbf)$ and $f(x,\tbf')$ collide), we will deduce that $f(x,\tbf)$ must have produced many different fields to begin with, that is, $y$ must grow at least like a small power of $T$.  
See also \cite[p. 137]{Ser97} for a hint at a similar philosophy applied to generating infinitely many $G$-extensions if one such extension is known.

In order to put this into action in precise terms, we require a quantitative version of the Hilbert irreducibility theorem, for which we cite \cite{CasDie16}:
\begin{letterthm}\label{thm_HIT}
Suppose  $f(X,T_1,\ldots, T_j) \in \Q[X,T_1,\ldots, T_j]$ is an irreducible polynomial with splitting field $K$ over $\Q(T_1,\ldots, T_j)$ such that $\Gal (K/\Q(T_1,\ldots, T_j))\isom G$. 
For any  subgroup $H \sub G$ set
\[ N_f(T;H) = \# \{ \tbf \in \Z^j : |\tbf |_\infty \leq T \; \text{and the splitting field of $f(X,\tbf)$ over $\Q$ has Galois group $\simeq H$}\}.\]
 Then for every $T \geq 1$ and every $\ep>0$,
$
N_f(T;H) \ll_{f,\ep} T^{j-1 + |G/H|^{-1} + \ep} .
$
\end{letterthm}

We also require the following key lemma:
\begin{lemma}\label{L:GtimesG}
Let $f(x,t_1,\dots,t_j)\in \Q(t_1,\dots,t_j)[x]$ be a polynomial with splitting field $K$ over $\Q(t_1,\dots,t_j)$ such that $\Gal(K/\Q(t_1,\dots,t_j))\isom G$.  Suppose that $f(x,t_1,\dots,t_j)$ is regular, i.e. $K$ does not contain a non-trivial finite extension of $\Q$.  
Then $f(x,t_1,\dots,t_j)f(x,s_1,\dots,s_j)$ has splitting field with Galois group $G\times G$ over $\Q(t_1,\dots,t_j,s_1,\dots,s_j)$.
\end{lemma}

\subsubsection{Proof of Lemma \ref{L:GtimesG}}
We will prove the lemma in the case $j=1$; a straightforward extension of this argument applies to the general case.
Let $F(x,t)\in \Q[t,x]$ be a monic irreducible polynomial of $x$ with a root $\theta$ that generates $K$ over $\Q(t)$.   
 We let all our splitting fields be in a fixed  algebraic closure of $\Q(s,t)$.  Then $K\Q(s)$ is the splitting field of $f(x,t)$ over $\Q(s,t)$.  
We will show below that if $G(s,x)\in \Q[s,x]$ is a monic polynomial irreducible over $\Q(s)$ that generates a Galois extension of $\Q(s)$ and does not contain a non-trivial finite extension of $\Q$,
 then $G(s,x)$ is irreducible over $K\Q(s)$.  We will see now that this will suffice to prove the lemma.  Applying this in the case where $F$ is trivial, we will see that $G(s,x)$ is irreducible over $\Q(s,t)$, and in particular, analogously we will see that $F(x,t)$ is irreducible over $\Q(s,t)$ and so $[K\Q(s):\Q(s,t)]=|G|$.  So if $L$ is the splitting field of $f(s,x)$ over $\Q(s)$, then $L$ is generated by $F(s,x)$, and applying the above with $G(s,x)=F(s,x)$, we see that $[KL:K\Q(s)]=|G|$.   Thus $\Gal(KL/\Q(s,t))$ has order $|G|^2$ and injects into $\Gal(K/\Q(t))\times\Gal(L/\Q(s))$, and so $\Gal(KL/\Q(s,t))\isom G\times G.$
Since $KL$ is the splitting field of $f(x,t)f(x,s)$ over $\Q(s,t)$, this proves the lemma.

Now we show that $G(s,x)$ with the assumptions above is irreducible over $K\Q(s)$.  Suppose that $G(s,x)$ factored into $a(x)b(x)$ over $K\Q(s)$.  We can write
$$
a(x)=\sum_{i=0}^{k} \frac{n_i(s,t,\theta)}{d_i(s,t)} x^i
$$
where $n_i(y,z,w)\in \Q[y,z,w]$ and $n_i(y,z)\in \Q[y,z].$  (Since $\theta$ is algebraic over $\Q(s,t)$, we can write elements of $K\Q(s)$ as polynomials in $\theta$ with coefficients in $\Q(s,t)$, and so we can arrange to have no $\theta$'s in the denominators.)
Let $n_{i,j}(z,w)$ be the coefficient of $y^j$ in  $n_i(y,z,w)$.  Let $I_a$ be the ideal in $\Q[z,w]$ generated by all the $n_{i,j}(z,w)$ for all $j$ and for $i\geq 1$, and define $I_b$ analogously.  We claim that, as ideals of $\Q[z,w]$, we have $(F(w,z)) \supset I_a I_b.$  Suppose not. Then there are infinitely many maximal ideals $m$ of $\Q[z,w]$ that contain $(F(w,z))$ but not $I_aI_b$.  Each such maximal ideal gives values $t_0,\theta_0 \in \overline{\Q}$ such that
$F(\theta_0, t_0)=0$, but upon substitution of $z\mapsto t_0$ and $w\mapsto \theta_0$, some element of $I_a$ and some element of $I_b$ remain non-zero, which gives a nontrivial factorization of $G(s,x)$ over $\overline{\Q}(s)$ unless some denominator $d_i(s,t_0)$ is identically zero (or similarly for the denominators in $b(x)$).    Since only finitely many $t_0$ can make a denominator zero, and each have a finitely many associated $\theta_0$, we conclude that $G(s,x)$ factors non-trivially over $\overline{\Q}(s)$, and thus over $E(s,x)$ for some Galois number field $E$.  Since $\Gal(E(s)/\Q(s))\ra\Gal(E/\Q)$ is an isomorphism, the subfields of $E(s)$ that contain $\Q(s)$ are $E'(s)$ for the subfields $E'$ of $E$.  If $M$ is the field generated by $G(s,x)$ over $\Q(s)$, then $[ME(s):E(s)]=[M:M\cap E(s)]$.  Since $G(s,x)$ factors non-trivially over $E(s),$
we have $[ME(s):E(s)]<[M:\Q(s)]$, and thus $M\cap E(s)$ is a non-trivial extension of $\Q(s)$ inside $E(s)$, and thus contains some number field $E'$.  In particular $M$ contains a non-trivial number field, which contradicts our assumption on $G(s,x)$.  Thus, we conclude that $(F(w,z)) \supset I_a I_b,$ and thus $(F(w,z)) | I_a I_b,$ and thus either$(F(w,z)) | I_a$ or $(F(w,z)) |  I_b,$ since $(F(w,z))$ is prime.  But this implies
that either $a$ or $b$ has all coefficients $0$ except the constant one, and thus we conclude $G(s,x)$ is irreducible over $K\Q(s)$. This concludes the proof of Lemma \ref{L:GtimesG}.

\subsubsection{Proof of Theorem \ref{thm_count_G_fields}}
With Lemma \ref{L:GtimesG} and Theorem \ref{thm_HIT} in hand, we may now prove Theorem \ref{thm_count_G_fields}.
Suppose $f(X,T_1,\dots,T_j)$ is a polynomial of total degree $d$ in the $T_i$ with Galois group $G$ over $\Q(T_1,\dots,T_j)$ (with degree $n$ in $X$).  
For $\tbf=(t_1,\dots,t_j)$, we define
 $|\tbf|_\infty = \max_{1 \leq \ell \leq j} |t_\ell|$, so that there are $\gg T^j$ possible values of $\tbf \in \Z^j$ with $|\tbf|_\infty\leq T$.
For each  $\tbf \in \Z^j$, let $L_{\tbf}$ be the splitting field of $f(X,t_1,\ldots, t_j)$ in $\bar{\Q}$.
Let $y$ be the size of the set $\{L_{\tbf} :  \tbf \in \Z^j, |\tbf|_\infty\leq T, \Gal(L_{\tbf}/\Q)\isom G\}$ (note it is possible that different $\tbf$  give the same $L_{\tbf}$), and we also write $L_1,\dots,L_y$ for the fields in this set.

For each $1\leq i \leq y$, suppose $A_i$ of the values $\tbf$ have $L_{\tbf}=L_i$.  So
$$
A_1+\dots +A_y =A,
$$
where $A$ is the total number of values of $|\tbf|_\infty\leq T$ with $\Gal(L_{\tbf}/\Q)\isom G$.
From Theorem \ref{thm_HIT} above, we have that $A \gg_{f} T^j$, since there are finitely many subgroups which each appear with an upper bound with exponent strictly smaller than $j$.

For each  $\tbf \in \Z^{2j}$, let $M_{\tbf}$ be the splitting field of $f(X,t_1,\ldots, t_j)f(X,t_{j+1},\ldots, t_{2j})$.
We ask how many $\tbf\in\Z^{2j}$ with $|\tbf|_\infty\leq T$ have $\Gal(M_{\tbf}/\Q)\isom G$?  By Lemma~\ref{L:GtimesG} and the assumption that $f$ is regular, we have that $f(X,T_1,\dots,T_j)f(X,T_{j+1},\dots, T_{2j})$ has Galois group $G\times G$ over $\Q(T_1,\dots,T_{2j}).$
Thus, by Theorem \ref{thm_HIT}, the number of  $\tbf \in \Z^{2j}$ with $|\tbf|_\infty\leq T$ and $\Gal(M_{\tbf}/\Q)\isom G$ is $\ll_{f,\ep} T^{2j-1+|G|^{-1}+\ep}$.
However, note that this occurs whenever $f(X,t_1,\dots,t_j)$ and $f(X,t_{j+1},\dots, t_{2j})$ have the same splitting field with Galois group $G$, and so
$$
A_1^2+\dots+A_y^2 \ll_{f,\ep} T^{2j-1+|G|^{-1}+\ep}. 
$$
By Cauchy-Schwarz, $(A_1+\cdots + A_y)^2 \leq y (A_1^2+\cdots + A_y^2)$, and we conclude that
$$
y \geq \frac{(A_1+\cdots + A_y)^2}{(A_1^2+\cdots + A_y^2)}  \gg_{f,\ep} \frac{T^{2j}}{ T^{2j-1+|G|^{-1}+\ep}}=T^{1-|G|^{-1}-\ep}.
$$
Thus there are $\gg_{f,\ep}T^{1-|G|^{-1}-\ep}$ different fields  with Galois group $G$ that come from specializations of 
$f(X,T_1,\dots,T_j)$ to some $\tbf$ with $|\tbf|_\infty\leq T$.
For $|\tbf|_\infty\leq T$, we have that 
$f(X,t_1,\dots,t_j)$
 is a degree $n$ polynomial in $X$ with coefficients $\ll_{f} T^d$ and thus with absolute discriminant $\ll_{f} T^{d(2n-2)}$.
Thus $L_{\tbf}$ has absolute discriminant $\ll_{f} T^{d(2n-2)}$.
In conclusion, there are  $\gg_{f,\ep}X^{(1-|G|^{-1}-\ep)/(d(2n-2))}$ degree $n$ $G$-fields with absolute discriminant  at most $X$, completing the proof of Theorem \ref{thm_count_G_fields}.

%%%%%%%%%%%%%%%%%%%%%%%%%%%%%%%%

\addcontentsline{toc}{part}{Part  II: Effective Chebotarev Theorems}

\vspace{.5cm}

\begin{center}
{\bf Part II: Effective Chebotarev Theorems}
\end{center}

\section{Quantitative statements of Chebotarev theorems for families}\label{sec_quant_thms}

We now state quantitative versions of our main Chebotarev theorems, starting with a quantitative version of Theorem \ref{thm_del_ex}.

\begin{thm}[Chebotarev conditional on zero-free region]\label{thm_assume_zerofree1}
Let $k$ be a fixed number field. 
Fix $A \geq 2$, $0 < \del \leq 1/(2A)$, and an integer $n \geq 1$. Let $G$ be a fixed transitive subgroup of $S_{n}$. Then there exists  $D_0 \geq 1$ and $\kappa_1,\kappa_2,\kappa_3>0$
  such that the following holds: for any Galois extension of number fields $L/k$ with $\Gal(L/k) \simeq G$ such that
 $D_L \geq D_0$ and such that
  the Artin $L$-function $\zeta_{L}(s)/\zeta_k(s)$ is zero-free in the region
\beq\label{region1_intro}
 [1-\del,1] \times [-(\log D_L)^{2/\del},(\log D_L)^{2/\del}],
 \eeq
we have that for any conjugacy class $\mathscr{C} \subseteq G$, 
\beq\label{asymp_log_intro}
\left| \pi_{\mathscr{C}}(x,L/k) -  \frac{|\mathscr{C}|}{|G|} \Li(x) \right| \leq \frac{|\mathscr{C}|}{|G|}\frac{x}{(\log x)^A} 
\eeq
for all
\beq\label{x_lower_bd1_intro}
 x \geq \kappa_1 \exp \{ \kappa_2 (\log \log (D_L^{\kappa_3}))^2\} .
\eeq
 If moreover $k=\Q$,  (\ref{asymp_log_intro}) holds for all 
\beq\label{x_lower_bd2_intro}
 x \geq \ka_1 \exp \{ \ka_2 (\log \log (D_L^{\ka_3}))^{5/3}(\log\log \log (D_L^2))^{1/3}\}  .
\eeq
\end{thm}
\begin{remark}
The parameters $D_0$ and $\kappa_1,\kappa_2,\kappa_3$ depend on $n,|G|,A, \del$, and the field $k$; they are precisely  specified in Remark \ref{remark_D0} and  (\ref{nu_dfn_final_ref}), respectively.
\end{remark}

We next state, in quantitative form, the cases of Theorem \ref{thm_gen} that are completely unconditional.

\begin{thm}\label{thm_cheb_main_quant_pt1} 
 For each family $Z_n^\Iscr(\Q,G)$ specified in items (1)--(5) of the list below,
there exist constants $\be,d$ with  $0< \be \leq d$ such that for all $X \geq 1$,
\beq\label{be_d_exist}
 X^\be \ll_{n,G,\Iscr} |Z_n^\Iscr(\Q,G;X)| \ll_{n,G,\Iscr} X^d.
 \eeq

 Moreover, there exists a constant $\tau_*$ with $ 0 \leq \tau_* < \be$, such that for every $\tau> \tau_*$  and every sufficiently small $\ep_0>0$,
 there exists a constant 
$ D_3$ and a constant 
\beq\label{del_final_thm}
\del = \del(\ep_0,m,|G|,d)
\eeq
such that for all $X \geq 1$, there are at most $D_3 X^{\tau + \ep_0}$ $\del$-exceptional fields in $Z_n^{\Iscr}(\Q,G;X)$; here $m$ is the maximum dimension of an irreducible representation of $G$.

Moreover, fix any $A \geq 2$. Then for any $\ep_0>0$ such that $\del$ as defined in (\ref{del_final_thm}) satisfies $\del \leq 1/(2A)$, there exists a constant 
$D_5 \geq 1$ and constants $\kappa_1,\kappa_2,\kappa_3 >0$ 
such that for all $X \geq 1$,  aside from a set $E(X)$ of at most $D_5X^{\tau + \ep_0}$ possible exceptions, each field $K \in Z_n^{\Iscr}(\Q,G;X)$ has the property that for its Galois closure $\tilde{K}$ over $\Q$, for every conjugacy class $\mathscr{C} \subseteq G$, 
\[
\left| \pi_{\mathscr{C}}(x,\tilde{K}/\Q) -  \frac{|\mathscr{C}|}{|G|} \Li(x) \right| \leq \frac{|\mathscr{C}|}{|G|}\frac{x}{(\log x)^A} 
\]
for all
$
 x \geq \ka_1 \exp \{ \ka_2 (\log \log (D_{\tilde{K}}^{\ka_3}))^{5/3}(\log\log \log (D_{\tilde{K}}^2))^{1/3}\} .
$

The families $Z_n^\Iscr(\Q,G)$ are defined by:
\begin{enumerate}
\item $G$ a cyclic group of order $n  \geq 2$, with $\Iscr$ comprised of all generators of $G$ (equivalently, every rational prime that is tamely ramified in $K$ is totally ramified). In this case 
$ |Z_n^\Iscr(\Q,G;X)| \sim c_n X^{1/(n-1)}$ and $\tau_*=0$. (Hence,  the density zero exceptional set $E(X)$ is at most of size $\ll_\ep X^{\ep}$ for every $\ep>0$.)
\item\label{item_S3} $n=3,$ $G \simeq S_3$ acting on a set of 3 elements, $\Iscr$ being the conjugacy class $[(1 \; 2)]$ of transpositions. In this case,  
$ |Z_3^\Iscr(\Q,S_3;X)| \sim c_3 X$ and $\tau_*=1/3$.  (Hence the density zero exceptional set $E(X)$ is at most of size $\ll_\ep X^{1/3+\ep}$ for every $\ep>0$.)
\item\label{item_S4} $n=4,$ $G \simeq  S_4$ acting on a set of 4 elements, $\Iscr$ being the conjugacy class $[(1 \; 2)]$ of transpositions. In this  case, 
$ |Z_4^\Iscr(\Q,S_4;X)| \sim c_4 X$ and $\tau_*=1/2$.  (Hence the density zero exceptional set $E(X)$ is at most of size $\ll_\ep X^{1/2+\ep}$ for every $\ep>0$.)
\item\label{item_Dp} $n=p$ an odd prime, $G = D_p$ the order $2p$ dihedral  group of symmetries of a regular $p$-gon, $\Iscr$ being the conjugacy class of order 2 elements. In this case, for all $X \geq 1$, 
\[ X^{2/(p-1)} \ll_p |Z_p^\Iscr(\Q,D_p;X)| \ll_{p,\ep} X^{3/(p-1) - 1/(p(p-1))+\ep} \]
 and $\tau_*=1/(p-1)$. (Hence   the  density zero exceptional set  $E(X)$ is at most of size $\ll_\ep X^{1/(p-1)+ \ep}$ for every $\ep>0$.)
\item\label{item_A4} $n=4$, $G \simeq A_4$ as a subgroup of $S_4$ acting on a set of 4 elements, $\Iscr$ comprised of the two conjugacy classes of order 3 elements.
 In this case, for all $X \geq 1$, 
 \[ X^{1/2} \ll |Z_4^\Iscr(\Q,A_4;X)| \ll_\ep X^{5/6+\ep},\]
  and $\tau_*=0.2784...$. (Hence  the density zero exceptional set $E(X)$ is at most of size $\ll_\ep X^{0.2784...  \ep}$ for every $\ep>0$.)
\end{enumerate}
\end{thm}

Note we have that $m\leq |G|^{.5}$ (see \cite{VK} for asymptotics when $G=S_n$).

\begin{remark}\label{remark_KM_intro}
 Kowalski and Michel's result \cite[Thm. 2]{KowMic02} leads to the choice
$ \del = \frac{\ep_0}{5m|G|/2 + 2d + 4\ep_0}.$
Within a fixed family $Z_n^\Iscr(\Q,G)$, note that as we choose $\ep_0$ smaller (so that $\del$ correspondingly decreases), the density of potential $\del$-exceptional fields decreases,  in accord with the fact that the requirement that $\zeta_{\tilde{K}}(s)/\zeta(s)$ be zero-free in a box to the right of $\Re(s) = 1-\del$ becomes less stringent, and fewer fields would be expected to violate it. Simultaneously, as $\del$ and accordingly the width of the zero-free region decreases,   the lower-bound threshold for $x$ increases, since the explicit expressions given for  the parameters $\ka_i$ grow with $1/\del$ as specified in (\ref{nu_dfn_final_ref}). This is also as expected.
\end{remark}

\begin{remark}[Cyclic fields of prime degree]\label{remark_cyclic_prime}
If $G$ is a cyclic group of prime order $p \geq 2$, then for each Galois extension $K/\Q$ with Galois group $\simeq G$, every ramified prime is totally ramified, so that for $\Iscr$ as in Theorem \ref{thm_cheb_main_quant_pt1}, $Z_p^\Iscr(\Q,G;X)=Z_p(\Q,G;X)$.
\end{remark}

\begin{remark}[degree $n$ $S_n$-fields with square-free discriminant]\label{remark_degn_Sn}
Recall from \S \ref{sec_sym_gp} that for each $n \geq 2$, the family $Z_n^\Iscr(\Q,S_n;X)$ with $\Iscr= [(1 \; 2)]$ includes all degree $n$ $S_n$-fields with square-free discriminant, which are known in the case of $n=3,4,5$ (and conjectured for $n \geq 6$) to be a positive proportion of all degree $n$ fields.
\end{remark}

\begin{remark}[degree $p$ $D_p$-fields]
It is conjectured that $|Z_p(\Q,D_p;X)| \sim c_{D_p} X^{2/(p-1)}$ for some $c_{D_p}>0$ (see \cite{Mal04}, \cite[p. 608]{Klu06}); assuming this is the true order, our family of degree $p$ $D_p$-fields exhibited in case (\ref{item_Dp}) is a positive proportion of all degree $p$ $D_p$-fields. 
\end{remark}

\begin{remark}[degree $4$ $A_4$-fields]
Based on heuristics as well as numerical evidence, it is conjectured that $|Z_4(\Q,A_4;X)| \sim c_{A_4} X^{1/2} \log X$ for some $c_{A_4}>0$ (see \cite[\S 2.7]{CDO02}, \cite[Ex. 3.2]{Mal04}); assuming this is the true order, our family of degree $4$ $A_4$-fields exhibited in case (\ref{item_A4}) of Theorem  \ref{thm_cheb_main_quant_pt1}  just fails to be a positive proportion of all degree 4 $A_4$-fields.
\end{remark}

Finally, we state the quantitative forms of Theorem \ref{thm_gen} that are conditional on the strong Artin conjecture, and in certain cases on hypotheses for counting number fields.

\begin{thm}[Quintic $S_5$-fields]\label{thm_cheb_main_quant_pt2} 
Consider the family $Z_5^\Iscr(\Q,S_5)$ for $\Iscr$ being the conjugacy class $[(1 \; 2)]$ of transpositions, in which case 
$ |Z_5^\Iscr(\Q,G;X)|  \sim c_5 X$.
 The conclusions of Theorem \ref{thm_cheb_main_quant_pt1} hold for $Z_5^\Iscr(\Q,G)$ if we assume
  the strong Artin conjecture holds for all irreducible Galois representations over $\Q$ with image $S_5$.
In this case, $\tau_* = 199/200$. (Hence the density zero exceptional set $E(X)$ is at most of size $\ll_\ep X^{199/200  +\ep}$ for every $\ep>0.$)
\end{thm}

\begin{remark}\label{remark_S5_Calegari}
An alternative formulation of Theorem \ref{thm_cheb_main_quant_pt2} uses the work of F. Calegari \cite{Cal13}. Let $Y_5^\Iscr(\Q,S_5)$ be the family of quintic $S_5$-fields $K$ such that
complex conjugation in $\Gal(\tilde{K}/\Q)$ has conjugacy class $(1 \; 2)(3 \; 4)$, $\tilde{K}/\Q$ is unramified at 5, and the Frobenius element  at 5 has conjugacy class $(1 \; 2)(3 \; 4)$.
By \cite[Thm. 1.3]{Bha14x}, $  |Y_5^\Iscr(\Q,S_5;X)| \sim c_5' X.$
For these fields, Calegari verifies the strong Artin conjecture for the dimension 4 and 6 irreducible representations of $S_5$, and reduces the verification for the dimension 5 irreducible representations to 
checking that a certain $L$-function is non-vanishing for $s \in [0,1]$. Precisely, for $K \in Y_5(\Q,S_5)$, let $E$ be the quadratic subfield of $\tilde{K}$, $F$ be a subfield of $\tilde{K}$ of degree 6 over $\Q$, and $H$ be the compositum of $E$ and $F$. Then by \cite[Thm. 1.2]{Cal13}, the strong Artin conjecture holds for the dimension 5 irreducible representations as long as $\zeta_H(s)$ is nonvanishing for $s \in [0,1]$.
(See \cite{Boo06}, \cite{Dwy14} for computational verification of this nonvanishing, in a finite number of cases with small discriminant.)
Thus we could alternatively state Theorem \ref{thm_cheb_main_quant_pt2} for the family $Y_5^\Iscr(\Q,S_5)$, assuming in place of the strong Artin conjecture that for each field $K \in Y_5^\Iscr(\Q,S_5)$ considered, the appropriate $L$-function $\zeta_H(s)$ is nonvanishing for $s \in [0,1]$.
\end{remark}

\begin{thm}[degree $n$ $S_n$-fields]\label{thm_cheb_main_quant_pt3} 
Consider for $n \geq 6$ the family $Z_n^\Iscr(\Q,S_n)$ with $\Iscr$ being the conjugacy class $[(1 \; 2)]$ of transpositions, in which case for all $X \geq 1$,
\[X^{1/2+1/n} \ll_n |Z_n^\Iscr(\Q,S_n;X)|   \ll_n X^{ \exp ( C \sqrt{\log n})} . \]
 The conclusions of Theorem \ref{thm_cheb_main_quant_pt1} hold for the family $Z_n^\Iscr(\Q,S_n)$ if we assume 
\begin{enumerate}[label=(\roman*)]
\item  the strong Artin conjecture holds for all irreducible Galois representations over $\Q$ with image $S_n$,
\item\label{Sn_varpi}
  for some $\varpi_n < 1/2+1/n$, for every fixed integer $D$, there are at most $\ll_{n} D^{\varpi_n}$ fields $K \in Z_n(\Q,S_n)$ with $D_K = D$.
\end{enumerate}
In this case, $\tau_*=\varpi_n$. (Hence the   density zero exceptional set $E(X)$ is at most of size $\ll_\ep X^{\varpi_n +\ep}$ for every $\ep>0$.)
\end{thm}
\begin{remark}
If it is known that $|Z_n^\Iscr(\Q,S_n;X)| \gg_n X^{\be_n}$,  then to deduce that the possible exceptional set  has density zero, we need only know \ref{Sn_varpi} for some $\varpi_n < \be_n$. 
\end{remark}

Similarly our results for simple groups are conditional on the strong Artin conjecture; for $A_n$, 
we additionally apply our new lower bound for the number of degree $n$ $A_n$-fields with bounded discriminant.

\begin{thm}[Alternating groups $A_n$, $n \geq 5$]\label{thm_cheb_main_quant_pt2'} 
For each $n \geq 5$, consider the family $Z_n(\Q,A_n)$ (with no restriction on inertia type, that is, $\Iscr = G$). In this case, there exists a positive exponent $\be_n >0$ such that for all $X \geq 1$,  
\[ X^{\be_n} \ll_n |Z_n^\Iscr(\Q,A_n;X)| \ll_n X^{ \exp ( C \sqrt{\log n})}      \]
for a certain absolute constant $C$.
In fact we may take $\be_n=(1-2/n!)/(4n-4).$
Then under the assumption that the strong Artin Conjecture holds for all irreducible Galois representations over $\Q$ with image $A_n$, the conclusions of Theorem \ref{thm_cheb_main_quant_pt1} hold  with $\tau_*=0$. (Hence the density zero exceptional set $E(X)$ is at most of size $\ll_\ep X^{\ep}$ for every $\ep>0$.) \end{thm}

Finally, we state a result for families of fields parametrized by a fixed simple group; here we simply assume that a lower bound that grows like a power of $X$ is known for the number of such fields
(as may be obtained by Theorem \ref{thm_count_G_fields} if an appropriate generating polynomial is known).
\begin{thm}[Simple groups]\label{thm_cheb_main_quant_pt4} 
For $n \geq 2$ and a fixed transitive simple group $G \subset S_n$, the conclusions of Theorem \ref{thm_cheb_main_quant_pt1} hold for the family $Z_n(\Q,G)$ with no restriction on inertia type (that is, $\Iscr = G$), if we assume 
\begin{enumerate}[label=(\roman*)]
\item the strong Artin conjecture for all irreducible representations over $\Q$ with image $G$,
\item a lower bound of the form $|Z_n(\Q,G;X)| \gg_{n,G} X^{\be}$ for some $\be>0$, for all $X \geq 1$.
\end{enumerate}
Then 
$X^{\be} \ll_n |Z_n(\Q,G;X)|   \ll_n X^{ \exp ( C \sqrt{\log n})} $
for an absolute constant $C$, 
and $\tau_*=0$. (Hence the density zero exceptional set $E(X)$ is at most of size $\ll_\ep X^{\ep}$ for  every $\ep>0$.)
\end{thm}
\begin{remark}
At present we do not treat families $Z_n(\Q,G)$ for $G$ a non-cyclic abelian group, or $Z_4(\Q,D_4)$; we remark on difficulties encountered in these settings in Remarks \ref{remark_G_abelian} and \ref{remark_D4}.
\end{remark}

We encapsulate two useful consequences in all the settings described above: 
\begin{cor}[Quantitative counts for small primes]\label{cor_primes_families}
Let $Z_n^\Iscr(\Q,G;X)$ be fixed to be one of the families of fields considered in Theorems \ref{thm_cheb_main_quant_pt1}, \ref{thm_cheb_main_quant_pt2}, \ref{thm_cheb_main_quant_pt3}, \ref{thm_cheb_main_quant_pt2'} and \ref{thm_cheb_main_quant_pt4}, and correspondingly assume the hypotheses (if any) of the relevant theorem. Recall the parameters $\tau_* < \be \leq d$ proved to exist for the family in (\ref{be_d_exist}), and for any sufficiently small $\ep_0>0$, let $\del \leq 1/4$ be defined as in (\ref{del_final_thm}).

(1)  For  any $\sig>0$, there exists a constant $D_6$
 such that  for every $X \geq 1$, every field $K\in Z_n^\Iscr(\Q,G;X)$ that has $D_K \geq D_6$  and is not $\del$-exceptional, has the property that for  any fixed conjugacy class (or finite union of conjugacy classes) $\Cscr$ in $G$,
\beq\label{many_primes}
\pi_\Cscr (D_K^\sig, \tilde{K}/\Q) \gg_{G,n,\sig} \frac{D_K^\sig}{\log D_K}.
\eeq
Here $\tilde{K}$ denotes the Galois closure of $K$ over $\Q$.

(2) For any $\sig >0$, there exists a constant 
$ D_7 $
such that for every $X \geq 1$, every  field $K\in Z_n^\Iscr(\Q,G,X)$ that has $D_K \geq D_7$ and is not $\del$-exceptional, has the property that for any conjugacy class $\Cscr$ of $G$,
\beq\label{pi_diff}
\pi_\Cscr (2D_K^{\sig}, \tilde{K}/\Q) - \pi_\Cscr (D_K^{\sig}, \tilde{K}/\Q)  \geq 1.
 \eeq
 Here $\tilde{K}$ denotes the Galois closure of $K$ over $\Q$.
 
Finally, in either case, recall that for every $\tau > \tau_*$ there exists a constant $D_3$ 
such that for every $X \geq 1$, at most $D_3 X^{\tau +\ep_0}$  fields $K \in Z_n^\Iscr(\Q,G;X)$ are $\del$-exceptional.

\end{cor}

%%%%%%%%%%%%%%%%%%%%%%%%%%%%%%%%%%%%%%%
\section{A Chebotarev density theorem conditional upon a zero-free region}\label{sec_ChebotarevProof}
The main goal of this section is to prove Theorem \ref{thm_assume_zerofree1}. 
A nice feature of Lagarias and Odlyzko's approach to the effective Chebotarev theorem is that it does not assume the Artin conjecture, so that Theorem \ref{thm_LagOdl_uncond} is completely unconditional. Similarly, Theorem \ref{thm_assume_zerofree1} is unconditional, aside from the assumed zero-free region. This is made possible by using Lagarias and Odlyzko's technical trick (originally due to Deuring) of expressing $\zeta_L$ as a product of Hecke $L$-functions, where $L/k$ is a Galois extension of number fields with $\Gal(L/k) \simeq G$.  Fixing an element $g \in G$ and letting $H = \langle g \rangle$ be the cyclic subgroup of $G$ generated by $g$, then upon setting $E$ to be the fixed field $L^H$, Lagarias and Odlyzko obtain the product expression on the left, in which $\chi$ varies over the irreducible characters of $H$: 
\beq\label{zeta_two_prod}
 \prod_{\chi \; \mathrm{irred}} L(s,\chi,L/E) = \zeta_L(s) =   \prod_{\rho_j} L(s,\rho_j,L/k)^{\dim \rho_j} .
\eeq
Each such factor is a Hecke $L$-function and hence is known to be entire if $\chi$ is nontrivial.
On the other hand, once we have Theorem \ref{thm_assume_zerofree1}, to deduce the assumed zero-free region via Kowalski-Michel, we will also factor $\zeta_L(s)$ as on the right-hand side, as a product of Artin $L$-functions, which we then need to show (or assume) are automorphic $L$-functions with certain properties.

If one is willing to assume the Artin conjecture, so that each factor on the right-hand side of (\ref{zeta_two_prod}) with $\rho_j$ nontrivial is entire, a Chebotarev density theorem with an effective error term is relatively quick to prove, since either a standard zero-free region (or the GRH zero-free region) may be applied to each of these Artin L-functions, obviating the alternative Hecke factorization; see for example, \cite[\S5.13 and Thm. 4.13]{IK}. (The conjugacy class $\Cscr$ of interest is picked out via trace functions, much as in Dirichlet's theorem on primes $p \con a \modd{q}$,  the residue class of interest is picked out via Dirichlet characters.) In our  application to families of fields we do indeed assume the strong Artin conjecture (or it is known). Nevertheless, to prove Theorem  \ref{thm_assume_zerofree1} we have used the Lagarias-Odlyzko approach, as we expect its unconditionality to be useful for other applications.

\subsection{Standard lemmas on zeroes}

We recall the currently best known zero free region for $\zeta(s)$, due to Vinogradov \cite{Vin58} and Korobov \cite{Kor}.

\begin{lemma}[Vinogradov-Korobov zero-free region for $\zeta(s)$]\label{lemma_zerofree_VK}
There exists an absolute constant $c_\Q>0$ such that $\zeta(s)$ has no zero $s = \sigma + it$ in the region
\beq\label{lemma_zerofree_VK_eqn}
\sigma \ge 1-\frac{c_\Q}{(\log(|t|+2))^{2/3}(\log\log(|t|+3))^{1/3}}.
\eeq
\end{lemma}

We will also use a standard zero-free region for any Dedekind zeta function \cite[Theorem 5.33]{IK}.

\begin{lemma}[ Standard zero-free region for $\zeta_k(s)$]\label{lemma_zerofree_std}
Let $k/\Q$ be a number field of degree $n_k \geq 1$ and with absolute discriminant $D_k$. There exists an absolute constant $c_k>0$ such that $\zeta_k(s)$ has no zero $s = \sigma + it$ in the region

\beq\label{lemma_zerofree_std_eqn}
\sigma \ge 1-\frac{c_k}{n_k^2\log (D_k(|t|+3)^{n_k})},
\eeq
except possibly a simple real ``exceptional'' zero $\beta_0^{(k)} <1$. 
\end{lemma}

We also recall a standard count for zeroes of Dedekind zeta functions at a fixed height:
 
\begin{lemma}[{\cite[Theorem 5.31, Proposition 5.7]{IK}}]\label{lemma_zetazeroes}
Let $k/\Q$ be a number field of degree $n_k \geq 1$ and with absolute discriminant $D_k$.
For a real variable $t$, let $n_k(t)$ denote the number of zeroes $\rho=\beta+i\gamma$ of $\zeta_k(s)$ with $0<\beta<1$ and $|\gamma-t|\le 1$. For all real $t$,
$
n_k(t) \ll \log D_k + n_k\log(|t|+4).
$
\end{lemma}

The corresponding result for Hecke $L$-functions is:
\begin{lemma}[{\cite[Lemma 5.4]{LagOdl75}}]\label{lemma_zetazeroes_Hecke}
Let $n_\chi(t)$ denote the number of zeroes $\rho=\beta+i\gamma$ of a Hecke $L$-function $L(s,\chi,L/E)$ with $0<\beta<1$ and $|\gamma-t|\le 1$. Let $F(\chi)$ denote the conductor of $\chi$, and $A(\chi)=D_E \Nm_{E/\Q}(F(\chi))$.  For all $t$,
$
n_\chi(t) \ll \log A(\chi) +n_E\log(|t|+2).
$
\end{lemma}

\subsection{Explicit description of assumed zero-free region} 
We now prove Theorem \ref{thm_assume_zerofree1}, using in particular the assumption, in the theorem statement, that $\zeta_L(s)/\zeta_k(s)$ is zero-free in the region 
\beq\label{region1'_0}
 [1-\del,1] \times [-(\log D_L)^{2/\del},(\log D_L)^{2/\del}].
 \eeq
 For use in potential computational applications, we specify the dependencies of all parameters on $k, \del$, etc., although we do not now optimize them (e.g. compared to recent work conditional on GRH in \cite{GreMol17}), as it is not relevant for our current applications. 
 
We may assume that $L$ has degree $n_L >1$ over $\Q$, since in the case $L=k = \Q$, $\pi_{\mathscr{C}}(x,L/k)$ is simply counting rational primes $p \leq x$.
Our proof will proceed in two stages: first, we  deduce from Theorem \ref{thm_LagOdl_uncond} of Lagarias and Odlyzko that the conclusion of Theorem \ref{thm_assume_zerofree1} is true if $x$ is sufficiently large. Second, for small $x$, we refine the method of Lagarias and Odlyzko, keeping track of the assumed zero-free region. (This manner of partitioning into large and small $x$ has appeared in the proof of the prime ideal theorem of \cite[Theorem 2.6]{ChoKim14a}.)

At each step, when we state that something holds for a number field $k$, it also applies to $k=\Q$; separately, we give refined statements so far applicable only to $k=\Q$.
 We do not rule out \emph{a priori} the possibility of an exceptional zero of $\zeta_L(s)$, say $\be_0$. Instead, in our application of Theorem \ref{thm_LagOdl_uncond}, the main idea is to assume that $D_L$ is sufficiently large that the real interval within the region (\ref{LO_region}) in Theorem \ref{thm_LagOdl_uncond} is contained inside the assumed zero-free region (\ref{region1'_0}), and thus $\zeta_L$ cannot have an exceptional zero $\be_0$.
In order to carry this out rigorously, we must be more careful, since (\ref{region1'_0}) is an assumed zero-free region for $\zeta_L/\zeta_k$ and not just $\zeta_L$. 

The function $\zeta_k(s)$ may have an exceptional (real) zero in the standard region \eqref{lemma_zerofree_std_eqn} given in Lemma \ref{lemma_zerofree_std}; we will denote this, if it exists, $\be_0^{(k)}$. (Of course when $k=\Q$, $\zeta_k(s) = \zeta(s)$, and no such exceptional zero exists.) 
Since $k$ is fixed, $\be_0^{(k)}$ is fixed. 
We now fix a new parameter $\del_0$ so that 
\beq\label{del0_dfn}
1 - \del_0 \geq 1- \del, \qquad \text{and} \qquad 1- \del_0 >\be_0^{(k)};
\eeq
we set $\del_0 = \del$ if $k=\Q$.  (Throughout this section we will use the notation $\del_0$; in the statement of Theorem \ref{thm_assume_zerofree1}, any dependence on $\del_0$ is equivalently a dependence on $\be_0^{(k)}$ and $\del$.)
From now on, instead of the zero-free region (\ref{region1'_0}), we work with the possibly smaller region 
\beq\label{region1'}
[1-\del_0,1] \times [-(\log D_L)^{2/\del},(\log D_L)^{2/\del}],
\eeq
which excludes the possible fixed zero $\be_0^{(k)}$. 

By our hypothesis, the Artin $L$-function $\zeta_L(s)/\zeta_k(s)$ has no zeroes in the region (\ref{region1'}), and it is an entire function by the Aramata-Brauer theorem; $\zeta_k(s)$ has no zeroes in the intersection of regions (\ref{lemma_zerofree_std_eqn}) and (\ref{region1'}) (respectively, no zeroes in the intersection of the regions (\ref{lemma_zerofree_VK_eqn}) and (\ref{region1'}) if $k=\Q$) and is holomorphic there. Thus $\zeta_L(s)$ has no zeroes  in the intersection of (\ref{lemma_zerofree_std_eqn}) and (\ref{region1'}) (respectively, (\ref{lemma_zerofree_VK_eqn}) and (\ref{region1'}) if $k=\Q$). 

\begin{figure}[h] %  figure placement: here, top, bottom, or page
   \centering
   \includegraphics[width=3.5in]{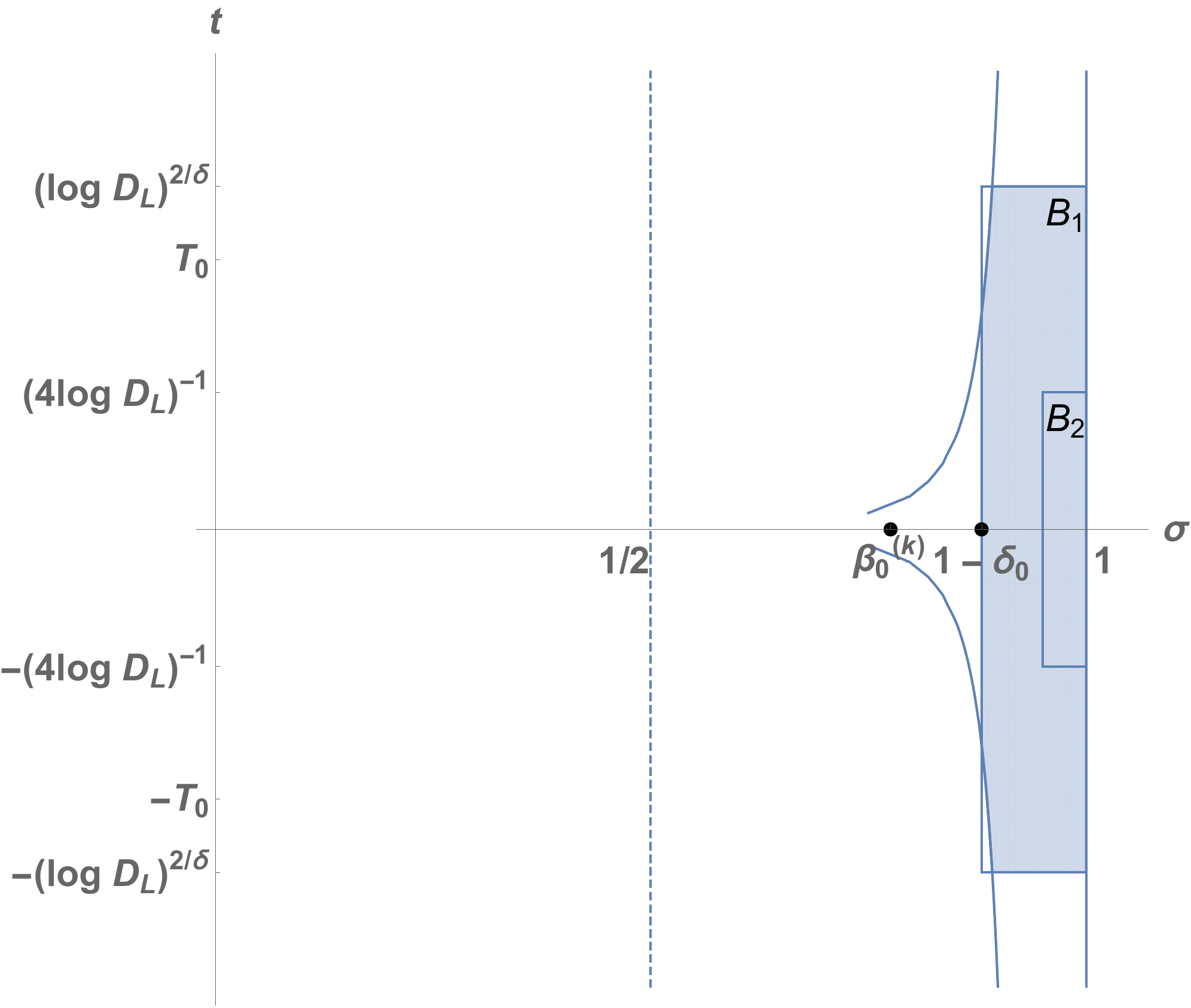} 
   \caption{The curved region represents the standard zero-free region for $\zeta_k$ and the point $\be_0^{(k)}$ denotes the possible (real) exceptional zero of $\zeta_k$. The larger box $B_1$ is the assumed zero-free region (\ref{region1'}) for $\zeta_L/\zeta_k$. The shaded region represents the consequent (assumed) zero-free region for $\zeta_L$,  determined by   the intersection of the known standard zero-free region for $\zeta_k$ and $B_1$. The box $B_2$ is the zero-free region (\ref{LO_region}) known to hold for $\zeta_L$, aside from a possible exceptional (real) zero; we will conclude no such zero can exist in $B_2$ as long as $D_L$ is sufficiently large. }\label{figure_zetaL}
   \label{fig:region}
\end{figure}
Thus we now specify (under the above hypotheses) the zero-free region of $\zeta_L(s)$ (see Figure \ref{figure_zetaL}): 
\begin{equation}\label{region}
\begin{cases}
\sigma \ge 1 - \delta_0 & {\text{if } |t| \le T_0},\\[.15in]
\displaystyle\sigma \ge 1 - \mathscr{L}(t) & {\text{if } T_0\le |t| \le (\log D_L)^{2/\del}},
\end{cases}
\end{equation}
where
\beq\label{L_t_dfn}
\mathscr{L}(t) = \begin{cases} \frac{c_k}{n_k^2 \log (D_k(|t| + 3)^{n_k})} & \text{general $k$} \\
						\frac{c_\Q}{(\log (|t|+2))^{2/3}(\log\log(|t|+3))^{1/3}} & \text{if $k  = \Q$,}
						\end{cases}
\eeq
and  $T_0$ is the height at which the zero-free region (\ref{lemma_zerofree_std_eqn})  for $\zeta_k$ (respectively (\ref{lemma_zerofree_VK_eqn}) for $\zeta$) intersects the line $\Re(s)=1-\del_0$. In our Chebotarev theorems we are interested in the range where $D_L \maps \infty$, so there is no harm in always assuming (for simplicity) that $D_L$ is sufficiently large that the left-hand boundary $\Re(s) = 1-\del_0$ of (\ref{region1'}) intersects the boundary of (\ref{lemma_zerofree_std_eqn}) (respectively (\ref{lemma_zerofree_VK_eqn}) if $k=\Q$) at a height $T_0 \leq (\log D_L)^{2/\del}$.  For example, for a field $k$ and the zero-free region (\ref{lemma_zerofree_std_eqn}), we compute that
\beq\label{T0}
 T_0 = D_k^{-1/n_k}\exp\left(\frac{c_k}{\delta_0 n_k^3}\right)-3 .  
 \eeq
A similar computation may be done to find $T_0$ in the case $k=\Q$ with the improved zero-free region (\ref{lemma_zerofree_VK_eqn}).
In either case, to have  $T_0 \leq (\log D_L)^{2/\del}$ it is sufficient to have 
 \beq\label{DL_lower_from_zerofree}
  D_L \geq \begin{cases} \exp\{\exp (c_k \del/\del_0)\} & \text{general $k$} \\
 	  \exp \{ (\exp \exp (c_\Q/\del))^{2/\del}\} & \text{if $k=\Q$};
	\end{cases}
	\eeq
we refer to this lower bound as $D_0'=D_0'(c_k,\del_0,\del)$.

\subsection{The proof of Theorem \ref{thm_assume_zerofree1} for large $x$}\label{sec_large_x}

With this zero-free region in mind, we dispatch the case of our Chebotarev theorem for large $x$, that is, for $x \geq \exp (10n_L (\log D_L)^2)$.
Recall the standard zero-free region (\ref{LO_region}) which is known to hold for $\zeta_L(s)$, aside from a possible real exceptional zero. 
We may define a constant $D_1=D_1(\del_0)$ so that 
\beq\label{delDdel}
 1 - \del_0 < 1 - (4 \log D_1(\del_0))^{-1}. 
 \eeq
For later purposes, we also assume $D_1(\del_0) \geq 4$. 
Our conclusion now is that for $D_L \geq D_1(\del_0)$, $\zeta_L$ can have no (real, exceptional) zero in the region  (\ref{LO_region}), and thus under the hypotheses of Theorem \ref{thm_assume_zerofree1}, the result of Theorem \ref{thm_LagOdl_uncond} holds without the $\be_0$ term. 

Now in order to show the remaining error term in Theorem \ref{thm_LagOdl_uncond} (with absolute constants $C_1,C_2$) is sufficiently small, as claimed in Theorem \ref{thm_assume_zerofree1}, we need only verify that there exists a constant $D_1'=D_1' (C_1,C_2,n_L,A)$ such that as long as $D_L \geq D_1'$, for all $x \geq \exp(10n_L(\log D_L)^2)$, 
\beq\label{bigx}
C_1 x \exp(-C_2 n_L^{-1/2} (\log x)^{1/2} ) \leq \frac{|\mathscr{C}|}{|G|} x (\log x)^{-A} .
 \eeq
 In fact it suffices that $D_L$ is sufficiently large that (\ref{bigx}) holds at the endpoint $x = \exp (10 n_L(\log D_L)^2)$, which is equivalent to requiring $D_L \geq c_2 (\log D_L)^{c_3}$ with $c_2 = (c_1^{-1/(2A)} (10 n_L)^{1/2})^{(2A)C_2^{-1} 10^{-1/2}}$ and $c_3 = 2AC_2^{-1} 10^{-1/2}$; this provides the necessary threshold $D_1'$.
As a consequence, the conclusion of Theorem \ref{thm_assume_zerofree1} holds for $x \geq \exp (10n_L (\log D_L)^2)$, as long as $D_L \geq   \max \{D_1,D_1'\}$.

\subsection{Small $x$}
In the remaining region of small $x$ (that is, for $x< \exp (10n_L (\log D_L)^2)$) we return to the original strategy of Lagarias and Odlyzko, which will be our focus for the remainder of \S \ref{sec_ChebotarevProof}.
As in the classical prime number theorem, it is convenient to work originally with a weighted prime-counting function, defined in this case by
\[ \psi_{\mathscr{C}} (x,L/k) = \sideset{}{'}\sum_{\tstack{\pfrak, m}{\Nm_{k/\Q} \pfrak^m \leq x}{\left[ \frac{L/k}{\pfrak} \right]^m = \mathscr{C}}} \log( \Nm_{k/\Q} \pfrak);\]
the final result for $\pi_{\mathscr{C}}(x,L/k)$ will then follow from partial summation.
Here $\Sigma'$ denotes  that the sum is restricted to those prime ideals $\pfrak$ in $\RI_{k}$ that are unramified in $\RI_{L}$. 
The notation $\left[ \frac{L/k}{\pfrak} \right]^m = \mathscr{C}$ denotes the requirement that if we pick any prime ideal $\qfrak \subset \Ocal_L$ lying above $\pfrak$, then $\mathscr{C}$ is the conjugacy class of the $m$-th power $(\sig_\qfrak)^m$ of the Frobenius element $\sig_\qfrak$ inside $G$. (This is well-defined no matter which prime $\qfrak$ is chosen above $\pfrak$, since if $\qfrak'= \tau(\qfrak)$ for some nontrivial automorphism $\tau \in G$, then $(\sig_{\qfrak'})^m = (\tau \sig_{\qfrak} \tau^{-1})^m = \tau (\sig_\qfrak)^m \tau^{-1}$, so that they lie in the same conjugacy class in $G$.)

Our main result for $\psi_{\mathscr{C}}$ in the region of small $x$ is as follows:
\begin{prop}\label{prop_assume_zerofree2}
Let $k$ be a fixed number field. 
Fix $A \geq 2$, $0 < \del \leq 1/(2A)$, and an integer $n \geq 1$. Let $G$ be a fixed transitive subgroup of $S_{n}$. Then
for any absolute constant $0<c_0 \leq 1$ of our choice, there exists  a constant $D_2$ and constants $\kappa_1',\kappa_2',\kappa_3'$ such that for any Galois extension of number fields $L/k$ with $\Gal(L/k) \simeq G$ such that
$D_L \geq \max \{ D_0',D_1, D_2 \},$
 and such that
  the Artin $L$-function $\zeta_{L}(s)/\zeta_k(s)$ is zero-free in the region
\beq\label{region1_intro2}
 [1-\del,1] \times [-(\log D_L)^{2/\del},(\log D_L)^{2/\del}],
 \eeq
we have for every conjugacy class $\Cscr$ in $G$ that
\[ \left| \psi_{\mathscr{C}}(x,L/k) - \frac{|\mathscr{C}|}{|G|} x  \right| \leq c_0 \frac{|\mathscr{C}|}{|G|} \frac{x}{(\log x)^{A-1}}, \]
as long as
\beq\label{x_range_k}
\kappa_1^{\prime} \exp \{ \kappa_2^{\prime} (\log \log (D_L^{\kappa_3^{\prime}}))^2\}  \leq x \leq \exp \{10 n_L (\log D_L)^2\}.
\eeq
If moreover $k  = \Q$ we may take $x$ in the range 
 \beq\label{x_range_Q}
   \ka_1^{\prime} \exp \Big\{ \ka_2^{\prime} (\log \log D_L^{\ka_3^{\prime}})^{5/3}(\log\log\log (D_L^2) )^{1/3}\Big\} \leq x \leq \exp \{10 n_L (\log D_L)^2\}.
   \eeq
\end{prop}

\begin{remark}
Recall that $D_0'$ was fixed by (\ref{DL_lower_from_zerofree}), $D_1$ was fixed by (\ref{delDdel}); we will construct $D_2$  explicitly in Lemma \ref{lemma_E1E2_bounds}. The constants $\kappa_1',\kappa_2',\kappa_3'$ depend on $c_0, D_k,n_k, n_L, \del_0,\del, A$ and are chosen in  (\ref{kappa_dfn_prime}).
\end{remark}

\subsection{The passage to sums over zeroes of Hecke $L$-functions}
To prove this proposition, we rebuild the argument of Lagarias and Odlyzko, inserting the zero-free region (\ref{region}) at a key point. With $\mathscr{C}$ the fixed conjugacy class of interest, we fix any element $g \in \mathscr{C}$ and let $H = \langle g \rangle$ be the cyclic group generated by $g$. Then $H$ defines a fixed field $E = L^H$ with $k \subseteq E \subseteq L$, and the cyclic group $H$ has an associated family of irreducible one-dimensional characters. For any such character $\chi$, we consider the Hecke $L$-function $L(s,\chi,L/E)$; in particular
if $\chi=\chi_0$ is the trivial character on $H$ then $L(s,\chi,L/E) = \zeta_E(s)$.
The following statement provides the key framework for proving Proposition \ref{prop_assume_zerofree2}:
\begin{prop}[Theorem 7.1 of \cite{LagOdl75}]\label{prop_main_errors}
For $L/k$ a finite Galois extension of number fields with $ \Gal(L/k) \simeq G$, cyclic subgroup $H \subseteq G$, and $k \subseteq E \subseteq L$ as described above,
there exists an absolute  constant $C_5 \geq 1$ such that if $x \geq2$ and $T \geq 2$, then
\beq\label{eqn_C5}
 \left|\psi_{\mathscr{C}}(x,L/k) - \frac{|\mathscr{C}|}{|G|} x \right| \leq C_5\frac{|\mathscr{C}|}{|G|}
\left\{ S(x,T) + E_1 + E_2 \right\},
\eeq
in which 
\[ S(x,T) =  \sum_{\chi} \overline{\chi}(g) \left( \sum_{\bstack{\rho = \be + i\ga}{|\ga|<T}} \frac{x^\rho}{\rho} - \sum_{\bstack{\rho = \be + i \ga}{|\rho|<1/2}} \frac{1}{\rho}\right),\]
where the sum is over irreducible characters $\chi$ of $H$, and for each character $\chi$ the inner sums are over nontrivial zeroes $\rho = \be + i \ga$ of the Hecke $L$-function $L(s,\chi,L/E)$, and
\begin{align}
 E_1 &=  x T^{-1} \log x \log D_L  + \log D_L + n_L \log x + n_L xT^{-1} \log x \log T, \label{E1_dfn} \\
 E_2 & =  \log x \log D_L + n_L x T^{-1} (\log x)^2. \label{E2_dfn}
 \end{align}
\end{prop}
\begin{remark}
Note that we may assume that $C_5 \geq 1$, by enlarging it if necessary.
As stated in (\ref{E2_dfn}), $E_2$ is slightly refined over \cite[Theorem 7.1]{LagOdl75}, which in place of $|\mathscr{C}||G|^{-1}E_2$ has
\[ E_2' =  \log x \log D_L + n_k x T^{-1} (\log x)^2,\]
(without a factor of $|\mathscr{C}||G|^{-1}$). As noted in \cite[Th\'{e}or\`{e}me 4]{Ser82}, the first term in $E_2'$ may be replaced by 
\[|G|^{-1} \log x \log D_L \leq |\mathscr{C}||G|^{-1} \log x \log D_L,\]
 by a refined estimate for a sum over prime ideals $\pfrak \subset \Ocal_k$ that ramify in $L$. For the second term in $E_2'$, we use the trivial observation that $n_k |G| = n_L$, so that 
\[ n_k x T^{-1} (\log x)^2 = |G|^{-1} n_L x T^{-1} (\log x)^2 \leq |\mathscr{C}||G|^{-1} n_L x T^{-1} (\log x)^2,\]
 as claimed.
\end{remark}

With Proposition \ref{prop_main_errors} in hand, Lagarias and Odlyzko use zero-free regions (either unconditional or on GRH) to deduce a bound for $S(x,T)$, which indicates an appropriate choice for the height $T$ that guarantees all the error terms are sufficiently small.
We proceed with a different zero-free region and  a different choice for $T$, namely 
\beq\label{T_choice}
T  = (\log D_L)^{2/\del},
\eeq
 where $\del$ is provided from our assumed zero-free region (\ref{region}). (In particular, we may assume  that $T \geq 2$ as long as $D_L \geq 3 > \exp(2^{\del/2})$, upon recalling $\del \leq 1/4$.)  

\subsection{Bounding the contribution of zeroes $|\rho| < 1/2$  in $S(x,T)$}

The contribution to $S(x,T)$ from $|\rho| < 1/2$ (so that certainly $|\ga| \leq T$ with $T$ as in (\ref{T_choice})) is bounded  by:
\beq\label{S_term1_bound}
 \sum_{\chi}  \sum_{\bstack{|\rho| < 1/2}{|\ga| \leq T}} \left\{\left| \frac{x^\rho}{\rho} \right|+\left| \frac{1}{\rho} \right|\right\} 
 \ll x^{1/2} \sum_{\chi}  \sum_{|\rho| < 1/2} \left| \frac{1}{\rho} \right|
  \ll  x^{1/2}n_L  (\log D_L)^2,
\eeq
 in which the implied constant is absolute.
The first inequality is clear; to prove the second inequality,
 recall the factorization (\ref{zeta_two_prod}) into Hecke $L$-functions,
\beq\label{zeta_prod2}
\zeta_L(s) = \zeta_E(s) \prod_{\chi \neq \chi_0} L(s,\chi,L/E),
  \eeq
  with the product over non-trivial irreducible characters of $H$. The Hecke $L$-functions are entire, and  $\zeta_E(s)$ and $\zeta_L(s)$ each have their only pole at $s=1$; thus (rigorously by multiplying both sides of the identity by $(s-1)$), it follows that none of the factors on the right hand side of  (\ref{zeta_prod2}) have a zero  in the region \eqref{region}. Recalling that (\ref{region})  contains the region (\ref{LO_region}) since $D_L \geq D_1(\del_0)$   we may conclude (by the functional equation) that each $L(s,\chi,L/E)$ is zero-free both in  (\ref{LO_region}) and in $0 \leq \sig \leq (4\log D_L)^{-1}, |t| \leq (4 \log D_L)^{-1}$. 
  Thus the only zeroes that can appear in (\ref{S_term1_bound}) must have $|\rho| \geq (4\log D_L)^{-1}$; recalling the notation of Lemmas \ref{lemma_zetazeroes} and \ref{lemma_zetazeroes_Hecke}, we then see that for each $\chi$, 
\beq\label{sum_zeroes_small}
\sum_{|\rho| < 1/2}\left| \frac{1}{\rho} \right| \le 4 (\log D_L)n_\chi(1) \ll (\log D_L)(\log A(\chi)+n_E\log 3) 
\eeq
with the implied constant being absolute.
The conductor-discriminant formula \cite[Ch. VII 11.9]{Neu}
 shows 
\beq\label{sum_A}
\sum_{\chi}\log A(\chi) 
 = \log\left[ D_E^{|H|}\Nm_{E/\Q} \left(\prod_{\chi}F(\chi)\right)\right]
 = \log \left[D_E^{[L:E]}\Nm_{E/\Q}(D_{L/E})\right]
=\log D_{L}.
\eeq
Thus, summing  (\ref{sum_zeroes_small}) over $\chi$ we have
\begin{align*}
\sum_{\chi}\sum_{|\rho| < 1/2}\left| \frac{1}{\rho} \right| 
&\ll   (\log D_L)^2 + n_E|H|\log 3 \ll n_L (\log D_L)^2,
\end{align*}
with an absolute implied constant, verifying (\ref{S_term1_bound}).

\subsection{Bounding the contribution of $|\ga| \leq T$ in $S(x,T)$}
Suppose that $\rho = \be + i\ga$ is a nontrivial zero of $L(s,\chi,L/E)$ with $|\ga| \leq T$ and $|\rho| >1/2$.  Recalling the definition (\ref{T0}) of the height $T_0$, by the assumption of the zero-free region (\ref{region}), we know that \emph{without exception,} all zeroes $\rho$ with $|\ga| \leq T_0$ have $\be \leq 1 - \del_0$, so that 
$
 |x^\rho| = x^\be \leq x^{1-\del_0}.
$
 Similarly, all zeroes $\rho$ with $T_0 \le |\ga| \leq T$ have $\be \leq 1-\mathscr{L}(T)$, so that 
$
 |x^\rho| = x^\be \leq x^{1-\mathscr{L}(T)}.
$
We also note that for any fixed $\chi$, by Lemma \ref{lemma_zetazeroes_Hecke},  
\begin{align*}
 \sum_{|\ga| \leq T_0}\left| \frac{x^\rho}{\rho} \right| 
& \ll x^{1-\del_0}	\sum_{j\le T_0}\frac{n_\chi(j)}{j} \nonumber \\
& \ll  x^{1-\del_0} (\log T_0)(\log A(\chi)+n_E\log T_0) \nonumber \\
&\ll  x^{1-\del_0} (\log T)(\log A(\chi)+n_E\log T); 
\end{align*}
similarly, \begin{align*} 
 \sum_{T_0 \leq|\ga| \leq T}\left| \frac{x^\rho}{\rho} \right| 
\ll  x^{1-\Lscr(T)} (\log T)(\log A(\chi)+n_E\log T).
\end{align*}
Summing over all $\chi$ as in (\ref{sum_A}), we see that 
\[x^{1-\del_0} \sum_\chi (\log T)(\log A(\chi) +n_E\log T)  \ll x^{1-\del_0} (\log T) \{ \log D_L + n_L\log T\},\]
and, likewise, 
\[x^{1-\Lscr(T)} \sum_\chi (\log T)(\log A(\chi) +n_E\log T)  \ll x^{1-\Lscr(T)} (\log T) \{ \log D_L + n_L\log T\}.\]
Combining these estimates with (\ref{S_term1_bound}), we may conclude
\beq\label{Sxt_terms}
|S(x,T)| \leq C_6 \left\{ E_3 + E_4 + E_5  \right\},
\eeq
for an absolute constant $C_6$ (which we may assume satisfies $C_6 \geq 1$), and 
\begin{align*}
E_3 &=x^{1/2}n_L  (\log D_L)^2 \\
E_4 & =  x^{1 - \del_0} (\log T )\log (D_L T^{n_L})\\
E_5 & =  x^{1-\Lscr(T)}(\log T )\log (D_L T^{n_L}).
\end{align*}

The proof of Proposition \ref{prop_assume_zerofree2} will then be complete, upon verification of two lemmas, which we record as Lemma \ref{lemma_Sxt_bounds} and Lemma \ref{lemma_E1E2_bounds} below.

\begin{lemma}\label{lemma_Sxt_bounds}
 Let $k$ be a fixed number field. Let $A \geq 2$ be fixed and let $0 < \del \leq 1/(2A)$ be a fixed positive constant; define $\del_0$ from $\del$ as in (\ref{del0_dfn}) according to whether or not $\zeta_k(s)$ has an exceptional zero. 
 Let $L/k$ be a Galois extension of number fields with $\Gal(L/k) \simeq G$, and assume that
  the Artin $L$-function $\zeta_{L}(s)/\zeta_k(s)$ is zero-free in the region
\beq\label{region1_intro3}
 [1-\del,1] \times [-(\log D_L)^{2/\del},(\log D_L)^{2/\del}].
 \eeq
For $D_L \geq D_1(\del_0)$ (as defined in (\ref{delDdel})), for any choice of absolute constant $0<c_1 \leq 1$, we have
\beq\label{SxT_bound}
 |S(x,T)| \leq 3c_1C_6  x (\log x)^{-(A-1)}
 \eeq
for all 
 \beq\label{x_lower_bound}
\kappa_1'' \exp \{ \kappa_2'' (\log \log (D_L^{\kappa_3''}))^2\} \leq x \leq \exp \{10 n_L (\log D_L)^2\},
 \eeq
 where 
 \begin{align}\label{kappa_dfn_primeprime}
 \kappa_1'' &:= (6c_1^{-1}10^{A-1} n_L^{A})^{1/\del_0} \del^{-2/\del_0}  \\
  \kappa_2'' & := \max \{ 2A \del_0^{-1}, 4A c_k^{-1} n_k^3 \del^{-1}\}  \nonumber \\
 \kappa_3'' & :=6 c_1^{-1/(2A)} D_k n_L \del^{-1/A}.  \nonumber
 \end{align}
Moreover, if $k  = \Q$ we may consider all
 \beq\label{x_lower_bound_Q}
  \ka_1'' \exp \Big\{ \ka_2'' (\log \log (D_L^{\ka_3''}) )^{5/3} (\log\log \log ( D_L^2) )^{1/3}\Big\} \leq x \leq \exp \{10 n_L (\log D_L)^2\}.
  \eeq
\end{lemma} 

It is in Lemma \ref{lemma_Sxt_bounds} that we fully utilize the fact that the zero-free region (\ref{region1'}) has a width that is independent of $D_L$; this is key to obtaining a small lower threshold on $x$.
\begin{proof}
The lemma is proved by simple computations. For any field $k$, we see that in the range $x \leq \exp (10n_L(\log D_L)^2)$, to guarantee $|E_3| \leq c_1 x (\log x)^{-(A-1)}$ it suffices that
\[ x  \geq c_1^{-2}(10)^{2(A-1)} n_L^{2A} (\log D_L)^{4A};\]
here we have explicitly used the upper bound $x \leq \exp (10n_L(\log D_L)^2)$.
Similarly for such an upper bound for $E_4$ it suffices to have
\[ x  \geq (6c_1^{-1} (10)^{A-1} n_L^{A})^{1/\del_0} \del^{-2/\del_0}(\log D_L)^{2A/\del_0}, \]
provided that $T  = (\log D_L)^{2/\del}$ and $x \leq \exp (10n_L(\log D_L)^2)$. Since $\del_0 \leq \del \leq 1/4$, both of the lower bounds for $x$ displayed above are satisfied if 
\beq\label{E3E4_both}
x \geq (6 c_1^{-1} 10^{A-1} n_L^{A})^{1/\del_0} \del^{-2/\del_0} \exp \{ 2A \del_0^{-1} \log \log D_L\}  .
\eeq
The distinction  of $k=\Q$ only appears in the treatment of  $E_5$; in the case of $k=\Q$, it suffices to find a lower bound on $x$ such that
\[x^{1-\frac{c_\Q}{(\log(T+2))^{2/3}(\log\log(T+3))^{1/3}}}(\log T)\log(D_LT^{n_L})  \leq c_1 x (\log x)^{-(A-1)}, \]
as long as $x \leq  \exp (10n_L(\log D_L)^2)$, $T  = (\log D_L)^{2/\del}$ and $D_L \geq D_1(\del_0)$.
Here one sees that it would suffice to have 
\[
x \geq \exp \{2A c_\Q^{-1}(\log(2(\log D_L)^{2/\delta}))^{2/3}(\log\log(2(\log D_L)^{2/\delta}))^{1/3} \log [ c_2^{1/(2A)} \del^{-1/A} (\log D_L)]\} ,
\]
with $c_2 = 6c_1^{-1}(10)^{A-1}n_L^A $, which can be simplified as the requirement that $x$ is at least
 \beq\label{SxT_E5aQ_main}
   \exp\left\{ 4A c_\Q^{-1}\del^{-2/3}(\log (2\del^{-1})+1)^{1/3} \left(\log\log D_L^{\max\{2^{\delta/2},c_2^{1/(2A)} \del^{-1/A}\}}\right)^{5/3}\left(\log\log\log D_L^{2^{\delta/2}}\right)^{1/3} \right\}.
 \eeq
 Note that $2^{\del/2} \leq 2$, $\log (2\del^{-1}) + 1 \leq \del^{-1}$ for all $\del \leq 1/4$, and $c_2^{1/2A}\del^{-1/A} \leq 6 c_1^{-1/(2A)} n_L \del^{-1/A}$. Thus upon comparing (\ref{E3E4_both}) and (\ref{SxT_E5aQ_main}), we see that (\ref{x_lower_bound_Q}) suffices with $\kappa_i''$ defined as above (specialized to the case $k=\Q$); the case of other fields $k$ follows from analogous computations.
 \end{proof}

\begin{lemma}\label{lemma_E1E2_bounds}
 Let $k$ be a fixed number field. Let $A \geq 2$ be fixed and let $0 < \del \leq 1/(2A)$ be a fixed positive constant.
 Let $L/k$ be a Galois extension of number fields with $\Gal(L/k) \simeq G$. Set $T = (\log D_L)^{2/\del}$.
 Given any absolute constant $0<c_1' \leq 1$, there exists a constant $D_2$ such that for $D_L \geq D_2$, 
\beq\label{E1E2_bound}
 |E_1| + |E_2| \leq 6c_1'   x (\log x)^{-(A-1)}
 \eeq
for all
\beq\label{x_lower_bound2}
(c_1')^{-1}  10^A n_L^{A+1}(\log D_L)^{2A+1}  \leq x \leq \exp( 10 n_L (\log D_L)^2).
\eeq
\end{lemma}
\begin{proof}
This is proved by simple computations checking error terms   in the range of  ``small'' $x$, that is  $x \leq \exp( 10 n_L (\log D_L)^2)$, and recalling $T=(\log D_L)^{2/\delta}$.
Writing $E_1 = E_{1,a} + \cdots + E_{1,d}$ and $E_2 = E_{2,a} + E_{2,b}$ we see that for example 
$|E_{1,a}| \leq c_1'  x (\log x)^{-(A-1)}$ for $x \leq \exp( 10n_L (\log D_L)^2)$ if $\del < 2/(2A+1)$
and 
\beq\label{E1a_D_main}
 D_L \geq  \exp \{( c_1'^{-1}(10n_L )^A)^{(1/\del - 2A-1)^{-1}} \}.
  \eeq
Similarly, $E_{1,b}, E_{1,c}, E_{2,a}$ are seen to be sufficiently small if $x$ is bounded below as in (\ref{x_lower_bound2}). 
The remaining two terms $E_{1,d}$ and $E_{2,b}$ impose (respectively) the constraints $\del < 2/(2A+1)$ and 
\beq\label{E1b_D_main}
 D_L \geq \exp \{(2 \cdot 10^{A}   n_L^{A+1} c_1'^{-1} )^{(2/\del - 2A-1)^{-1}}  \},
 \eeq
 and $\del< 1/(A+1)$, 
 \beq\label{E2b_D_main}
D_L \geq \exp \{ (10^{A+1} n_L^{A+2} c_1'^{-1})^{(2/\del - 2(A+1))^{-1}}\}.
\eeq
It suffices to assume $\del \leq 1/(2A)$ and to take $D_2=D_2(c_1',\del,n_L,A)$ to be the maximum of (\ref{E1a_D_main}), (\ref{E1b_D_main}) and (\ref{E2b_D_main}).
\end{proof}

\subsection{Proof of Proposition  \ref{prop_assume_zerofree2}}

To deduce Proposition  \ref{prop_assume_zerofree2}, for a fixed absolute constant $c_0$, from these lemmas, we will apply Lemma \ref{lemma_Sxt_bounds}  and Lemma \ref{lemma_E1E2_bounds} with the respective choices
\beq\label{c1_choice}
c_1 = c_0/(6C_5C_6), \qquad  c_1' = c_0/ (12C_5),
\eeq
  where $C_5$ and $C_6$ are the absolute constants arising in (\ref{eqn_C5}) and (\ref{Sxt_terms}) from the Lagarias-Odlyzko argument. After this choice in Lemma \ref{lemma_E1E2_bounds}, we  could denote the dependencies of $D_2(c_1',\del,n_L,A)$ by $D_2(c_0,C_5,\del,n_L,A)$.
The last step of the proof of  Proposition \ref{prop_assume_zerofree2} is to check that we can ensure that the parameters   are such that  (\ref{x_lower_bound2}) holds whenever (\ref{x_lower_bound}) (or (\ref{x_lower_bound_Q}) respectively) is satisfied. 
Note that the lower bound in (\ref{x_lower_bound2})  will hold if we have
\[ x \geq \exp \{ (2A+1) \log \log (D_L^{(c_1')^{-1/(2A+1)} 10^{1/2} n_L}) \}
	\geq \exp \{ (2A+1) \log \log (D_L^{((c_1')^{-1} 10^A n_L^{A+1})^{1/(2A+1)}}) \}.
\]
Thus either for general $k$ or $k=\Q$ it suffices to set
\beq\label{kappa_dfn_prime}
\kappa_1' = \kappa_1'' \geq 1, \quad \kappa_2' = \kappa_2''= \max \{ \kappa_2'', 2A+1\}, \quad \kappa_3' = (c_1')^{-1/(2A+1)}  \kappa_3''.
\eeq

\subsection{Partial summation back to prime counting}\label{partialsummation}

There are two remaining steps to pass from Proposition  \ref{prop_assume_zerofree2} to Theorem \ref{thm_assume_zerofree1} (in the regime of small $x$). 
First, we define the function 
\[ 
\theta_{\mathscr{C}} (x,L/k) = \sideset{}{'}\sum_{\tstack{\pfrak}{\Nm_{k/\Q} \pfrak \leq x}{\left[ \frac{L/k}{\pfrak} \right] = \mathscr{C}}} \log( \Nm_{k/\Q} \pfrak) 
= \sideset{}{'}\sum_{\bstack{\pfrak}{\Nm_{k/\Q} \pfrak \leq x}} \mathbf{1}_{\mathscr{C}} (\sig_\qfrak ) \log( \Nm_{k/\Q} \pfrak) ,
\]
in which the sum is restricted to prime ideals $\pfrak \subset \Ocal_k$ that are unramified in $L$, and we fix any prime ideal $\qfrak$ in $\Ocal_L$ above $\pfrak$ and let $\mathbf{1}_{\mathscr{C}}$ detect whether the conjugacy class of the Frobenius element $\sig_\qfrak$ is $\mathscr{C}$. 
A  Chebyshev argument shows that $\theta_{\mathscr{C}} (x,L/k)$ is well-approximated by $\psi_{\mathscr{C}}(x,L/k)$ and then partial summation passes from $\theta_{\mathscr{C}}(x,L/k)$ back to $\pi_{\mathscr{C}}(x,L/k)$; we only mention the highlights. 
We note that
\[
\psi_{\mathscr{C}} (x,L/k) - \theta_{\mathscr{C}} (x,L/k) = \sideset{}{'}\sum_{\bstack{\pfrak, m \geq 2}{ \Nm_{k/\Q} \pfrak^m \leq x}} \mathbf{1}_{\mathscr{C}} ( \sig_\qfrak^m ) \frac{1}{m} \log( \Nm_{k/\Q} ( \pfrak^m)),
\]
so that upon setting $m$ to be the smallest integer such that $x^{1/m} \geq 2$ (so in particular $m \leq \log x / \log 2$), the above difference is at most 
\[
	 \frac{\log x}{\log 2} \left( 	\frac{1}{2} \pi(x^{1/2},L/k)	+ \cdots + \frac{1}{m} \pi (x^{1/m},L/k) \right)
		\leq \frac{3}{2 \log 2} n_k x^{1/2} \log x,
	\]
	where we have denoted by $\pi(x,L/k)$ the  counting function for prime ideals (unramified in $L$) with $\Nm_{k/\Q} \pfrak \leq x.$
 Thus we see that the statement of Proposition \ref{prop_assume_zerofree2} holds for $\theta_{\mathscr{C}} (x,L/k)$ in place of $\psi_{\mathscr{C}}(x,L/k)$, with an additional error term of size at most $3 n_k x^{1/2} \log x$, 
which is no bigger than $c_0|\mathscr{C}||G|^{-1} x(\log x)^{-(A-1)}$ (for an absolute constant $c_0\leq 1$ we will choose later) as soon as the sufficient condition  
$ 3 |G|n_k=3 n_L \leq c_0 x^{1/2} (\log x)^{-A}$
 is met. It is simple to check that this holds in the regimes (\ref{x_range_k}) or (\ref{x_range_Q}) we consider in Proposition \ref{prop_assume_zerofree2}, with the parameters $\kappa_i'$ as already defined. 
 Thus for $x$ in either range we have 
\beq\label{theta_expansion}
 \left| \theta_{\mathscr{C}}(x,L/k) - \frac{|\mathscr{C}|}{|G|} x  \right| \leq 2c_0\frac{|\mathscr{C}|}{|G|} \frac{x}{(\log x)^{A-1}} .
 \eeq

Let $x_0$ denote the lower bound for $x$ in (\ref{x_range_k}) for general $k$ and for $x$ in the range (\ref{x_range_Q}) for $k=\Q$, respectively.
To pass from $\theta_{\mathscr{C}}(x)$ to $\pi_{\mathscr{C}}(x)$ (temporarily suppressing the notational dependence on $L/k$ for simplicity), we let $\lam_n$ be an increasing sequence of positive real numbers running over the norms $\Nm_{k/\Q}(\pfrak)$ attained by prime ideals of $k$ (unramified in $L$). By partial summation, for any $x_0 \leq x \leq \exp\{ 10 n_L (\log D_L)^2\}$, 
\beq\label{pi_theta}
 \pi_{\mathscr{C}}(x)
	= \sum_{\lam_n \leq x} \left(\sideset{}{'}\sum_{\Nm_{k/\Q} \pfrak = \lam_n}  \mathbf{1}_{\mathscr{C}} (\sig_\qfrak)  \log \lam_n\right) (\log \lam_n)^{-1}
= \int_{\lam_1}^x \frac{\theta_\Cscr(t)dt}{t \log ^2 t}  + \frac{\theta_\Cscr(x)}{\log x}.
\eeq
We split the integral into the region $\lam_1 \leq t \leq x_0$, in which the asymptotic (\ref{theta_expansion}) has not been verified, and the region $x_0 \leq t \leq x$, in which it has. For the first portion of the integral we apply the trivial bound $\theta_\Cscr(t) \leq n_k t \log t$ to see that this integral contributes at most $n_k \Li (x_0)$. 
In the remaining contributions to (\ref{pi_theta}), we may replace $\theta_\Cscr(t)$ by $|\Cscr||G|^{-1}t$ as in (\ref{theta_expansion}) (deferring the error terms for a moment), and similarly for $\theta_\Cscr(x)$; this main contribution becomes after integration by parts
\beq\label{pi_theta_t}
 \frac{|\Cscr|}{|G|} \left[ \int_{x_0}^x t \frac{d}{dt} \left( - \frac{1}{\log t} \right) dt + \frac{x}{\log x} \right]=  \frac{|\Cscr|}{|G|}  \left[ \Li(x) - \left(\Li(x_0) - \frac{x_0}{ \log x_0}\right) \right].
\eeq
The error terms accrued via this replacement are (in absolute value) at most
\beq\label{pi_theta_error2}
2c_0 \frac{|\Cscr|}{|G|}\int_{x_0}^{x} \frac{dt }{ (\log t)^{A+1}}  + 2c_0 \frac{|\Cscr|}{|G|} \frac{x}{(\log  x)^A}.
\eeq
In the first term of (\ref{pi_theta_error2}) we may bound the contribution from, say, $x_0 \leq t \leq x^{1/2}$ trivially by $2c_0 |\Cscr||G|^{-1} x^{1/2}$ while in the remaining portion we have $\log t \geq (1/2) \log x$, yielding a total contribution of at most $2^{A+2} c_0 |\Cscr||G|^{-1} x(\log x)^{-(A+1)}$; we trivially dominate this from above by $2^{A+2}c_0 |\Cscr||G|^{-1} x(\log x)^{-A}$ so that we may combine it with the second term in (\ref{pi_theta_error2}). Finally, we crudely bound the last two terms in (\ref{pi_theta_t}), in absolute value, by $2\Li(x_0)$. 	
In total, we have represented 
$ \pi_\Cscr(x)$ as $|\Cscr||G|^{-1} \Li(x)  + E $
where
\begin{align}
 |E| &  \leq (n_k +2) \Li(x_0) + 2c_0 |\Cscr||G|^{-1} x^{1/2}  + (2^{A+2} +2)c_0 |\Cscr||G|^{-1} x ( \log  x)^{-A} \nonumber \\
	& \leq (n_k+2) \Li(x_0) + (2^{A+2} +4)c_0|\Cscr||G|^{-1}x (\log  x)^{-A}. \label{PNT_2terms}
	\end{align}
Here we have used that $x^{1/2} \leq x (\log x)^{-A}$ in the regime of $x \leq \exp\{ 10 n_L (\log D_L)^2\}$ as soon as 
$x \geq \exp \{4A \log \log (D_L^{10^{1/2} n_L^{1/2}})\}$, which holds for all $x \geq x_0$.
The first term  on the right-hand side of (\ref{PNT_2terms}) is certainly dominated by the second as long as
\beq\label{x_lower_partialsum}
\frac{x}{(\log x)^A} 
  \geq \frac{|G|(n_k+2)}{(2^{A+2}+4)|\Cscr| c_0} \Li(x_0),
\eeq
for which it suffices to have $x \geq n_L c_0^{-1} x_0 (\log x)^A.$
Of course, we are already assuming that $x \geq x_0$; recalling that we presently only consider 
$x \leq \exp\{ 10 n_L (\log D_L)^2\}$
we see that (\ref{x_lower_partialsum}) holds as long as
\beq\label{x_lower_partialsum2}
x
  \geq  10^A n_L^{A+1} c_0^{-1} x_0 (\log D_L)^{2A} = 10^A n_L^{A+1} c_0^{-1}  \exp \{ 2A \log \log D_L \} \cdot x_0.
\eeq
Under this condition, we have shown that 
\[|E| \leq 2 (2^{A+2} +4) c_0 |\Cscr| |G|^{-1} x ( \log x)^{-A} \leq |\Cscr| |G|^{-1} x ( \log x)^{-A},\]
 upon making the choice 
 \beq\label{c_0_choice}
 c_0 = ( 2^{A+3} + 8)^{-1}.
 \eeq
We may accommodate the requirement (\ref{x_lower_partialsum2}) simply by enlarging the parameters $\kappa_i'$ by setting
$
 \kappa_1 = c_0^{-1} \kappa_1', $ $\kappa_2 = \kappa_2' + 2A$, $\kappa_3= \kappa_3' \geq 1.
$
We record the definitions here, with $c_0$ as in (\ref{c_0_choice}):
\begin{align}\label{nu_dfn_final_ref}
\kappa_1 &=  c_0^{-1} (6 (\frac{c_0}{12C_5C_6})^{-1}  10^{A-1} n_L^{A})^{1/\del_0} \del^{-2/\del_0} 	 \\
\kappa_2 &=  \max \{ 2A \del_0^{-1}, 4A c_k^{-1} n_k^3 \del^{-1}\}  +2A	\nonumber \\
\kappa_3  &= 6 (\frac{c_0}{12C_5})^{-1/(2A+1)}(\frac{c_0}{12C_5C_6})^{-1/(2A)} D_k n_L \del^{-1/A}.	\nonumber
\end{align}
To conclude, for $x$ in the ranges (\ref{x_range_k}) and (\ref{x_range_Q}) with $\kappa_i'$ replaced by $\kappa_i$, we have verified the effective error term in the asymptotic for $\pi_{\mathscr{C}}(x,L/k)$.
This completes the treatment of small $x$, and combining this with the result of \S \ref{sec_large_x} for large $x$, we may  conclude that Theorem \ref{thm_assume_zerofree1} holds.
\begin{remark}\label{remark_D0}
The threshold $D_0=D_0(\del, c_k, \be_0^{(k)}, n_L,C_1,C_2,A)$ appearing in Theorem \ref{thm_assume_zerofree1} is the maximum of $D_0'$ in (\ref{DL_lower_from_zerofree}), $D_1$  in (\ref{delDdel}),  $D_1'$ defined in  \S \ref{sec_large_x}, and  $D_2$  defined as the most restrictive of (\ref{E1a_D_main}), (\ref{E1b_D_main}), (\ref{E2b_D_main}) (with the  imposed choices $c_1'  = c_0/12C_5$ and $c_0 = ( 2^{A+3} + 8)^{-1}$).
\end{remark}

\subsection{Remark: A Chebotarev theorem for fields without quadratic subfields}\label{sec_app_Cheb_quad}
In the introduction, we stated that one of our two goals was to remove the $\be_0$ term in Theorem \ref{thm_LagOdl_uncond}.
As an aside, we note that for certain fields, the existence of an exceptional zero can already be ruled out, so that an immediate application of Theorem \ref{thm_LagOdl_uncond} yields:
\begin{thm}\label{thm_Stark_Cheb}
Let $k$ be a number field such that $\zeta_k(s)$ has no real zeroes. Let $L/k$ be a Galois extension of relative degree at least $3$ such that $L/k$ contains no quadratic extension of $k$. Then there exist absolute effectively computable constants $C_1, C_2$ such that   for all $
x \geq \exp(10n_L(\log D_L)^2),
$
\beq\label{LO_asymp_app1}
  \left| \pi_{\mathscr{C}}(x,L/k) - \frac{|\mathscr{C}|}{|G|} \Li(x) \right| \leq  C_1 x \exp(-C_2 n_L^{-1/2} (\log x)^{1/2} ).
  \eeq

\end{thm}
\begin{remark}
In particular, if $k=\Q$ this theorem holds unconditionally for any $L/\Q$  such that $G = \Gal(L/\Q)$ with $|G| \geq 3$ has no subgroup of index 2  (for example, $G \simeq C_p$ for $p$ an odd prime).
\end{remark}

Theorem \ref{thm_Stark_Cheb} is an application of a  nice idea of Stark \cite[Theorem 3]{Sta74}, in turn a refinement of a theorem of Heilbronn \cite{Hei72}. (See also further work on eliminating Siegel zeroes in towers of fields in \cite{Mur99,OdlSki93}.)

\begin{letterthm}[{\cite[Theorem 3]{Sta74}}]\label{thm_Stark}
Let $L$ be a Galois extension of $k$ with $\Gal(L/k) \simeq G$ and let $\chi$ be a character of $G$. Suppose $\rho$ is a simple zero of $\zeta_L(s)$. Then $L(s,\chi, L/k)$ is analytic at $s=\rho$. Furthermore, there is a field $F$ with $k \subseteq F \subseteq L$ such that $F/k$ is cyclic and for any field $E$ with $k \subseteq E \subseteq L$, $\zeta_E(\rho)=0$ if and only if $F \subseteq E$. If in particular $\rho$ is real, then either $F=k$ or $F$ is quadratic over $k$.
\end{letterthm}

By Theorem \ref{thm_LagOdl_uncond}, we need only consider a possible real zero of $\zeta_L(s)$, which by Theorem \ref{thm_Stark} (and the assumption that $\zeta_k(s)$ has no real zero) can only occur if there is a quadratic extension $F$ of $k$ contained in $L$. No such $F$ can exist  if $\Gal(L/k)$ has no index 2 subgroup.  Nevertheless, as remarked before, the lower bound on $x$ in Theorem \ref{thm_Stark_Cheb} is too large for our ultimate application to $\ell$-torsion, a problem which  Theorem  \ref{thm_assume_zerofree1} alleviates via careful attention to the assumed box-shaped zero-free region.

\section{A zero density result for families of Dedekind zeta functions}\label{sec_zero_density}

We have proved a Chebotarev density theorem conditional on a box-shaped zero-free region for $\zeta_{L}(s)/\zeta_k(s)$. Now we restrict our attention to $k=\Q$ and show that within appropriate families of Galois extensions of $\Q$, except for a possible exceptional subfamily of density zero within the family, each $\zeta_L(s) / \zeta(s)$ is in fact zero-free in the desired region.
To do so we will build on the result of Kowalski and Michel \cite[Thm. 2]{KowMic02} on the density of zeroes among a family of cuspidal automorphic $L$-functions. We describe our approach somewhat generally to facilitate future applications, and then specialize to our present setting.

\subsection{The Kowalski-Michel zero density estimate}\label{sec_KM}

Let $m \geq 1$ be fixed. For any cuspidal automorphic representation $\rho$ of $\GL(m)/\Q$, define the zero-counting function for the corresponding automorphic $L$-function $L(s,\rho)$ in a region with $\al \in [1/2,1]$, $T \geq 0$ by 
\[ N(\rho;\al,T) = | \{ s  = \be + i \ga : \be \geq \al, |\ga| \leq T, L(s,\rho)=0\}|,\]
counting with multiplicity. For an isobaric representation $\pi = \rho_1 \boxplus \cdots \boxplus \rho_r$ with $\rho_j$ cuspidal, define
\beq\label{zeroes_sum_up}
 N(\pi;\al,T) = N(\rho_1;\al,T) + \cdots + N(\rho_r;\al,T),
 \eeq
again counting each zero with multiplicity.

The main outcome of \cite{KowMic02} is a bound for $N(\rho;\al,T)$ that holds on average for an appropriate family of cuspidal representations $\rho$. Our innovation is to develop a means to apply their results to the case when $\pi$ varies over an appropriate family of isobaric representations, in our case, obtained from Dedekind zeta functions.
We first recall the original setting for cuspidal representations, which assumes the following conditions hold:

\begin{condition}\label{KM_condition}
For each $X \geq 1$ let $S(X)$ be a finite (possibly empty) set of cuspidal automorphic representations $\rho$ of $\GL(m)/\Q$ such that the following properties hold for $(S(X))_{X \geq 1}$:
\begin{enumerate}[label=(\roman*)]

\item\label{item1} Every $\rho \in S(X)$ satisfies the Ramanujan-Petersson conjecture  at the finite places.
\item\label{item2} There exists $A>0$ and a constant $M_0$ such that for all  $X \geq 1$, for all $\rho \in S(X)$, 
$\Cond(\rho)\leq M_0 X^{A}.$
\item\label{item3} There exists $d>0$ and a constant $M_1$ such that for all $X \geq 1$,
$ |S(X)| \leq M_1 X^{d}.$
\item\label{item4}
 For any $\ep >0$ there exists a constant $M_{2,\ep}$ such that for all $\rho \in S(X)$ we have the convexity bound
  \[|L(s, \rho)| \leq M_{2,\ep} (\Cond(\rho)(|t|+2)^{m})^{(1-\Re(s))/2+\varepsilon}, \qquad \text{ for $0 \le \Re(s) \le 1$.}\] 
For any $\ep >0$ there exists a constant $M_{3,\ep}$ such that for all $\rho \not\simeq \rho' \in S(X)$ 
  we have the convexity bound
  \[|L(s, \rho \otimes \rho')| \leq M_{3,\ep} (\Cond(\rho \otimes \rho')(|t|+2)^{m^2})^{(1-\Re(s))/2+\varepsilon}, \qquad \text{ for $0 \le \Re(s) \le 1$.}\] 
\end{enumerate}
\end{condition}
\begin{remark}
Kowalski and Michel call  $(S(X))_{X \geq 1}$  a \emph{family} of automorphic representations, with associated automorphic $L$-functions; following their convention we will call the associated collection of constants $\{m,A,d,M_0,M_1,M_{2,\ep},M_{3,\ep}\}$ the family parameters. 
\end{remark}

\begin{remark}
It is worth comparing precisely Condition \ref{KM_condition} to the hypotheses originally stated in the work of \cite{KowMic02}.
We note that  the above criteria \ref{item1} -- \ref{item3} reduce to exactly the criteria of  \cite[Thm. 2]{KowMic02}; 
Condition \ref{item4} above replaces their assumption that all the $L$-functions in $(S(X))_{X \geq 1}$ have the same gamma factors at infinity. That condition is only used in order to attain the uniform convexity bounds of \cite[Lemma 10]{KowMic02} (Kowalski, personal communication), and thus we merely assume the relevant uniform convexity bounds directly. 
 \end{remark}

In this context, we recall Kowalski and Michel's original theorem: 
\begin{letterthm}[{\cite[Theorem 2]{KowMic02}}]\label{thm_KM_original}
Let $(S(X))_{X \geq 1}$ be a family of cuspidal automorphic representations of $\GL(m)/\Q$ satisfying Condition \ref{KM_condition}. Let $\al \geq 3/4$ and $T \geq 2$. Then there exists a constant $c_0' = c_0'(m,A,d)$, in particular
\beq\label{c0'_dfn}
 c_0' = \frac{5mA}{2} + d,
 \eeq
and a constant $B \geq 0$, depending only on the family parameters, such that for every choice of $c_0 > c_0'$ we have that there exists a constant $M_{4,c_0}$ depending only on $c_0$ such that for all $X \geq 1$,
\[ \sum_{\rho \in S(X)} N(\rho;\al,T) \leq M_{4,c_0} T^B X^{c_0 \frac{1 - \al }{2\al -1}}.\]
\end{letterthm}

\subsection{Defining a family of automorphic representations}\label{sec_families_aut}
Fix $n \geq 2$  and a transitive subgroup $G \subseteq S_n$. Let $\mathscr{F}(\Q,G) \subset Z_{|G|}(\Q,G)$ be a set  of Galois extensions $L/\Q$ with $\Gal(L/\Q) \simeq G$, 
and let $\mathscr{F}(\Q,G;X)$ denote the finite subset comprised of those fields  
 with
$D_L =| \Disc L/\Q| \leq X$. (Momentarily we will construct such a set from each of the families $Z_n^\Iscr(\Q,G)$ of degree $n$ fields considered in our main theorems.)

 Denote the irreducible representations of $G$ by $\rho_0, \rho_1,\ldots ,\rho_s,$ with $\rho_0$ being the trivial representation. 
For each field $L \in \mathscr{F}(\Q,G;X)$, the Dedekind zeta function may be written as
 \beq\label{zeta_factor_text}
  \zeta_L(s) = \zeta(s)  \prod_{j=1}^s L(s,\rho_j,L/\Q)^{m_j}.
  \eeq
 The regular representation, of total dimension $|G|=1 + \sum_{1 \leq j \leq s} m_j^2$, may be written as an isobaric sum
\[ \reg_G = \rho_0 \boxplus (\rho_{1} \boxplus \cdots \boxplus \rho_{1}) \boxplus \cdots \boxplus
	(\rho_{s} \boxplus \cdots \boxplus \rho_{s})
\]
 in which $\rho_j$ appears $m_j = \dim \rho_j$ times.  
Thus for each field $L \in \mathscr{F}(\Q,G;X)$ we may consider the Artin $L$-function
$L(s,\pi)= \zeta_{L}(s)/\zeta(s) $
for the representation  
\beq\label{L_pi}
\pi = (\rho_{1} \boxplus \cdots \boxplus \rho_{1}) \boxplus \cdots \boxplus
	(\rho_{s} \boxplus \cdots \boxplus \rho_{s})
	\eeq
in which $\rho_j$ appears $m_j = \dim \rho_j$ times.
Additionally, assuming the Strong Artin Conjecture (see \S \ref{sec_SC_KM}), to each field $L \in \mathscr{F}(\Q,G;X)$ and each representation $\rho_j$, there is an associated cuspidal automorphic representation $\pi_{L,j}$ of $\GL(m_j)/\Q$; we then have 
\[ L(s,\pi_{L,j}) = L(s,\rho_j, L/\Q).\]

Now fix $1 \leq j \leq s$. For each $ X \geq 1$, let $\Lscr_j(X)$ be the set of cuspidal automorphic representations $\pi_{L,j}$ of $\GL(m_j)/\Q$ associated by the Strong Artin Conjecture to the fields $L \in \mathscr{F}(\Q,G;X)$ and the  representation $\rho_j$.

The main  result of this section, and the key result underlying our effective Chebotarev theorem in families, relates to the following phenomenon. For each $j$, under appropriate assumptions, we show that   Theorem \ref{thm_KM_original} applies to the family $(\Lscr_j(X))_{X \geq 1}$, so that for each $X \geq 1$, aside from very few possible ``bad'' exceptional representations, for each representation $\pi \in \Lscr_j(X)$ the associated $L$-function $L(s,\pi)$ possesses a certain zero-free region. Now a key difficulty arises: in general, depending on the group $G$ and the family $\Fscr(\Q,G;X)$, it could happen that a given $L$-function $L(s,\pi)$ corresponding to a representation $\pi \in \Lscr_j(X)$ occurs as a factor in $\zeta_{L}(s)/\zeta(s)$ for ``many'' fields $L \in \Fscr(\Q,G;X)$, indeed even possibly a positive proportion of such fields (see \S \ref{sec_rationale_ram}). We need to rule out this possibility that a ``bad'' exceptional representation in $\Lscr_j(X)$ could lead to an $L$-function that ``contaminates'' $\zeta_L/\zeta$ for a positive proportion of fields in $\Fscr(\Q,G;X)$. 
In this section, we state appropriate conditions on a set  $\Fscr(\Q,G;X)$ of Galois extensions that allow us to rule out this problem (see in particular the condition (\ref{item5_eqn}) below). In \S \ref{sec_app_Dedekind}, we show that the families of fields that we consider in our main theorems obey these conditions.

Now we state the conditions we assume on the set   $\Fscr(\Q,G)$ of Galois extensions and the associated families $(\Lscr_j(X))_{X \geq 1}$ of automorphic representations ($1 \leq j \leq s$), building on Condition \ref{KM_condition}. (Note that we explicitly assume the Strong Artin Conjecture below, but for certain groups $G$ it is known; see \S \ref{sec_SC_KM}.)
\begin{condition}\label{KM_condition_composed}
Let $\Fscr(\Q,G)$ be a set of Galois extensions as specified above. 
For each $1 \leq j \leq s$ and each $X \geq 1$, define the set $\Lscr_j(X)$ of automorphic representations as above, assuming the Strong Artin Conjecture.

Assume that for each $1 \leq j \leq s$, the family $(\Lscr_j(X))_{X \geq 1}$ satisfies Condition \ref{KM_condition}, with corresponding parameters $\{m_j,A_j,d_j,M_{0,j},M_{1,j},M_{2,j,\ep},M_{3,j,\ep}\}$. In particular, for $1 \leq j \leq s$, $(\Lscr_j(X))_{X \geq 1}$ is a family in the sense of Theorem \ref{thm_KM_original}.

Let $A  \geq 0$, $M_0$ be such that for all $X \geq 1$, for every 
field $L \in \Fscr(\Q,G;X)$, the representation $\pi$ associated to $L(s,\pi) = \zeta_L(s)/\zeta(s)$ has $\Cond(\pi) \leq M_0 X^A$.

Let   $d,M_1$ be such that for all $X \geq 1$, $|\Fscr(\Q,G;X)| \leq M_1 X^d$.

We assume that for each $1 \leq j \leq s$, 
there exists $0 \leq \tau_j < d$ and a constant $M_{5,j}$   such that for all $X \geq 1$, for any fixed $\pi \in  \Lscr_j(X)$, 
\beq\label{item5_eqn}
 \#\{ L \in \Fscr(\Q,G;X) : \pi_{L,j} = \pi \} \leq M_{5,j} X^{\tau_j}.
 \eeq
\end{condition}
We will call  $\{M_0,M_1,A,d\}$ and $\{ m_j,A_j,d_j,M_{1,j},M_{2,j,\ep},M_{3,j,\ep}, M_{5,j}\}$ for $1 \leq j \leq s$ the family parameters for $\Fscr(\Q,G)$.

\subsection{A zero density theorem for $L$-functions associated to the family $\Fscr(\Q,G)$}

To bound on average the number of zeroes of $L$-functions $\zeta_L(s)/\zeta(s)$ in a certain region, as the field $L$ varies over the family $\Fscr(\Q,G)$, we will apply Theorem \ref{thm_KM_original} repeatedly,
under the assumption of Condition \ref{KM_condition_composed}. 
\begin{thm}\label{thm_KM_badfamily}
Let $\Fscr(\Q,G)$ be a set of Galois extensions as specified above. 
For each $1 \leq j \leq s$ and each $X \geq 1$, define the set $\Lscr_j(X)$ of automorphic representations as above, assuming the Strong Artin Conjecture. 
 
Assume that $\Fscr(\Q,G)$ and the families $(\Lscr_j(X))_{X \geq 1}$ for $j=1,\ldots, s$, satisfy Condition \ref{KM_condition_composed}.

 Set $\tau = \max_j \tau_j$ and $m= \max m_j$. 
 Then for any $0<\Del<1$ sufficiently small that $\Del < 1-\tau/d$, and for any $\eta<1/4$, there exists  $B$ depending only on the family parameters for $\Fscr(\Q,G)$, and $0<\del \leq 1/4$ depending only on $A,m,d,\Del,\tau$, such that for all $X \geq 1$, at most $O(X^{ (1 - (1-\eta)\Del)d})$ fields $L \in \Fscr(\Q,G;X)$ can have the property that $\zeta_{L}(s)/\zeta(s)$  has a zero in the region 
\[
[1-\del,1] \times[-X^{\eta \Delta d/B}, X^{\eta \Delta d/B}].
\]
The implied constant in the $O(\cdot)$ notation depends only on $A,m,d,\Del,\tau$ and $s$ (the number of nontrivial irreducible representations of $G$).

\end{thm}

\begin{remark}\label{remark_Del_0}
 We see that in the hypotheses there is a non-empty range of $0<\Del< 1 - \tau/d$ since each $\tau_j < d$.
 \end{remark}

To deduce Theorem \ref{thm_KM_badfamily} we first apply Theorem \ref{thm_KM_original} to the family $(\Lscr_j(X))_{X \geq 1}$ for each $1 \leq j \leq s$.   Let $1 \leq j \leq s$ be fixed. By Theorem \ref{thm_KM_original}, for any $\al_j \geq 3/4$ and $T_j \geq 2$, for all $X \geq 1$,
\beq\label{Nsum}
  \sum_{\pi \in \Lscr_j(X)} N(\pi;\al_j,T_j) \ll_{c_{j,0}}  T_j^{B_j} X^{c_{j,0} \frac{1-\al_j}{2\al_j-1}},
\eeq
in which  we may choose any $c_{j,0} > c_{j,0}'$, with $c_{j,0}' = c_{j,0}'(m_j,A_j,d_j)$ as shown to exist in Theorem \ref{thm_KM_original}; the particular form is not critical, but we may for example take
\[ c_{j,0}' = \frac{5m_jA_j}{2} + d_j.\]
In the spirit of \cite[Remark 3]{KowMic02}, we pause to observe that although the parameter $d_j$ assumed to exist  in the upper bound \ref{item3} of Condition \ref{KM_condition} may not provide a sharp upper bound, this does not cause any contradictions in terms of its role in $c_{j,0}'$; if $d_j$ is an over-estimate, then the right-hand side of (\ref{Nsum}) is similarly an overestimate (and similarly with respect to the possibly non-sharp parameter $A_j$). 
Indeed,  for convenience we may choose $c_{j,0}  =  c_{j,0}'' + \ep_1$ (for a certain $\ep_1$ to be chosen later) with 
\beq\label{cj0''}
 c_{j,0}'' = \frac{5m_jA}{2} + d.
 \eeq
Note that $A \geq \max_j A_j, d \geq \max_j d_j$ so that this choice is valid.

 Set $\tau = \max_{1 \leq j \leq s} \tau_j$. Recalling that $\Del$ is given, we fix $\al_j$ to be such that 
 \[ \frac{c_{j,0} (1-\al_j)}{(2\al_j-1)} = (1-\Delta)d - \tau.\]
We see that the right-hand side is positive, so that $\al_j<1$,  since $\Del < 1 - \tau/d$. Theorem \ref{thm_KM_original} applies when $\al_j \geq 3/4$; if necessary one could simply impose this using monotonicity of the estimates, but in fact it is simple to check that this holds in our scenario. (This will also easily be satisfied in our ultimate applications, in which we will be working very close to the line $\Re(s)=1$.) We compute that 
 \[ \al_j = \frac{c_{j,0} + (1-\Del) d - \tau}{c_{j,0} + 2((1-\Del)d - \tau)},\]
 so that $\al_j \geq 3/4$ as long as 
 \beq\label{al_c_cond}
 c_{j,0} \geq 2((1-\Del)d - \tau).
 \eeq
By assumption, $\Del < 1 - \tau/d$; let $\ep_2>0$ be such that 
\beq\label{Del_choice1}
\Del =  1 - \tau/d - \ep_2/2d.
\eeq
Then (\ref{al_c_cond}) is equivalent to the requirement that $c_{j,0} \geq \ep_2$, which will always hold as long as we choose $\ep_1 \geq \ep_2$, according to the definition (\ref{cj0''}), upon recalling that $A,d \geq 0$.
Upon setting $T_j = X^{\eta \Delta d/B_j}$, we conclude that 
 \beq\label{Tj_orig}
   \sum_{\pi \in \Lscr_j(X)} N(\pi;\al_j,T_j) \ll_{c_{j,0}} X^{\eta \Delta d} X^{(1 - \Delta)d - \tau}
 	\ll_{c_{j,0}} X^{ (1 - (1-\eta)\Del)d - \tau}.
\eeq

Now we assemble these results together for $1 \leq j \leq s$. For notational convenience, given an $L$-function $L(s)$ (which could be an Artin $L$-function $L(s,\rho,L/\Q)$ or an automorphic $L$-function $L(s,\pi)$ corresponding to an automorphic representation $\pi$), we will let $N'(L(s);\al,T)$ denote the number of zeros $ \be+i\ga$ with $L(\be+i\ga)=0$, and $\be \geq \al, |\ga| \leq T$.
Set $\al = \max_j \al_j$ and $T = \min_j T_j$. 
(Note that $\al \geq 3/4$.) Then for each $X \geq 1$,  assuming the Strong Artin Conjecture,
\begin{align*}
 \sum_{L \in \mathscr{F}(\Q,G;X)} N'(\zeta_L/\zeta;\al,T) 
& = \sum_{L \in \mathscr{F}(\Q,G;X)}  \sum_{j=1}^s m_jN'(L(s,\rho_j,L/\Q);\al,T)  \\
& = \sum_{L \in \mathscr{F}(\Q,G;X)}  \sum_{j=1}^s m_jN'(L(s,\pi_{L,j},L/\Q);\al,T)  \\
& = \sum_{j=1}^s m_j \sum_{\pi \in \Lscr_j(X)} N'(L(s,\pi);\al,T)  \sum_{\bstack{L \in \mathscr{F}(\Q,G;X)}{\pi_{L,j} =\pi}} 1.
\end{align*}
Using condition (\ref{item5_eqn}), we can bound the right-hand side from above by
\[ \ll \sum_{j=1}^s m_j X^{\tau_j} \sum_{\pi \in \Lscr_j(X)} N'(L(s,\pi);\al,T) .\]
Thus by applying (\ref{Tj_orig}) for each $1 \leq j \leq s$, we see that
\[ \sum_{L \in \mathscr{F}(\Q,G;X)} N'(\zeta_L/\zeta;\al,T)  \ll_{c_0,s,m}X^{ (1 - (1-\eta)\Del)d},\] 
where $c_{0} = \max_j c_{j,0}$. 
From this we conclude that at most $O_{c_0,s,m}(X^{ (1 - (1-\eta)\Del)d})$ fields  $L \in  \mathscr{F}(\Q,G;X)$ 
can have the property that $\zeta_L(s)/\zeta(s)$ has a zero in the region
$
[\al,1] \times[-X^{\eta \Delta d/B}, X^{\eta \Delta d/B}],
$
where $B= \max B_j$. The implied constant depends on $c_0,s,m$, and hence on $A,m,d,\tau,\Del,s,\ep_1$.
Now from (\ref{Del_choice1}), $\ep_2$ is defined, and then we can choose $\ep_1 = \ep_2$ in the definition of $c_{j,0}$. Then we may compute that upon setting $\del = 1 - \al = 1 - \max_j \{ \al_j\}$ (which we have therefore verified satisfies $0<\del \leq 1/4$), we have
\beq\label{describe_del}
\del = \frac{\ep_2}{5 \max_j \{m_j\} A  + 2d + 4\ep_2} = \frac{\ep_2}{5 m A  + 2d + 4\ep_2}
\eeq
as an allowable choice. Since $\ep_2$ is determined by $\Del,\tau,d$ we can write the dependencies in terms of these parameters.
This yields the result of Theorem \ref{thm_KM_badfamily}, moreover with a specific description of $\del$.

\begin{remark}\label{remark_del_choice_1}
This argument shows that although the parameters $A,d$ are only assumed to yield valid upper bounds (not necessarily sharp) in Condition \ref{KM_condition_composed}, it is advantageous to make them as small as possible.
In a similar vein, it is worth asking why, if making $1-\Del$ smaller gives better control on the exceptional set, we do not in (\ref{Del_choice1}) artificially inflate the size of $d$. The reason is that $1-\Del$ only controls the density (roughly $O(X^{(1-\Del)d})$) of the exceptional set relative to the assumed upper bound $O(X^d)$ for the family; thus  in this instance also, it is advantageous to make $d$ as sharp as possible. 
\end{remark}

\begin{remark}
We see that the size of $\Del$, and hence of the possible exceptional set of bad fields in $\mathscr{F}(\Q,G;X)$ depends on the largest value of $\tau_j$ with $1 \leq j \leq s$ coming from the condition (\ref{item5_eqn}). The larger $\max_j \tau_j$ is, the smaller we must take $\Del$, and the less savings we have for the possible exceptional set in $\mathscr{F}(\Q,G;X)$.
\end{remark}

\begin{remark}\label{remark_Cho_Kim}
We recall that Cho and Kim (e.g. \cite[Theorem 3.1]{ChoKim12} and other works) have also applied \cite{KowMic02} to certain families of isobaric representations, say $\pi = \pi_1 \boxplus \cdots \boxplus \pi_r$  of $\GL(m)/\Q$, with $m=m_1 + \cdots + m_r$, and each $\pi_j$  a cuspidal automorphic representation of $\GL(m_j)/\Q$. Let us momentarily call the family of such $\pi$ by $S(X)$ and for each $j$ the family of such $\pi_j$ by $S_j(X)$. In their work, item \ref{item4} of Condition \ref{KM_condition} is replaced by the requirement that for each $1 \leq j \leq r$, for all $\rho_j \in S_j(X)$  the gamma factor of $L(s,\pi_j)$ is of the form $\prod_{i=1}^{m_j}\Gamma(s+\alpha_i)$, where $\alpha_i \in \mathbb{R}$ are fixed; this is a special case of the version of \ref{item4} stated here.
More importantly, instead of the key item (\ref{item5_eqn}) in Condition \ref{KM_condition_composed}, Cho and Kim assume that for any two inequivalent $\pi, \pi' \in S(X)$ with $\pi = \pi_1 \boxplus \cdots \boxplus \pi_r$ and $\pi' = \pi_1' \boxplus \cdots \boxplus \pi_r'$, they have $\pi_j \not\simeq \pi_k'$ for all $1 \leq j,k \leq r$.  Relative to (\ref{item5_eqn}), this would be the statement that for each $j$, for any fixed $\rho \in S_j(X)$, precisely one $\pi \in S(X)$ has $\pi_j \simeq \rho$, which in our notation is even stronger than the case $\tau_j=0$ for all $1 \leq j \leq r$. Cho and Kim used this to deduce that $|S_j(X)| = |S(X)|$ for each $j$, which was crucial to their proof, but also limited the types of families $S(X)$ they could consider. 
\end{remark}

%%%%%%%%%%%%%%%%%%%
\section{Verifying the conditions of the zero density theorem for families of Dedekind zeta functions}\label{sec_app_Dedekind}

The main result of this section is that Theorems \ref{thm_cheb_main_quant_pt1}, \ref{thm_cheb_main_quant_pt2}, \ref{thm_cheb_main_quant_pt3}, \ref{thm_cheb_main_quant_pt2'} and \ref{thm_cheb_main_quant_pt4} may be deduced from Theorem \ref{thm_KM_badfamily} by verifying that for each of the families of fields considered in these theorems, Condition \ref{KM_condition_composed} is satisfied. 
Accordingly, in this section we fix $Z_n^\Iscr(\Q,G)$ to be one of the families specified in the above theorems, 
under the associated hypotheses of the theorem (if any).

\subsection{Passage to a family of Galois closures}\label{sec_construct_Galois_family}

We now pass from considering the original set of the degree $n$ fields in $Z_n^\Iscr(\Q,G)$ to considering the set of Galois closures $\tilde{Z}_n^\Iscr(\Q,G) = \{ \tilde{K}: K \in Z_n^{\Iscr}(\Q,G)\}$;  each Galois closure corresponds to a constant number of fields in $Z_n^\Iscr(\Q,G)$ (only depending on $G$ as a permutation group).
We now recall the notation of \S \ref{sec_families_aut}. Using that notation, we define $\mathscr{F}(\Q,G) = \tilde{Z}_n^\Iscr(\Q,G)$ to be the set of Galois extensions we consider, and we accordingly define the sets $\Lscr_j(X)$ for each $1 \leq j \leq s$ and every $X \geq 1$, and thereby the corresponding families $(\Lscr_j(X))_{X \geq 1}$ of automorphic representations.

\subsection{Verification of Condition \ref{KM_condition} \ref{item1} -- \ref{item4}}\label{sec_SC_KM}
Now that we have constructed the appropriate families $\mathscr{F}(\Q,G) = \tilde{Z}_n^\Iscr(\Q,G)$ and $(\Lscr_j(X))_{X \geq 1}$ for each $1 \leq j \leq s$, we must verify that Condition \ref{KM_condition_composed} is satisfied.
We first note that for each family $\tilde{Z}_n^\Iscr(\Q,G)$ we consider, either the strong Artin conjecture is known to apply to all the Galois representations considered  (this is the case in Theorem  \ref{thm_cheb_main_quant_pt1}) or it is explicitly assumed (this is the case in Theorems  \ref{thm_cheb_main_quant_pt2}, \ref{thm_cheb_main_quant_pt3}, \ref{thm_cheb_main_quant_pt2'} and \ref{thm_cheb_main_quant_pt4}). To be precise,
let us write the Euler product of an Artin $L$-function as $L(s,\rho) = \prod_v L(s,\rho_v)$, and 
the Euler product for an automorphic $L$-function  as $L(s,\pi)=\prod_vL(s,\pi_v)$.

\begin{letterconjecture}[Strong Artin Conjecture] Let $L$ be a finite Galois extension of a number field $k$, with $\Gal(L/\Q) \simeq G$. Let $\rho$ be an $m$-dimensional complex representation of $G$. There exists an automorphic representation $\pi(\rho)$ of $\GL(m)/\Q$ such that the $L$-functions $L(s,\rho)$ and $L(s,\pi)$ agree almost everywhere, i.e. except at a finite number of places $v$, $L(s,\rho_v) = L(s,\pi_v)$. Moreover, if $\rho$ is irreducible, then $\pi$ is cuspidal. 
\end{letterconjecture}
This is known to hold for:
  1-dimensional representations $\rho$, due to Artin \cite{Art30};
 nilpotent Galois extensions $L/k$, due to Arthur and Clozel \cite{ArtClo};
 $A_4$ and $S_4$, due to Langlands \cite{Lan80} and Tunnell \cite{Tun81}, respectively;
 dihedral groups (and in particular $S_3$), due to Langlands \cite{Lan80}.
We also note that in the setting we will work in, a stronger identity is known. (See, for example, \cite[Theorem 4.6]{DelSer}, \cite[Proposition 2.1]{Martin}, \cite[Appendix A]{MarThesis}, and \cite[Proposition 1.5]{MarRam}.)
\begin{letterthm}
If $\pi$ is cuspidal and $L(s,\pi_v)=L(s,\rho_v)$ for almost all $v$, then in fact $L(s,\pi) = L(s,\rho)$.
\end{letterthm}
 These considerations guarantee that in the settings we consider (with the relevant hypotheses we assume), each $\Lscr_j(X)$ is a set of cuspidal automorphic representations.

We next confirm that for each $1 \leq j \leq s$, the family $(\Lscr_j(X))_{X \geq 1}$ satisfies the four items in  Condition \ref{KM_condition}. 
For item \ref{item1}, since the Ramanujan-Petersson conjecture holds for automorphic $L$-functions associated to Artin $L$-functions once they are known to exist (see e.g. the comment following \cite[Thm. 5]{KowMic02}), under the assumption (or known truth) of the strong Artin conjecture, the Ramanujan-Petersson conjecture holds for all the cuspidal automorphic $L$-functions in each set $\Lscr_j(X)$.

For item \ref{item2},  note that if $K \in Z_n^\Iscr(\Q,G;X)$ then by construction $D_K \leq X$. The following standard lemma
relates to discriminants of a field and its Galois closure.
\begin{lemma}[Discriminant comparisons]\label{lemma_compare_disc}
Let $G$ be a transitive subgroup of $S_n$. 
There exist constants $C_1=C_1(G)$ and $C_2=C_2(G)$ such that for every field $K \in Z_n(\Q,G)$,
\[ C_1 D_K^{|G|/n} \leq D_{\tilde{K}} \leq C_2  D_K^{|G|/2}.\]
\end{lemma}
(The lemma follows from Lemmas~\ref{lemma_power_p} and \ref{lemma_wild_primes}, recorded below, and for the left-hand inequality,  the fact that every cycle length in a permutation is at most the order of the permutation.)

Recall that in general for an Artin $L$-function $L(s,\rho,L/k)$, if  $F(\chi)$ denotes the Artin conductor of $\chi = \mathrm{Tr}(\rho)$, then the conductor of $L(s,\rho,L/k)$ is given by
$ A(\chi) = D_k^{\chi(1)}\Nm_{k/\Q}F(\chi).
$
According to the multiplicativity relation 
$D_L = D_k \prod_{\chi_j} A(\chi_j)^{\chi_j(1)}$ 
for the conductors in the identity (\ref{zeta_factor_text}), we see that for each $1 \leq j \leq s$, the conductors of $L(s,\rho_j,L/\Q)$ are bounded by $\ll_{n,G} X^{|G|/2}$ and we may take $A_j = A = |G|/2$ for all $j$.

For \ref{item3}, to control the size of the family of fields $\tilde{Z}_n^\Iscr(\Q,G;X)$ it suffices 
to control the sizes of the families $Z_n^\Iscr(\Q,G;X)$ (and moreover it suffices to bound from above the sizes of the families $Z_n(\Q,G;X)$ without the ramification restriction). 
Thus we may apply the following known unconditional upper bounds to show the existence of $d_j = d$ for all $j$: 
$G$ cyclic, Proposition \ref{prop_counting_cyclic_fields}; $G \simeq S_n$ see (\ref{EllVen_upper}); $G \simeq D_p$ see (\ref{Thorne_upper}); $G \simeq A_4$ see (\ref{Wong_upper}); $G \subseteq S_n$ simple, we simply embed $Z_n(\Q,G;X)$ in the family of all fields of degree $n$ and apply (\ref{EllVen_upper}).

For item \ref{item4}, we use the known convexity bounds for automorphic $L$-functions, which apply to our Artin $L$-functions under the strong Artin conjecture. Briefly, to be precise, we recall for $t\in \R$ the analytic conductor of $L(s,\pi)$ (in terms of the arithmetic conductor $\Cond(\pi)$ and the local parameters at infinity, $\mu_\pi(j)$),
\[Q_\pi(t)=\Cond(\pi)\prod_{j=1}^{m}(1+|it-\mu_\pi(j)|).
\]
 Then via the functional equation, Stirling's formula, and an application of the Phragmen-Lindel\"{o}f principle, one may derive the classical convexity bound (see e.g. \cite[page 5]{Harcos}):
\[L(s,\pi) \ll_{\pi,\varepsilon} Q_\pi(t)^{\frac{1-\Re(s)}{2}+\varepsilon}, \qquad 0 \le \Re(s) \le 1.
\]
For $\pi, \pi'$ unitary cuspidal automorphic representations of $\GL(m)/\Q$, $\GL(m')/\Q$, the Rankin-Selberg $L$-function $L(s,\pi \otimes \tilde{\pi})$ (see e.g. \cite[\S 1.1.2]{Mic07})  has a corresponding arithmetic conductor $\Cond(\pi \otimes \pi')$  and
analytic conductor,  given for $t \in \R$ by
\[
Q_{\pi\otimes\pi'}(t)=\Cond(\pi\otimes \pi')\prod_{j=1}^{mm'}(1+|it-\mu_{\pi\otimes\pi'}(j)|).
\]
The convexity bound for $L(s,\pi\otimes\pi')$ in the critical strip is known:
\[
L(s,\pi\otimes\pi') \ll_{ \pi, \pi', \ep} Q_{\pi\otimes\pi'}(t)^{\frac{1-\Re(s)}{2}+\varepsilon}, \qquad 0 \le \Re(s) \le 1.
\]

 \begin{remark}
Note that for each $1 \leq j \leq s$, the uniformity of the convexity bounds assumed in Condition \ref{KM_condition} \ref{item4} with respect to $m_j$ is critically reliant on the fact that within a family $Z_n^\Iscr(\Q,G;X)$, all fields share a fixed degree and a fixed Galois group of the Galois closure. 
 \end{remark}

\subsection{Verification of condition (\ref{item5_eqn}): controlling the propagation of bad $L$-function factors}
 Now we turn to the most difficult task: verifying that for each choice of  the family $Z_n^\Iscr(\Q,G;X)$ that we consider in our main theorems, condition (\ref{item5_eqn}) of Condition \ref{KM_condition_composed} is satisfied.

\subsubsection{Reframing the question in terms of subfields}

Let $Z_n^\Iscr(\Q,G;X)$ be a fixed family of fields, for a fixed transitive group $G \subseteq S_n$, and let $\rho$ be an irreducible representation of $G$. 
Let $L_1,L_2 \in \tilde{Z}_n^\Iscr(\Q,G;X)$. Then $\Gal(L_i/\Q) \simeq G$ while $\Gal(L_i^{\Ker(\rho)}/\Q) \simeq  G/\Ker(\rho)$.
 The following proposition transforms the property of identical $L$-functions into a property of identical fixed fields.
\begin{prop}\label{prop_KluNic_new}
Let $\rho$ be a fixed representation of  a fixed transitive subgroup $G\subseteq S_n$. For $L_1/\Q$ and  $L_2/\Q$ with $\Gal(L_1/\Q) \simeq \Gal(L_2/\Q) \simeq G$, then if
\beq\label{L_coincide}
L(s,\rho,L_1/\Q) = L(s,\rho,L_2/\Q)
\eeq
it follows that $L_1^{\Ker(\rho)} = L_2^{\Ker(\rho)}$.
\end{prop}

We recall a standard lemma. 
 \begin{lemma}
 Suppose  for two Galois extensions $F_1/\Q$ and $F_2/\Q$, that, aside from finitely many exceptions, the set of rational primes that  split completely in $F_1$ is the same as the set of rational primes that  split completely in $F_2$. Then $F_1=F_2$. 
 \end{lemma}
 \begin{proof}
By the Chebotarev density theorem, the density of rational primes that are split completely in $F_1$, $F_2$, or $F_1F_2$ are, respectively $[F_1:\Q]^{-1}$,$[F_2:\Q]^{-1}$,
$[F_1F_2:\Q]^{-1}$.  Since a prime is split completely in $F_1F_2$ if and only if it is split completely in $F_1$ and $F_2$, we have $[F_1:\Q]=[F_2:\Q]=[F_1F_2:\Q]$ and so $F_1=F_2$.
 \end{proof}

Thus to prove Proposition \ref{prop_KluNic_new}, it suffices to show that for each fixed representation $\rho$ of $G$, aside from finitely many exceptions, the set of rational primes that split completely in $L_1^{\Ker(\rho)}$ is the same as the set of rational primes that split completely in  $L_2^{\Ker(\rho)}$, under the assumption that $L(s,\rho,L_1/\Q) = L(s,\rho,L_2/\Q)$. First we assume that $p$ is a rational prime that is unramified in $L_1,L_2$ (and hence is unramified in $L_1^{\Ker(\rho)}, L_2^{\Ker(\rho)}$) and splits completely in $L_1^{\Ker(\rho)}$.  In particular, this means that for any $\pfrak_1$ in $L_1^{\Ker(\rho)}$ that lies above $p$,   the conjugacy class of the Frobenius element $\sig_{\pfrak_1}$ is trivial in $\Gal(L_1^{\Ker(\rho)}/\Q) \simeq G/\Ker(\rho)$, that is to say,  
$\rho(\sig_{\pfrak_1})$ is the identity matrix $I$. 
 
 Now letting $\pfrak_2 \in L_2^{\Ker(\rho)}$ be any prime lying  above $p$, by the assumption that the $L$-functions are equal, we have that the factors corresponding to $p$ are equal as functions of $s$ and therefore 
 \beq\label{factor_equal}
 \det (I - \rho (\sig_{\pfrak_2}) p^{-s})^{-1}=
   \det (I - \rho (\sig_{\pfrak_1}) p^{-s})^{-1}=
   \det (I - Ip^{-s})^{-1} .
   \eeq
   Now recall that the Frobenius element $\sig_{\pfrak_2}$ is necessarily finite order. 
We recall a simple observation. Suppose $M$ is an $n \times n$ matrix over $\C$ of finite order, say $M^k=I$ for some $k$, such that $\det(I -Mx)=\det(I-Ix)=(1-x)^n$ for a formal variable $x$. Then we claim $M=I$.  Indeed, since $M$ is finite order, $M$ is diagonalizable, for the minimal polynomial of $M$ divides $x^k-1$ and so has no
repeated roots.  By our second assumption, all the roots of the characteristic polynomial of $M$ are equal to 1, so that all the eigenvalues of $M$ are 1 and $M=I$. 
We apply this in (\ref{factor_equal}) to conclude that $\rho(\sig_{\pfrak_2})=I$ as well. 
Thus the conjugacy class of the Frobenius element $\sig_{\pfrak_2}$ is trivial in $\Gal(L_2^{\Ker(\rho)}/\Q) \simeq G/\Ker(\rho)$ and $p$ must split completely in $L_2^{\Ker(\rho)}$.

 In this fashion we see that any prime that is unramified in $L_1,L_2$ and splits completely in $L_1^{\Ker(\rho_1)}$ must split completely in $L_2^{\Ker(\rho)}$. Starting from primes   unramified in $L_1,L_2$ that split completely in $L_2^{\Ker(\rho)}$ we can similarly show that they must split completely in $L_1^{\Ker(\rho)}$, and this concludes the proof of Proposition \ref{prop_KluNic_new}.

 \begin{remark}\label{remark_other_k}
 Proposition \ref{prop_KluNic_new} can alternatively be deduced from  \cite[Theorem 5]{KluNic16}, which also includes a converse, which we do not require in our application. To apply  \cite[Theorem 5]{KluNic16} in our setting, one first passes as in \cite[p. 162]{KluNic16} to the case of a faithful representation $\varpi (\sig \Ker (\rho)) = \rho(\sig)$ acting on $H=G/\Ker(\rho)$. 
 Kl\"{u}ners and Nicolae present a counterexample to the characterization deduced in Proposition \ref{prop_KluNic_new} when working over $k \neq \Q$ \cite[p. 167]{KluNic16}, but see their relative version \cite[Thm. 6]{KluNic16}. It is possible that certain other families $Z_n^\Iscr(k,G;X)$ with $k \neq \Q$ and certain choices of $G$ can be treated by an adaptation of our methods with such a relative result. (When working over $k \neq \Q$ one would also need to take into account the more nuanced situation that arises with regards to arithmetic equivalence.)
 \end{remark}

We now apply Proposition \ref{prop_KluNic_new}.
As before, let $G$ be a fixed transitive subgroup of $S_n$ and let $\rho_1,\ldots, \rho_s$ be the nontrivial irreducible representations of $G$. For each $1 \leq j \leq s$, consider the set of fields
\[
\{ L^{\Ker(\rho_j)} : L \in \tilde{Z}_n^\Iscr(\Q,G;X) \}.
\]
(Note that we define this as a set, not a multi-set.) 
Philosophically, we would like to show that the cardinality of this set is ``large,'' or equivalently very few of the fields $L$ share the same fixed field, which would imply that ``few'' collisions $L(s,\rho_j, L_1/\Q) = L(s,\rho_j, L_2/\Q)$ could occur for $L_1 \neq L_2 \in \tilde{Z}_n^\Iscr(\Q,G;X). $

Formally, recall the definition of the set $\Lscr_j(X)$  in \S \ref{sec_families_aut} according to the family of fields $\mathscr{F}(\Q,G;X) =  \tilde{Z}_n^\Iscr(\Q,G;X)$.
Let us first consider the special case in which $\rho_j$ is faithful so that $\Ker(\rho_j)$ is trivial.
Then by Proposition \ref{prop_KluNic_new}, for two fields $L_1 \neq L_2 \in  \tilde{Z}_n^\Iscr(\Q,G;X)$, we cannot have
 $L_1^{\Ker(\rho_j)} = L_2^{\Ker(\rho_j)}$ and so we cannot have $L(s,\rho_j, L_1/\Q) = L(s,\rho_j, L_2/\Q)$, and so in this case 
\beq\label{Lj_L}
 |\mathscr{L}_j(X)| = |\tilde{Z}_n^\Iscr(\Q,G;X)| .
 \eeq
Thus if $\rho_j$ is faithful, we have verified  (\ref{item5_eqn}) of Condition \ref{KM_condition_composed} with $\tau_j=0$, which certainly suffices.

More generally, even if $\rho_j$ is not a faithful representation, Proposition \ref{prop_KluNic_new} shows that the number of fields $L_i \in \tilde{Z}_n^\Iscr(\Q,G;X)$ for which $L(s,\rho_j,L_i/\Q)$ is identical to a specific $L$-function is bounded above by the number of fields  $L_i \in \tilde{Z}_n^\Iscr(\Q,G;X)$ for which $L_i^{\Ker(\rho_j)}$ is identical to a specific field. Thus we have translated the problem of verifying (\ref{item5_eqn}) for a particular family $\Lscr_j(X)$ to a problem of counting fields. 

Precisely, we summarize the implications of Proposition \ref{prop_KluNic_new} as the following statement:
\begin{prop}\label{prop_count_equiv}
Let $Z_n^\Iscr(\Q,G)$ be a set of fields considered in Theorem \ref{thm_cheb_main_quant_pt1}, \ref{thm_cheb_main_quant_pt2}, \ref{thm_cheb_main_quant_pt3}, \ref{thm_cheb_main_quant_pt2'} or \ref{thm_cheb_main_quant_pt4} under the associated hypotheses, if any. Let $\tilde{Z}_n^\Iscr(\Q,G)$ be the corresponding set of Galois closures. Let $\rho_1,\ldots, \rho_s$ be the nontrivial irreducible representations of $G$. Define the families $(\Lscr_j(X))_{X \geq 1}$ for $1 \leq j \leq s$ accordingly, as in \S \ref{sec_families_aut}. Then $\tilde{Z}_n^\Iscr(\Q,G)$ and the families $(\Lscr_j(X))_{X \geq 1}$ for $1 \leq j \leq s$ satisfy (\ref{item5_eqn}) of Condition \ref{KM_condition_composed} with parameters $\{ \tau_j \}_{1 \leq j \leq s}$ if the following holds: 
 for each irreducible representation $\rho_j$ of $G$, given any  field $F \in Z_u(\Q,G/\Ker(\rho_j);X)$ (where $u=|G/\Ker(\rho_j)|$), at most $O_{n,G,\Iscr}(X^{\tau_j})$ fields $L \in  \tilde{Z}_n^\Iscr(\Q,G;X) $ have $L^{\Ker(\rho_j)} = F$.
\end{prop}

\subsubsection{Rationale for the restriction on ramification types of tamely ramified primes}\label{sec_rationale_ram}

For $G$ not a simple group, Proposition \ref{prop_count_equiv} spurs us to quantify, for each proper normal subgroup $H$ of $G$ that appears as the kernel of at least one (non-faithful, non-trivial) irreducible representation of $G$, how often a particular field occurs as a fixed field $L^H$, as $L$ varies over a relevant family of Galois extensions of $\Q$ with Galois group $G$. 

For certain groups $G$, fixed fields could collide with high repetition. 
For example, taking $G = \Z/4\Z$, then for any fixed quadratic field such as $F = \Q(e^{2\pi i /3})$, a positive proportion of quartic Galois fields $K \in Z_4(\Q,\Z/4\Z;X)$ have $K^{\Z/2\Z} = F$. This can be seen for example via a counting argument similar to that of \S \ref{sec_count_cyclic}. (See also comments in Remarks \ref{remark_G_abelian} and  \ref{remark_D4}.)

To eliminate such possibilities, we will critically use our restrictions on the ramification types of the tamely ramified primes in the fields in $Z_n^\Iscr(\Q,G;X)$.
Given $G$, we will select $\Iscr$ so that it has two properties: 
\begin{enumerate}
\item For $Z_n^{\Iscr}(\Q,G)$ to be infinite,  we need the elements in $\Iscr$ to generate $G$.

\item  We need $\Iscr$ to have the property that for each proper normal subgroup $H$ in $G$ that is the kernel of a non-faithful irreducible representation of $G$, given any field $K \in Z_n^\Iscr(\Q,G;X)$ with associated Galois closure $\tilde{K}/\Q$, then $p|D_K$ implies $p|D_F$, where $F=\tilde{K}^H$. 
\end{enumerate}
(Of course the primes that appear in $D_K$ are the same that appear in $D_{\tilde{K}}$ but this need not \emph{a priori} be true of $D_{\tilde{K}}$ and $D_F$.) 
Property (2) will enable us to obtain the information we seek in Proposition \ref{prop_count_equiv}, that is, to count the number of $\tilde{K} \in \tilde{Z}_n^\Iscr(\Q,G;X)$ sharing the same fixed field $F=\tilde{K}^H$, by applying quantitative information about $\Dbf_n(G;\varpi)$ (Property \ref{property_D}).
This is one of the most novel features of this paper.

\subsection{The counting problem}\label{sec_counting_problem}

We now define the counting problem that is the heart of the matter. 

\begin{property}[Property $\mathrm{Mult}_n(G,\Iscr;\tau)$]
Let  $\mathrm{Mult}_n(G,\Iscr;\tau)$ denote the property that for every $X \geq 1$, for each  irreducible representation $\rho$ of $G$, given any particular field $F \in Z_u(\Q,G/\Ker(\rho))$ (with $u=|G/\Ker(\rho)|$) that arises as a fixed field $\tilde{K}^{\Ker(\rho)}$ for at least one field $K \in Z_n^\Iscr(\Q,G;X)$, for every $\ep>0$, at most $O_{n,G,\Iscr,\ep}(X^{\tau+\ep})$ fields $K \in Z _n^\Iscr(\Q,G;X) $ have $\tilde{K}^{\Ker(\rho)} = F$.
\end{property}

Given a family $Z_n^\Iscr(\Q,G;X)$, if we can prove that $\mathrm{Mult}_n(G,\Iscr;\tau)$ holds for a sufficiently small $\tau$, then 
by Proposition \ref{prop_count_equiv}, the relevant effective Chebotarev density theorem for the family  $Z_n^\Iscr(\Q,G;X)$ will follow (that is, either Theorem \ref{thm_cheb_main_quant_pt1}, \ref{thm_cheb_main_quant_pt2}, \ref{thm_cheb_main_quant_pt3}, \ref{thm_cheb_main_quant_pt2'} or \ref{thm_cheb_main_quant_pt4}). Quantitatively, lowering the size of $\tau$ for which we can prove $\mathrm{Mult}_n(G,\Iscr;\tau)$ will allow us to better control the size of the possible exceptional set of fields.

 \begin{prop}[Counting problem]\label{ques_count}

We can prove the following:
\begin{enumerate}
\item $\mathrm{Mult}_n(G,\Iscr;0)$ for $G$ a simple group, $\Iscr$ imposing no restriction.
\item $\mathrm{Mult}_n(G,\Iscr;0)$  for $G$ cyclic, $\Iscr$ specifying totally ramified.
\item $\mathrm{Mult}_n(S_n,\Iscr;\varpi_n)$ for $n  \geq 3$, $\Iscr$ the conjugacy class $[(1 \; 2)]$, where $\varpi_3 = 1/3$, $\varpi_4 = 1/2$, $\varpi_5=199/200$, and for $n \geq 6,$ $\varpi_n = \varpi$ if we assume Property $\Dbf_n(S_n,\varpi)$.
\item $\mathrm{Mult}_p(D_p,\Iscr; \tau_p)$ holds for $\tau_p=1/(p-1)$, $p$ an odd prime, $\Iscr$ the conjugacy class of order 2 elements.
\item $\mathrm{Mult}_4(A_4,\Iscr;0.2784...)$, $\Iscr$ the two conjugacy classes of order 3 elements.
\end{enumerate}
 \end{prop}
 As observed above, $\mathrm{Mult}_n(G,\Iscr;0)$ is tautologically true when $G$ is a simple group ($\Iscr$ imposing no restriction), since all the nontrivial irreducible representations are faithful and (\ref{Lj_L}) applies. All the other cases of the counting problem require work.
 We first explicitly prove this for $S_n$, $n \neq 4$; in particular, to aid the reader, we include our argumentation for choosing $\Iscr = [(1 \; 2)]$.

\subsubsection{Background lemmas on inertia groups and discriminants}\label{sec_prelim_disc_lemma}

We recall standard results on the powers of primes dividing $D_K$.

\begin{lemma}[Powers of tamely ramified primes in discriminants]\label{lemma_power_p}
Let $K \subset \tilde{K} \subset \overline{\Q}$ with $\Gal(\tilde{K}/\Q) \simeq G$ and $H=\Gal(\tilde{K}/K)$.
Let $p$ be a rational prime that is tamely ramified in $K$ and $\tilde{K}$, and has an inertia group in $\Gal(\tilde{K}/\Q)$ generated by $\pi \in G$. The power $\al$ such that $p^\al || D_K$ is  
\beq\label{order_field}
 [G:H] - \text{number of orbits of $\pi$ acting on the cosets $G/H$}.
 \eeq
 
\end{lemma}
\begin{proof}
We have that $D_K$ is the Artin conductor of $\tilde{K}/\Q$ for the permutation representation  of $G$ on $G/H$ \cite[Ch. VII, Corollary 11.8]{Neu}.
By definition, the Artin conductor of $\tilde{K}/\Q$ for a representation $V$ of $\Gal(\tilde{K}/\Q)$ is $\prod_p p^{f_p(V)}$, where the product is over rational primes and 
$$
f_p(V)=\sum_{i\geq 0} \frac{g_{p,i}}{g_{p,0}}\operatorname{codim} V^{G_{p,i}},
$$
for $G_{p,i}$ an $i$th ramification group for $p$ in $\Gal(\tilde{K}/\Q)$ and $g_{p,i}:=|G_{p,i}|$.  Recall that $G_{p,0}$ is the inertia group $I_p$ and
that for tamely ramified $p$, we have $G_{p,i}=1$ for $i\geq 1$. So for tamely ramified $p$, we have
$
f_p(V)=\operatorname{codim} V^{I_p}.
$
The lemma follows, since for a permutation representation $V$, the dimension of the fixed subspace $V^{\pi}$ is the number of orbits of $\pi$.
\end{proof}

\begin{lemma}[Maximum contribution of wild primes]\label{lemma_wild_primes}
Let $G$ be a transitive subgroup of $S_n$. Then for all fields $K \in Z_n(\Q,G)$,  the total contribution to $D_K$ from the rational primes that are wildly ramified in $K$ is at most  a certain finite constant $C_{G}$ depending only on $G$.
\end{lemma}
This lemma follows from \cite[Ch. III, Thm. 2.6]{Neu} and the fact that all wildly ramified  primes divide $|G|$.

In order to consider only the tame part of the discriminant in our investigations below, it will be convenient to use the following notation. 
Given a finite set of primes $\Omega$,  define $D_K^{(\Omega)}$ to denote the 
contribution to the discriminant from primes $p \not\in \Omega$, i.e. $D_K^{(\Omega)}$ is the maximal positive divisor of $D_K$ that is not divisible by any prime in $\Omega$. We will apply this in particular when $\Omega$ is comprised of the primes dividing $|G|$.

\subsubsection{Exemplar case: $G = S_n$, $n=3$ or $n \geq 5$}\label{sec_countSn}

Recall that when $n=3$ or $n\ge 5$, $S_n$ has one nontrivial, proper normal subgroup, namely $A_n$, which certainly appears as the kernel of the sign representation. Thus we must specify a ramification type $\Iscr$ so that the counting problem for fixed fields $\tilde{K}^{A_n}$ can be handled.
We wish, for a fixed quadratic field $F \in Z_2(\Q,C_2)$, to count the number of degree $n$ fields $K \in Z_n(\Q,S_n;X)$ such that $\tilde{K}^{A_n}=F$.
\begin{center}
\begin{tikzpicture}[-,>=stealth',shorten >=1pt,auto,node distance=2cm, thick,main node/.style]
\node[main node] (Kt) {$\widetilde K$};
\node[main node] (K) [below of=Kt] {$K$};
\node[main node] (F) [right  of=K] {$F = \widetilde K^{A_n}$};
\node[main node] (Q) [below of=K] {$\Q$};

\path
(Kt)	edge [bend right] node [left]{$\Gal(\widetilde{K}/\Q)\simeq S_n$} (Q)
      	edge node {} (K)
        edge node {} (F)
(K)	edge node {$n$} (Q)
(F) 	edge [bend right]  node [right=1mm]{$\Gal(\widetilde{K}/ F)\simeq A_n$} (Kt)
        	edge node{$2$} (Q);
\end{tikzpicture}
\end{center}

Using Lemma \ref{lemma_power_p}, we can compute for fields $K,\tilde{K},F$ in such a constellation the exact power of $p$ that appears in the absolute discriminants $D_K, D_{\tilde{K}}, D_F$, for each prime $p \ndiv |G|$. 
We show these exponents in Table \ref{table_Sn}: the leftmost column specifies the conjugacy class of the generating element $\pi$ of the (cyclic) inertia group for $p$, while the other columns specify the exact power of $p$ appearing in the discriminants.  We only list a few of the $p(n)$ conjugacy classes of $S_n$;  we set $\ep_n=0$ if $n$ is odd and $\ep_n=1$ if $n$ is  even.

{\renewcommand{\arraystretch}{1.5}
\begin{table}[H]
\begin{center}
\begin{tabular}{|C{10em} |C{3cm} C{3cm}  C{3cm}|}
 \hline
&  \multicolumn{3}{c|}{Exponent of $p$ appearing in the discriminant of} \\
 Inertia type of $p$ & $K$  & $ \widetilde{K}$ &$F = \widetilde{K}^{A_n}$\\

 \hline
$[()]$ 	&  0 	&  0	&  0	\\
 $[(1 \; 2)]$	&  $1$	&  $ n! - n!/2$	  & 1	\\
  $[(1\; 2\; 3)] = [(1\;2)(2\;3)]$	&  2  	&  $ n! - n!/3$	  & 0	\\
    $[(1\;2)(3\;4)] $	&  2   	&  	 $ n! - n!/2$ & 	0\\
$\vdots$				&  $\vdots$  	& $\vdots$	& $\vdots$	\\
$[(1\;2\;3\ldots n)]		$		& $ n-1  $	& $ n! - n!/n$	& $\ep_n$	\\
 \hline
\end{tabular}\vspace{.15in}
\caption{Table of exponents for $p$ when $\Gal(\widetilde{K}/\Q) \simeq S_n$, for each $p \ndiv n!$}\label{table_Sn}
\end{center}
\end{table}
}

From Table \ref{table_Sn} we observe that every $p \ndiv |G|$ that  has inertia group generated by a transposition has $p\|D_K, p^{n!/2} \| D_{\tilde{K}}, p \|D_F$. This will allow us to control, for a fixed field $F$,  how many $K$ can yield  a constellation including $F$. 
These observations from Table \ref{table_Sn} motivated our choice of $\Iscr=[(1 \; 2)]$ for $G \simeq S_n$ ($n=3$, $n\geq 5$).

Now we come to the crux of the argument. Suppose that $F$ is fixed, and hence $D_F \geq 1$ is fixed. 
Set $\Omega = \{ p : p|n!\}$ and recall the notation $D_K^{(\Omega)}$ defined above.  Our discussion above shows that any degree $n$ extension $K \in Z_n^\Iscr(\Q,S_n;X)$ such that $\tilde{K}^{A_n} = F$ must have
\beq\label{DKDF}
 D_K^{(\Omega)} = D_F^{(\Omega)}.
\eeq
Assuming Property $\Dbf_n(S_n,\varpi)$ is known, then since the power of any $p\in \Omega$ dividing $D_K$ is bounded in terms of $n$, for a given $F$ there are at most $\ll_{n,G,\ep} D_F^{\varpi+\ep} \ll_{n,G,\ep} X^{\varpi+\ep}$ such $K$ satisfying (\ref{DKDF}), for every $\ep>0$. 
Now to obtain the conclusion on $\mathrm{Mult}_n(S_n,\Iscr;\varpi_n)$ of Proposition \ref{ques_count} for $S_n$, $n =3, n \geq 5$, we simply apply the currently best known upper bounds for Property $\Dbf_n(S_n,\varpi)$ in these cases, as stated in \S \ref{sec_sym_gp}.

Having completed this exemplar case of $G=S_n$ ($n \neq 4$) in some detail, we are now more brief with the remaining groups $G$, which follow similarly by using  Lemma \ref{lemma_power_p} in order to fix an appropriate choice of ramification type $\Iscr$ for which the counting problem can be resolved.

\subsubsection{$G \simeq S_4$} Recall that $S_4$ has four nontrivial irreducible representations (see e.g. \cite[p. 43]{Ser77}): two three-dimensional faithful representations (the standard representation and the product of the standard representation with the sign representation) and two non-faithful representations. The subgroup $A_4$ is the kernel of the one-dimensional sign representation, and $K_4 \simeq C_2 \times C_2$ the Klein four group is the kernel of the irreducible two-dimensional representation of $S_4$. We thus have two counting problems to consider.

Relevant to the counting problem for fixed fields under $A_4$, by choosing $\Iscr$ to be the conjugacy class $[(1 \; 2)]$ of transpositions, 
we may conclude that triads $K, \tilde{K}, F=\tilde{K}^{A_4}$ behave exactly as in the case of $S_n$ in \S \ref{sec_countSn}, so that for any $p \ndiv 4!$, 
$p\|D_K, p^{12} \| D_{\tilde{K}}, p \|D_F$, and hence upon setting $\Omega = \{2,3\}$, we have 
$D_K^{(\Omega)} = D_F^{(\Omega)}.$
 Thus arguing as in \S \ref{sec_countSn}, any fixed $F$ corresponds to at most $\ll_\ep D_F^{\varpi +\ep} \ll_{n,G,\ep} X^{\varpi + \ep}$ possibilities for $K \in Z_4^\Iscr(\Q,S_4;X)$, if Property $\Dbf_4(S_4,\varpi)$ is known. 

Relevant to the counting problem for fixed fields under $K_4$, still choosing $\Iscr$ to be the conjugacy class $[(1 \; 2)]$ of transpositions, for triads $K, \tilde{K}, F=\tilde{K}^{K_4}$ the exponents are different: for every $p \ndiv 4! $ we have $p \| D_K, p^{12}\|D_{\tilde{K}}$, and $p^{3}\|D_F$. Thus upon setting $\Omega = \{2,3\}$,  we have 
\beq\label{DKDF_S4}
D_K^{(\Omega)} = (D_F^{(\Omega)})^{1/3},
\eeq
 and so  
 any fixed $F$ corresponds to at most $\ll_\ep D_F^{\varpi/3 +\ep} \ll_{n,G,\ep} X^{\varpi+\ep}$ possibilities (for every $\ep>0$) for $K \in Z_4^\Iscr(\Q,S_4;X)$, if Property $\Dbf_4(S_4,\varpi)$ is known. 
(Here we have used the fact that if $K \in Z_n(\Q,G;X)$ and (\ref{DKDF_S4}) holds, then $ D_F \ll_{n,G} X^3$.)
We conclude that $\mathrm{Mult}_4(S_4,\Iscr;1/2)$ holds since Property $\Dbf_4(S_4,1/2)$ is known.

\subsubsection{$G \simeq A_4$}

Recall (see \cite[Section 5.7, page 41]{Ser77}) that $A_4$ has four nontrivial irreducible representations: two faithful representations and two one-dimensional non-faithful representations,  each with kernel $K_4 \simeq C_2 \times C_2$ the Klein-four group. Thus we need only complete the counting problem for triads $K, \tilde{K}, F=\tilde{K}^{K_4}$.
We will require all tamely ramified primes to have inertia type belonging to either of the conjugacy classes $\Cscr_1,\Cscr_2$ of order 3 elements (specified e.g. in Proposition \ref{prop_counting_A4_fields}).

Suppose we restrict to primes of inertia type in the conjugacy class $\mathscr{C}_1$. The image of this inertia type in $A_4/K_4 \simeq C_3$ is nontrivial, and we see that for any $p\ndiv |A_4|$, $p^2 \| D_K, p^8 \| \tilde{K}, p^2 \| D_F$.  Thus upon setting $\Omega = \{2,3\}$, within the triad we have 
$ D_K^{(\Omega)} = D_F^{(\Omega)},
$
 and so,  
  any fixed $F$ corresponds to at most $\ll_\ep D_F^{\varpi +\ep} \ll_{n,G,\ep}X^{\varpi+\ep}$ possibilities for $K \in Z_4^\Iscr(\Q,A_4;X)$, if Property $\Dbf_4(A_4,\varpi)$ is known. The computation for primes of inertia type in the conjugacy class $\mathscr{C}_2$ is identical. Recalling our result of Proposition \ref{prop_counting_A4_fields}, we conclude that $\mathrm{Mult}_4(S_4,\Iscr;0.2784...)$ holds.

\subsubsection{$G \simeq D_p$, $p$ an odd prime}
We think of $D_p$ (with $p$ an odd prime) as the group of order $2p$ of symmetries on a regular $p$-gon, acting in the usual way. Thus $D_p$  has one nontrivial, proper normal subgroup, namely $C_p$; this subgroup certainly appears as the kernel of the (one-dimensional) sign representation. Thus we must consider the corresponding counting problem for fixed fields $\tilde{K}^{C_p}$. We restrict the inertia type $\Iscr$ to the conjugacy class $[(2 \; p)(3 \; \; (p-1))\cdots(\frac{p+1}{2} \; \frac{p+3}{2})]$, that is the conjugacy class of reflections (each with with $(p+1)/2$ orbits acting on $p$ elements). 
 
For a triad $K, \tilde{K}, F= \tilde{K}^{C_p}$ we then have for every prime $\ell \ndiv 2p$ that $\ell^{(p-1)/2}  \| D_K, \ell^p \| D_{\tilde{K}},  \ell \| D_F$. 
Thus upon setting $\Omega =\{2,p\}$ we have 
\beq\label{DK_DF_Dp}
D_K^{(\Omega)} = (D_F^{(\Omega)})^{\frac{p-1}{2}},
\eeq
 and so  
 any fixed $F$ corresponds to at most $\ll_{p,D_p,\ep} D_F^{(p-1)\varpi/2 +\ep}$ possibilities for $K \in Z_p^\Iscr(\Q,D_p;X)$, if Property $\Dbf_p^\Iscr(D_p,\varpi)$ is known.
 Now if (\ref{DK_DF_Dp}) is known and $K \in Z_p^\Iscr(\Q,D_p;X)$ then $D_F \ll_{p,D_p} X^{2/(p-1)}$, so we have at most
 $\ll_{p,D_p,\ep} X^\varpi$ choices for such $K$ if Property $\Dbf_p^\Iscr(D_p,\varpi)$ is known.
 We conclude from Proposition \ref{prop_counting_Dp_fields} that $\mathrm{Mult}_p(D_p,\Iscr; 1/(p-1))$ holds unconditionally.

\subsubsection{$G$ a cyclic group}
Finally, for $G$ a cyclic group of order $n$, note that $Z_n(\Q,G;X)$ already is comprised of Galois fields, so we do not need to pass to the Galois closures. (As a special case, if $G \simeq C_p$ with $p$ prime, then $G$ has no nontrivial proper (normal) subgroups, so all nontrivial representations are faithful, without the need to artificially impose a ramification restriction. But  in this case, every ramified prime is naturally totally ramified, so we still group this with the general case below.)
In general, consider $G$ an arbitrary cyclic group of order $n$, say $G \simeq C_{p_1^{e_1}} \times  \cdots \times C_{p_k^{e_k}}$ with distinct primes $p_1,\ldots, p_k$. 
We restrict to $\Iscr$ specifying that every tamely ramified prime must be totally ramified, that is, its inertia group must be generated by an element of full order in $G$; in particular, such an element does not belong to any proper, nontrivial subgroup $C_m$ of $C_n$.
By Lemma \ref{lemma_power_p} the following properties hold:
 \begin{enumerate}
 \item for every prime $\ell \ndiv n$ we have $\ell^{n-1} \| D_K = D_{\tilde{K}}$;
 \item for every nontrivial proper (normal) subgroup $C_m$ of $C_n$ (corresponding to a proper divisor $m|n$) there exists an integer $1 \leq \al_m \leq n-1$ (depending on $m$ and $C_n$) such that $\ell^{\al_m} \| D_F$ where $F = \tilde{K}^{C_m}$. 
 \end{enumerate}
 As a result, upon setting $\Omega = \{ p : p |n\}$,  for each nontrivial proper subgroup $C_m$ of $G$, parametrized by divisors $m$, we have that
 \[D_K^{(\Omega)} = (D_F^{(\Omega)})^{\frac{n-1}{\al_m}}\]
 when $F = K^{C_m} = \tilde{K}^{C_m}$.
  Thus %by Lemma \ref{lemma_propertyDn_tame}, 
  any fixed $F$ corresponds to at most $\ll_{n,m,C_n,\ep} D_F^{\varpi (n-1)/\al_m+\ep} \ll_{n,m,C_n,\ep} X^{\varpi+\ep}$ possibilities (for any $\ep>0$) for $K \in Z_n^\Iscr(\Q,C_n;X)$ if Property $\Dbf_n(C_n,\varpi)$ is known. By  Proposition \ref{prop_counting_cyclic_fields} we have $\Dbf_n(C_n,0)$, so that we have verified $\mathrm{Mult}_n(C_n,\Iscr;0)$.
 
\begin{remark}[Non-cyclic abelian groups]\label{remark_G_abelian}
The above arguments show that we are able to pick an appropriate ramification restriction to control the propagation of bad $L$-function factors  if there exists a set that generates $G$ and such that none of them lies in any (nontrivial, proper, normal) subgroup $H$ of $G$ that appears as the kernel of at least one nontrivial irreducible representation of $G$. 
We may already observe the difficulty of adapting this general strategy to a non-cyclic abelian group by considering the simple case of $G\simeq C_{p^{e}} \times C_{p^{f}}$ for a prime $p$.  
 Consider an element in the generating set of  the form $(a,b)$ with $a\ne 0$. Let $p^k$ be the highest power of $p$ that divides both $a$ and $b$. Then for $\zeta_p = e^{2\pi i/p}$, the map $C_{p^e} \times C_{p^f} \ra \C$ given by $(i,j)\mapsto \zeta_p^{ib/p^k  -j a /p^k}$  is a non-trivial irreducible representation of $C_p^e \times C_p^f$, and our generator is in the kernel of this map.
\end{remark}

\begin{remark}[Quartic $D_4$-fields]\label{remark_D4}
Difficulties also arise for quartic $D_4$-fields: there are irreducible representations of $D_4$ with kernels $K_4$, $K_4'$ (two different subgroups isomorphic to the Klein-four group) and $C_4$, but no set of generators of $D_4$ that avoid all three of these subgroups, and hence no choice of ramification type $\Iscr$ for which the three counting problems can simultaneously be resolved. It may be possible to apply our method to a particular subfamily of quartic $D_4$-fields generated from a fixed biquadratic field; in this case the counting problems will be trivial, although proving a lower bound that grows with $X$ for such a family may not be.
\end{remark}

\subsection{Deduction of Theorems \ref{thm_cheb_main_quant_pt1}, \ref{thm_cheb_main_quant_pt2}, \ref{thm_cheb_main_quant_pt3}, \ref{thm_cheb_main_quant_pt2'} and \ref{thm_cheb_main_quant_pt4}  from Theorem \ref{thm_KM_badfamily}}
We have verified Condition \ref{KM_condition_composed} for each family $Z_n^\Iscr(\Q,G;X)$ considered in the above theorems; now we apply
Theorem \ref{thm_KM_badfamily}.
 The family parameters notated in Condition \ref{KM_condition_composed}, namely $\{M_0,M_1,A,d\}$  and $\{ m_j,A_j,d_j,M_{1,j},M_{2,j,\ep},M_{3,j,\ep}, M_{5,j}\}$ for $1 \leq j \leq s$, all depend only on $n,G,\Iscr$, and thus in the following statements we can replace any dependence on family parameters by dependence on $n,G,\Iscr$.

\begin{prop}\label{thm_KM_familyZn_prop}
Fix a family $Z_n^\Iscr(\Q,G;X)$ considered in Theorem \ref{thm_cheb_main_quant_pt1}, \ref{thm_cheb_main_quant_pt2}, \ref{thm_cheb_main_quant_pt3}, \ref{thm_cheb_main_quant_pt2'} or \ref{thm_cheb_main_quant_pt4} under the associated hypotheses (if any).
If it is known that $X^\be \ll_{n,G,\Iscr} |Z_n^\Iscr(\Q,G;X)| \ll_{n,G,\Iscr} X^d$
and
that $\mathrm{Mult}_n(G,\Iscr;\tau_*)$ holds for some $\tau_* < \be$, then the conclusions of the relevant theorem hold for those values of $\tau_*, \be, d$. 
\end{prop}

Let $\tau_*< \be \leq d$ be as assumed in the proposition. Fix $\tau = \tau_* + \ep_3$ for some sufficiently small $\ep_3$ (in particular so that $\tau < d$) and fix $\ep_0 \leq \min \{1/2,2(d-\tau)\}$ sufficiently small. We apply Theorem \ref{thm_KM_badfamily} with
\beq\label{Del_tau_d}
\Del = 1- \frac{\tau}{d} - \frac{\ep_0}{2d},
\eeq
$\del$ chosen as in (\ref{describe_del}) (according to $A=|G|/2$ and $\ep_2 = \ep_0$ so that we obtain the expression  for $\del$ in Remark \ref{remark_KM_intro}),
and $\eta = \ep_0/2d$. Then 
\[  (1 - (1-\eta)\Del)d =  \tau + \frac{\ep_0}{2}  + \frac{\ep_0}{2} (1 - \tau/d-\ep_0/2d) 
	 \leq \tau + \ep_0.\]
 Then there exists $B$ depending only on $n,G,\Iscr$ such that for all $X \geq 1$, at most 
\beq\label{exceptional1}
O_{n,G,\Iscr,\tau,d,\ep_0}(X^{ \tau + \ep_0})
\eeq
 fields $K \in Z_n^\Iscr(\Q,G;X)$ are such that $\zeta_{\tilde{K}}/\zeta$ can have a zero in the region 
\beq\label{zerofree_del_ep}
[1-\del,1] \times [-X^\be,X^\be],
\eeq
where $\be = \ep_0(1 - \tau/d - \ep_0/2d)/(2B).$

 Our goal now is to express this in terms of how many $\del$-exceptional fields there can be.
It is temporarily convenient to work in terms of families of fields with discriminant in a dyadic range; thus we set $Z_n^{\Iscr,\sharp}(\Q,G;X)$ to be the subset of $Z_n^\Iscr(\Q,G;X)$ with $X/2<D_K \leq X$.  
We next verify that for $X$ sufficiently large, for every $K \in Z_n^{\Iscr,\#}(\Q,G;X)$ the region (\ref{zerofree_del_ep}) contains the region (\ref{region1_intro}), which we write now in the notation
\beq\label{region1_endgame}
[1-\del,1] \times [ - (\log D_{\tilde{K}})^{2/\del}, (\log D_{\tilde{K}})^{2/\del}].
\eeq
If $K \in Z_n^{\Iscr,\sharp}(\Q,G;X)$ then by Lemma \ref{lemma_compare_disc}, $C_1(n,G)(X/2)^{|G|/n} \leq D_{\tilde{K}} \leq C_2(n,G) X^{|G|/2}$, for certain constants $C_i(n,G)$. Thus it suffices to show that  there exists a threshold $D_3=D_3(n, G, \Iscr,\tau,d,\del,\ep_0)$ such that if $X \geq D_3$ then
\beq\label{D_lowerbound2}
 (\log (C_2(n,G)X^{|G|/2}))^{2/\del} \leq X^{\be}.
 \eeq
 This is the claim that a fixed power of $X$ is larger than any fixed power of $\log X$, as long as $X$ is sufficiently large; thus an appropriate threshold $D_3$ exists. 
 
 We have shown that for every $X \geq 1$ there are at most 
 $O_{n,G,\Iscr,\tau,d,\ep_0}(X^{ \tau + \ep_0})$ fields  $K \in Z_n^\Iscr(\Q,G;X)$ such  that $\zeta_{\tilde{K}}/\zeta$ can have a zero in (\ref{zerofree_del_ep}); consequently if $X /2\geq D_3$, at most 
 $O_{n,G,\Iscr,\tau,d,\ep_0}(X^{ \tau + \ep_0})$ fields in $K \in Z_n^{\Iscr,\#}(\Q,G;X)$ are such  that $\zeta_{\tilde{K}}/\zeta$ can have a zero in 
(\ref{region1_endgame}), that is, can be $\del$-exceptional.
 Now we suppose that $A \geq 2$ has been fixed, and we recall the threshold 
 $D_0$
  from Theorem \ref{thm_assume_zerofree1}.
As long as
\beq\label{D_lowerbound2a}
X/2 \geq  D_0,
\eeq
any $K \in Z_n^{\Iscr,\#}(\Q,G;X)$ that is not $\del$-exceptional satisfies the hypothesis of Theorem \ref{thm_assume_zerofree1}, and therefore 
for every conjugacy class $\mathscr{C} \subseteq G$ yields (\ref{asymp_log_intro}) for all $x$ sufficiently large as in (\ref{x_lower_bd2_intro}).
Upon taking
\[ D_4=D_4(n,G,\Iscr,\tau,d,\del,\ep_0,c_\Q,C_1,C_2,A)  := \max \{ D_0, D_3\}\]
we have shown that for any $X$ such that $X/2 \geq D_4$ we have that at most $O_{n,G,\Iscr,\tau,d,\ep_0}(X^{ \tau + \ep_0})$ fields in $K \in Z_n^{\Iscr,\#}(\Q,G;X)$ can be $\del$-exceptional, and for all remaining fields, (\ref{asymp_log_intro}) holds for all $x$ satisfying (\ref{x_lower_bd2_intro}).
We may in fact omit the dependence on $\del$ in the notation, as it is defined in terms of the other parameters.

The final step to complete the proof of Proposition \ref{thm_KM_familyZn_prop} is to sum over dyadic ranges of discriminants.
Now for any $X \geq 1$ (say using $\log_2$ temporarily), 
\[ Z_n^\Iscr (\Q,G;X) \subseteq \Union_{j=0}^{1+\log X} Z_n^{\Iscr,\sharp}(\Q,G;2^j) .\]
We may dissect this into two pieces: those for which $j$ is such that $2^{j-1} \geq D_4$, in which case our work above applies, and we conclude that the number of $\del$-exceptional  fields in 
\[  \Union_{2^{j-1} \geq D_4}^{1+\log X} Z_n^{\Iscr,\sharp}(\Q,G;2^j) \]
is at most 
\beq\label{exceptional_quant}
\sum_{2^{j-1} \geq D_4}^{1+\log X} O_{n,G,\Iscr,\tau,d,\ep_0}((2^j)^{ \tau + \ep_0})
	= O_{n,G,\Iscr,\tau,d,\ep_0}(X^{ \tau + \ep_0}).
\eeq
For those $j$ such that $2^{j-1} \leq D_4$, we count all the fields as possible exceptions, noting that 
\[
 \left| \Union_{1 \leq 2^{j-1} \leq D_4} Z_n^{\Iscr,\sharp}(\Q,G;2^j) \right|
	 	 \leq |Z_n^\Iscr (\Q,G;2D_4)| \\
	 \ll_{n,G} D_4^{d}.\]
We enlarge the implied constant in (\ref{exceptional_quant}) to include this constant, and call the resulting implied constant $D_5$, as appears in the theorem statements.  
This completes the proof of Proposition \ref{thm_KM_familyZn_prop}, and in combination with the values of $\tau_*$ supplied by Proposition \ref{ques_count},
we have proved Theorems \ref{thm_cheb_main_quant_pt1}, \ref{thm_cheb_main_quant_pt2}, \ref{thm_cheb_main_quant_pt3}, \ref{thm_cheb_main_quant_pt2'} and \ref{thm_cheb_main_quant_pt4} (and the non-quantitative Theorem \ref{thm_gen}).

\subsection{Proof of Corollary \ref{cor_primes_families}}

Let $Z_n^\Iscr(\Q,G;X)$ be a specified family, with corresponding parameters $\tau_*< \be \leq d$, set $A=2$ and let $\ep_0$ (sufficiently small) be fixed, with corresponding choice $\del \leq 1/4$. 
First, we verify that for $\sig>0$ fixed, there is a threshold $D_6'=D_6'(n,G,\Iscr,d,c_\Q,C_5,C_6,\ep_0,\sig)$ such that for $D_K \geq D_6'$, 
 \[ D_K^\sig \geq \ka_1 \exp \{ \ka_2 (\log \log (D_{\tilde{K}}^{\ka_3}))^{5/3}(\log\log \log (D_{\tilde{K}}^2))^{1/3}\} ,\]
where this lower bound is as stated in (\ref{x_lower_bd2_intro}), and the parameters $\ka_i$ have the dependencies
$ \ka_i =  \ka_i(n,G,d,c_\Q,C_5,C_6,\ep_0)$
(dropping the notational dependence on $A=2$).
In fact it suffices to compute a threshold above which
 \[ D_K^\sig \geq \ka_1 \exp \{ \ka_2 ( \log \log (D_{\tilde{K}}^{\ka_5}))^2\}  \]
where we set $\ka_5 = \max \{\ka_3,2\}$.
By Lemma \ref{lemma_compare_disc}, $D_{\tilde{K}} \leq C_2(n,G) D_K^{|G|/2}$ for a certain constant $C_2(n,G)$, so that it further suffices to show
 \[ D_K^\sig \geq \ka_1 \exp \{ \ka_2 ( \log \log (\ka_6 D_K^{\ka_7}))^2\}  \]
 where $\ka_6 = C_2(n,G)^{\ka_5}$ and $\ka_7= \ka_5 |G|/2$.
This will hold when $D_K$ is sufficiently large that
\[ \sig \geq \frac{ \log \ka_1 }{ \log D_K} +\frac{ \ka_2 (\log \log (\ka_6 D_K^{\ka_7}))^2}{ \log D_K} , \]
and we denote this threshold by $D_6'=D_6'(n,G,\Iscr,d,c_\Q,C_5,C_6,\ep_0,\sig)$.
Finally, recall the parameter $D_0$ provided in Theorem \ref{thm_assume_zerofree1}. While this is used as a constraint $D_{\tilde{K}} \geq D_0$, we apply Lemma \ref{lemma_compare_disc} to see that $D_{\tilde{K}} \geq C_1(n,G) D_K^{|G|/n}$ for a certain constant $C_1(n,G)$. Then $D_{\tilde{K}} \geq D_0$ is certainly satisfied if $D_K \geq D_0'$ with 
\beq\label{D0'_dfn}
D_0' :=  (C_1(n,G)^{-1} D_0)^{n/|G|}.
\eeq

Now for part (1) of Corollary \ref{cor_primes_families}, we may conclude from Theorem \ref{thm_assume_zerofree1} with $A=2$
that for every $X \geq 1$, for every field in $Z_n^\Iscr(\Q,G;X)$ that has $D_K \geq \max \{ D_0', D_6'\}$ and is not $\del$-exceptional, 
\beq\label{asymp_DK_sig}
\left| \pi_\Cscr (D_K^\sig,\tilde{K}/\Q)  - \frac{|\Cscr|}{|G|} \Li(D_K^{\sig}) \right| \leq \frac{|\Cscr|}{|G|} \frac{D_K^\sig}{(\log D_K^\sig)^2}.\eeq
Finally, we enlarge $\max \{ D_0', D_6'\}$ if necessary  to a parameter $D_6$, so that for all $D_K \geq D_6$, 
 the error term in (\ref{asymp_DK_sig}) is  at most 
$(1/2)|G|^{-1} \Li (D_K^\sig) \leq  (1/2)|\Cscr||G|^{-1} \Li (D_K^\sig).$
Then $\pi_\Cscr (D_K^\sig,\tilde{K}/\Q) \geq (1/2)|\Cscr||G|^{-1} \Li (D_K^\sig) \geq (1/2)|G|^{-1} \Li (D_K^\sig) $, and we can further enlarge $D_6$ if necessary to write the lower bound as in (\ref{many_primes}).
 
For part (2) of Corollary \ref{cor_primes_families}, we may follow e.g. Vaaler and Widmer \cite[Lemma 5.1]{VaaWid13} (but without assuming GRH, as they do).  
Suppose that $K$ is not $\del$-exceptional and furthermore that $D_K\geq D_0'$ with parameter $D_0'$ as above in (\ref{D0'_dfn}). Then for every $x$ satisfying the lower bound (\ref{x_lower_bd2_intro}), we apply (\ref{asymp_log_intro}) with $A=2$ to both 
$\pi_{\mathscr{C}} (x,\tilde{K}/\Q)$ and $\pi_{\mathscr{C}} (2x,\tilde{K}/\Q)$. 
If  the (non-negative) difference 
\beq\label{pi_diff_pf}
\pi_\Cscr (2x, \tilde{K}/\Q) - \pi_\Cscr (x, \tilde{K}/\Q) 
\eeq
 were zero, this in combination with (\ref{asymp_log_intro}) would imply that 
\beq\label{Li_cont}
  \mathrm{Li}(2x) -  \mathrm{Li}(x)   \leq \frac{2x}{(\log 2x)^2} + \frac{x}{ (\log x)^2} \leq \frac{3x}{(\log x)^2}.
\eeq
Yet certainly for $x \geq 2$,
\[ \int^{2x}_x \frac{dt}{\log t} \geq \frac{x}{\log 2x}  \geq \frac{x}{2\log x}.\]
Thus (\ref{Li_cont}) fails (so the difference in (\ref{pi_diff_pf}) must be $\geq 1$)  as soon as $x \geq \max \{ 2, e^6\}$. Given $\sig>0$, we apply this to $x= D_K^{\sig}$, in which case we require $D_K \geq D_7 = \max \{D_0', D_6', 2, e^6\}$ with the parameter $D_6'$ (depending on $\sig$) as above.
This completes the verification of Corollary \ref{cor_primes_families}.

\vspace{1cm}

\addcontentsline{toc}{part}{Part III: Applications and $\ell$-torsion in class groups}
\begin{center}
{\bf Part III: Applications}
\end{center}

\section{Bounding $\ell$-torsion in class groups}

For a finite extension $K/\Q$, the ideal class group $\Cl_{K}$  is a finite abelian group
that encodes information about arithmetic in $K$, and interest in  the class number $|\Cl_K|$ has a long history, going back to the Gauss class number conjecture, early attempts at proving Fermat's Last Theorem, and Dirichlet's development of the class number formula, which unites class numbers with $L$-functions. 
We  focus on the $\ell$-torsion subgroup of $\Cl_K$, defined for any integer $\ell \geq 1$ by
\beq\label{torsion_subgroup_dfn}
 \Cl_{K}[\ell]:= \{ [\mathfrak{a}]\in\Cl_{K} : [\mathfrak{a}]^{\ell} = \mathrm{Id} \}.
\eeq
 For any number field $K/\Q$ of degree $n$ and absolute discriminant $D_K = |\Disc K/\Q |$,  we may trivially bound the $\ell$-torsion subgroup by the full class group, which admits the following bound  (see  \cite[Theorem 4.4]{Nar80}):
\beq\label{trivial_bound}
1 \leq |\Cl_K[\ell]| \leq |\Cl_K| \ll_{n,\ep} D_K^{1/2+\ep},
 \eeq
 for any integer $\ell \geq 1$, and $\ep>0$ arbitrarily small. We will refer to this as the trivial bound for $|\Cl_K[\ell]|$.
 
  Our work on $\ell$-torsion is inspired by  the following well-known  conjecture (e.g. see \cite[``Question $\mathrm{CL}(\ell,d)$'']{BruSil96},  \cite{Duk98}, \cite[Conjecture 3.5]{Zha05}):

\begin{conjecture}[$\ell$-torsion Conjecture]\label{conj_class}
\item Let $K/\Q$ be a number field of degree $n$. Then for every integer $\ell \geq 1$ and every $\ep>0$,
 $$|\Cl_{K}[\ell]| \ll_{n,\ell,\ep}D_K^\ep.$$
\end{conjecture}
 
 Now, with our new effective Chebotarev theorems for families of fields, we can make new progress toward this conjecture: we improve on the trivial bound (\ref{trivial_bound}) and in fact do as well as previous bounds that assumed GRH, for all but a possible  density zero subfamily of fields.
In particular, we prove the first unconditional nontrivial upper bounds for $\ell$-torsion, for all $\ell \geq 1$, for almost all fields in infinite families of fields of arbitrarily high degree.

\begin{thm}\label{thm_class_main}
Let $Z_n^\Iscr(\Q,G)$ be fixed to be one of the families of fields considered in Theorems \ref{thm_cheb_main_quant_pt1}, \ref{thm_cheb_main_quant_pt2}, \ref{thm_cheb_main_quant_pt3}, \ref{thm_cheb_main_quant_pt2'} and \ref{thm_cheb_main_quant_pt4}, and correspondingly assume the hypotheses (if any) of the relevant theorem.  Let the parameters $\tau_* < \be \leq d$ be those proved to exist for that family in (\ref{be_d_exist}). For every $\tau> \tau_*$ sufficiently close to $\tau_*$, every $\ep_0>0$ sufficiently small,
and  every integer $\ell \geq 1$, there exists a constant
$D_8$
such that for for every $X \geq 1$, aside from at most $D_8 X^{\tau +\ep_0}$ exceptions, every field $K \in Z_n^\Iscr(\Q,G;X)$ satisfies 
\beq\label{EllVen_bd0_corr1}
|\Cl_K[\ell]| \ll_{n,\ell,G,\ep} D_K^{\frac{1}{2} - \frac{1}{2\ell (n-1)} + \ep}
\eeq
for all $\ep>0$.
\end{thm}
Recalling that for each family considered we have shown that $|Z_n^\Iscr(\Q,G;X)| \gg_{n,G,\Iscr} X^{\be}$ with $\be > \tau_*$, the exceptional family has density zero once $\tau$ is sufficiently close to $\tau_*$ and $\ep_0$ is taken to be sufficiently small.
(In \S \ref{sec_class_avg}, we re-state Theorem \ref{thm_class_main} in terms of \emph{averages} of $\ell$-torsion.)

 The deduction of Theorem \ref{thm_class_main} follows a general approach codified by Ellenberg and Venkatesh for bounding $\ell$-torsion in $\Cl_K$ by finding many small rational primes that split completely in $K$:

\begin{letterthm}[{\cite[Lemma 2.3]{EllVen07}}]\label{lemma_M_primes}
Suppose $K/\Q$ is an extension of degree $n$, and let $\ell$ be a positive integer. Set $0<\del< \frac{1}{2\ell(n-1)}$ and suppose that there are at least $M$ rational primes with $p \leq D_K^\del$ that are unramified and split completely in $K$. Then for any $\ep>0$,
\[ |\Cl_K[\ell]| \ll_{n,\ell,\ep} D_K^{\frac{1}{2} + \ep}M^{-1}.\]
\end{letterthm}

To find small primes that split completely in $K$ it is sufficient to find small primes that split completely in the Galois closure $\tilde{K}$ of $K$ over $\Q$, and to do so Ellenberg and Venkatesh applied Lagarias and Odlzyko's conditional Theorem \ref{thm_LO} to obtain: 

\begin{letterthm}[{\cite[Prop. 3.1]{EllVen07}}]\label{thm_EV_cor}
Let $K/\Q$ be a number field of degree $n$ and $\ell \geq 1$ an integer. Assuming GRH, then for any $\ep>0$,
\beq\label{EllVen_bd0}
|\Cl_K[\ell]| \ll_{n,\ell,\ep} D_K^{\frac{1}{2} - \frac{1}{2\ell (n-1)} + \ep}.
\eeq
\end{letterthm}
The argument in \S \ref{sec_thm_class_main_proof} will show that any quantitative improvement to the exponent obtained in Theorems \ref{lemma_M_primes} and \ref{thm_EV_cor} is expected to be similarly reflected in the exponent obtained in (\ref{EllVen_bd0_corr1}).

As $n,\ell$ grow large, to produce the primes required in Theorem \ref{lemma_M_primes},  we must be allowed to count primes as small as any fixed positive power of $D_K$. 
 This in particular illuminates why previously known lower bounds for $\pi_{\Cscr}(x,L/k)$, such as obtained in the recent work of
Thorner and Zaman \cite{ThoZam16a}, \cite[Eqn. 1.6]{ThoZam16b}, or even the result of Theorem \ref{thm_LagOdl_uncond} (assuming no exceptional zero $\be_0$ exists, or in the setting of Theorem \ref{thm_Stark_Cheb}), do not suffice for our application. 
The new results in Theorem \ref{thm_class_main} show that the following fields satisfy (\ref{EllVen_bd0_corr1}) unconditionally, for all integers $\ell \geq 1$:
\begin{enumerate}[label=(\roman*)]
\item  almost all degree $p$ cyclic extensions of $\Q$ ($p$ prime)
\item  almost all totally ramified  cyclic extensions of $\Q$
\item  almost all degree $p$ $D_p$-extensions ($\Iscr$ the conjugacy class of order 2 elements, odd prime $p$)
\item almost all degree $4$ $A_4$-extensions ($\Iscr$ the two conjugacy classes of order 3 elements).
\end{enumerate} 
Furthermore, Theorem \ref{thm_class_main} shows that for every $n \geq 2$, almost all degree $n$ $S_n$-extensions of $\Q$ with square-free discriminants satisfy (\ref{EllVen_bd0_corr1}) for all $\ell \geq 1$, where this result is 
\begin{enumerate}[label=(\roman*)]
\setcounter{enumi}{4}
\item unconditional if $n=2,3,4$
\item  if $n=5$, conditional on the strong Artin conjecture
\item  if $n \geq 6$, conditional on the strong Artin conjecture and $\Dbf_n(S_n,\varpi_n)$ for  some $\varpi_n < 1/2  +1/n$.
\end{enumerate}
Finally,   Theorem \ref{thm_class_main} shows (among other results for simple groups) that (\ref{EllVen_bd0_corr1}) holds 
\begin{enumerate}[label=(\roman*)]
\setcounter{enumi}{6}
\item for every $n \geq 5$, almost all degree $n$ $A_n$-extensions
of $\Q$ satisfy (\ref{EllVen_bd0_corr1}) for all $\ell \geq 1$, conditional on the strong Artin conjecture.
\end{enumerate}

\begin{remark}
In fact, our proof of Theorem \ref{thm_class_main} works as well if we replace any of our families of fields with the family of their Galois closures.  
\end{remark}

 \subsection{Previous results toward Conjecture \ref{conj_class}}\label{sec_class_previous_results}
To situate our results, we briefly review previous results in the literature toward Conjecture \ref{conj_class} in terms of a property we now define.

   \begin{property}[${\bf C}_{n,\ell}(\Del)$]\label{property_C}
Given integers $n, \ell \geq 1$ and a fixed real number $\Del \geq 0$, we say that property ${\bf C}_{n,\ell}(\Del)$ holds if it is known that for every $\ep>0$ there is a constant $C_{\Del,n,\ell,\ep}$ such that for all fields $K/\Q$ of degree $n$,
 $$|\Cl_{K}[\ell]| \leq C_{\Del,n,\ell,\ep} D_K^{\Del+\ep }.$$
\end{property} 
 Thus in particular, (\ref{trivial_bound}) shows that $\Cbf_{n,\ell}(1/2)$ is trivially true for all $n, \ell \geq 1$.
 The strongest type of result holds for \emph{all} fields of a fixed degree. In this vein,
Gauss \cite{Gau} genus theory shows $\Cbf_{2,2}(0)$ holds.
This is the only case in which Conjecture \ref{conj_class} is known to hold, for a certain prime $\ell$, for all fields of a fixed degree. 
The only other known pointwise bounds for prime $\ell$ are:
 $n=2$ and $\ell=3$, where initial progress occurred in \cite{HelVen06}, \cite{Pie05}, \cite{Pie06}, and \cite{EllVen07} holds the record $\Cbf_{2,3}(1/3)$;
 $\Cbf_{3,3}(1/3)$ due to \cite{EllVen07};
 $ \Cbf_{4,3}(1/2-\del)$ due to \cite{EllVen07} where $\del=1/168$ if $K$ is non-$D_4$;
$\Cbf_{n,2}(0.2784...)$ for $n=3,4$ and $\Cbf_{n,2}(1/2 - 1/2n)$ where $n \geq 5$, due to \cite{BSTTTZ17}. 
Also in \cite{EllVen07}, there is a proof of pointwise bounds for $\ell$-torsion for certain families of fields of arbitrarily high degree, where these fields always contain $\zeta_{\ell}+\zeta_{\ell}^{-1}$.
Conditional on the Birch--Swinnerton-Dyer conjecture and GRH,  Wong \cite{Won99} has observed that $\Cbf_{2,3}(1/4)$ holds.

 For $n=2,3,4,5$,  bounds for $\ell$-torsion at least as strong as (\ref{EllVen_bd0_corr1}) were already known to hold, unconditionally, for \emph{almost all} degree $n$ $S_n$-fields (without any ramification condition). 
 For imaginary quadratic fields, Soundararajan \cite{Sou00} showed that for each prime $\ell$, the nontrivial bound $|\Cl_K[\ell]| \ll_{\ell,\ep} D_K^{1/2-1/2\ell + \ep}$ holds for all but a possible family of exceptional  fields of density zero. 
 Furthermore,
Heath-Brown and the first author \cite{HBP14a} obtained for each prime $\ell \geq 5$ the unconditional bound 
$|\Cl_K[\ell]| \ll_{\ell,\ep}  D_K^{1/2 - 3/(2\ell+2)+\ep}
$ for all but a possible density zero family of imaginary quadratic fields; their methods also yield upper bounds for higher moments of $\ell$-torsion for all $\ell \geq 3$. 
 For each degree $n \leq 5$, Ellenberg and the first and third authors \cite{EPW16} proved the bound (\ref{EllVen_bd0_corr1}) holds unconditionally
 for all but a possible density zero exceptional family of degree $n$ extensions of $\Q$. (In the case $n=4$, this work had the additional requirement that the fields be non-$D_4$ quartic fields and $\ell \geq 8$ and for $n=5$, the requirement $\ell \geq 25$.) In both \cite{HBP14a} and \cite{EPW16}, the upper bound for the possible exceptional family becomes weaker as $\ell$ increases (e.g. in \cite{EPW16} the number of exceptional fields is at most $O_{n,\ep}(X^{1 - 1/(2\ell(n-1)) + \ep})$ for $\ell$ large); this  is noticeably different from the bound for the exceptional set  in Theorem \ref{thm_class_main}.

\begin{remark}
At the time of posting, the authors learned of the works of Frei and Widmer \cite{FreWid17} and Widmer \cite{Wid17}.
Frei and Widmer  obtain the upper bound (\ref{EllVen_bd0_corr1}) for $\ell$-torsion for almost all totally ramified cyclic extensions of $\Q$ (see our case (ii) above), albeit with a larger upper bound for the possible exceptional family of fields, analogous to that in \cite{EPW16}. Frei and Widmer use the sieve method of Ellenberg and the first and third authors \cite{EPW16} combined with new counts for the number of totally ramified cyclic extensions with a finite number of specified local conditions. Notably, their method also works for totally ramified cyclic extensions of any fixed number field $F$. Moreover they remark, building on \cite{Wid17}, on the possibility of sharpening to $1/2 - 1/(\ell(n+1))$ the exponent in (\ref{EllVen_bd0}) for almost all fields in a family $Z_n(\Q,G;X)$ that is sufficiently dense (e.g. $|Z_n(\Q,G;X)| \gg X$). Of the families we consider, the latter strategy could conceivably similarly improve the exponent in (\ref{EllVen_bd0_corr1}) only for the family $Z_n^\Iscr(\Q,S_n;X)$, conditional on such a lower bound being known for the family. We thank Frei and Widmer for sharing their preprint \cite{FreWid17}.
\end{remark}

\subsection{Proof of Theorem \ref{thm_class_main}}\label{sec_thm_class_main_proof}
Theorem \ref{thm_class_main} is an immediate consequence of
Corollary \ref{cor_primes_families}. 
We suppose that a family $Z_n^\Iscr(\Q,G;X)$ and a sufficiently small $\ep_0>0$ have been fixed. We let $0<\del \leq 1/4$ be defined as in (\ref{del_final_thm}).
We set $\Cscr=\{ \mathrm{id}\}$, in which case we are counting primes that split completely in $\tilde{K}$ and hence  in $K$. 
For any integer $\ell \geq 1$,  we take $\tau> \tau_*$ sufficiently close and a sufficiently small $\ep_1>0$ and we set $\sig = 1/(2\ell(n-1)) - \ep_1$. Then for every $X \geq 1$,  for any field $K \in Z_n^\Iscr(\Q,G;X)$ with $D_K \geq D_6$ that is not one of the at most $D_3X^{\tau+\ep_0}$ $\del$-exceptional fields in $Z_n^\Iscr(\Q,G;X)$,  there are $\gg_{G,n,\ell,\ep_1} D_K^{1/2(\ell(n-1)) - \ep_1} / \log D_K$ primes $p \leq D_K^{1/2(\ell(n-1)) - \ep_1}$ that split completely in $K$. Thus for such a  $K$ that is not $\del$-exceptional, by
Theorem \ref{lemma_M_primes},
\beq\label{CL_bound_endgame}
 |\Cl_K[\ell]| \ll_{n,G,\ell,\ep_1,\ep_2}  D_K^{\frac{1}{2} - \frac{1}{2\ell(n-1)} - \ep_1 + \ep_2 },
 \eeq
for all sufficiently small $\ep_1,\ep_2>0$.
Now we count all those fields that are $\del$-exceptional and all those fields in $Z_n^\Iscr(\Q,G;X)$ that have discriminant smaller than $D_6$, of which there are at most $\ll_{n,G,\Iscr} D_6^d$, by the definition of the parameter $d$. Defining $D_8 = D_8(n, \ell, G,\Iscr,d,\tau,\ep_0)$ to be an appropriate maximum of $D_3$ and the above multiple of $D_6^d$, we see  that for every $X \geq 1$ we may say that (\ref{CL_bound_endgame}) holds for each field in $Z_n^\Iscr(\Q,G;X)$, apart from at most $D_8X^{\tau + \ep_0}$ fields. 
This completes the proof of Theorem \ref{thm_class_main}.

\subsubsection{Averages of $\ell$-torsion}\label{sec_class_avg}
The results of Theorem \ref{thm_class_main} can alternatively be stated in terms of averages of $\ell$-torsion over a fixed family of degree $n$ extensions. If 
$|Z_n^\Iscr(\Q,G;X)| \ll_{n,G,\Iscr}X^d,$
 Theorem \ref{thm_class_main}  shows that for all $X \geq 1$, $\ell \geq 1$, $\tau > \tau_*$ sufficiently close, $\ep_0$ sufficiently small,
\[ \sum_{K \in Z_n^\Iscr(\Q,G;X)} |\Cl_K[\ell]| \ll X^{ d + \frac{1}{2} - \frac{1}{2\ell(n-1)} + \ep}
	+ X^{\tau + \frac{1}{2} + \ep_0 + \ep},\]
	for every $\ep>0$,
	with an implied constant depending on $n,\ell,G,\Iscr,d,\tau,\ep_0,\ep$.
For $ \tau_*< \tau < d$, for $\ell$ sufficiently large we will obtain $\tau \leq d - 1/(2\ell(n-1))$, so that
\[ \sum_{K \in Z_n^\Iscr(\Q,G;X)} |\Cl_K[\ell]| \ll X^{ d + \frac{1}{2} - \frac{1}{2\ell(n-1)} + \ep}
	\]
	for every $\ep>0$.
The ``trivial bound'' would be $ \ll_{n,G,\Iscr,\ep} X^{ d + \frac{1}{2}+\ep}$ for all $\ep>0$.

\section{Number fields with small generators}\label{sec_VaaWid_app}

For our second application, we turn to a question of whether all number fields have a ``small'' generator.
Given a number field $K/\Q$ of degree $n$ (inside our fixed algebraic closure $\overline{\Q}$), one can ask for the element $\al \in K$ of smallest height $H(\al)$ such that $K = \Q(\al)$; here $H(\al)$ denotes the absolute multiplicative Weil height. Precisely, for  an element $\al \in K$, 
\[ H(\al) = \prod_{v} \max \{ 1 ,|\al|_v\}^{\frac{d_v}{n}} ,\]
in which $v$ runs over the places of $K$ and for each place $v$, $|\cdot |_v$ is the unique representative that either extends the Archimedean absolute value on $\Q$ or a $p$-adic absolute value on $\Q$, while $d_v = [K_v:\Q_v]$ denotes the local degree at $v$. 
(By Northcott's theorem \cite[Thm.1]{Nor49}, there are finitely many elements in $K$ with height at most any fixed real number, and thus a generator of smallest height does exist.)
 
In terms of lower bounds, it is known by Silverman \cite[Thm. 1]{Sil84}  that for each $n \geq 2$, for all fields $K/\Q$ of degree $n$, 
for any element $\al \in K$ such that $K=\Q(\al)$,
\beq\label{Silverman_lower}
 H(\al) \geq B_1 D_K^{\frac{1}{2 n(n-1)}},
 \eeq
where we may take $B_1=B_1(n)= n^{-\frac{1}{2(n-1)}}$.
In fact, this lower bound  led to the numerology of the savings in the exponent in Theorem \ref{thm_EV_cor}. (See \cite[Lemma 2.2]{EllVen07}, with the lower bound now further explored in the recent preprints \cite{FreWid17,Wid17}, where it is shown that improving on (\ref{Silverman_lower}) for a sufficiently dense class of fields can improve on Theorem \ref{thm_EV_cor} in an average sense.)

On the other hand, regarding upper bounds, Ruppert asked two questions \cite{Rup98} of increasing strength:
\begin{question}\label{ques_Ruppert}
 Does there exist for each $n \geq 2$: 
 \begin{enumerate}
 \item\label{Rup_weak} a positive constant $B_2=B_2(n)$ such that for every field $K/\Q$ of degree $n$ there exists an element $\al \in K$ such that $K=\Q(\al)$ and 
$H(\al) \leq B_2 D_K^{\frac{1}{2n}}?$
\item\label{Rup_strong}
a positive constant $B_3=B_3(n)$ such that for every field $K/\Q$ of degree $n$ there exists an element $\al \in K$ such that $K=\Q(\al)$ and 
$ H(\al) \leq B_3 D_K^{\frac{1}{2n(n-1)}}?$
\end{enumerate}
\end{question}
(Ruppert posed these questions in terms of the naive height, but up to constants this is equivalent to the form given here, for which we cite the presentations of \cite{VaaWid13,VaaWid15}.) 
 The second question is effectively asking whether the exponent in Silverman's lower bound (\ref{Silverman_lower}) is sharp.
For degree $n=2$ the two questions are equivalent, and Ruppert \cite[Prop. 2]{Rup98} answered them in the affirmative. Moreover, \cite[Prop. 3]{Rup98}  verified (\ref{Rup_weak}) for totally real fields $K$ of prime degree. Recently, Vaaler and Widmer \cite[Thm. 1.2]{VaaWid13} verified (\ref{Rup_weak}) for all number fields with at least one real embedding, with a constant $B_2(n) \leq 1$. 
In contrast, they provided in \cite{VaaWid15}, for each composite degree $n$, an infinite family of fields violating (\ref{Rup_strong}). 
Furthermore, in \cite[\S 3 and \S 4]{Wid17}, Widmer shows that for $n \geq 4$, the number of degree $n$ fields satisfying the bound in case (\ref{Rup_strong}) of Question \ref{ques_Ruppert} is $o(X)$, so that the answer to this case must be no. (For clarity, note  that Widmer works in terms of the relative Weil height.) 

This leaves the question of whether case (\ref{Rup_weak}) is true.
As an application of our effective Chebotarev density theorem, we show that within appropriate families of fields, (\ref{Rup_weak}) is true for ``almost all'' fields. 
\begin{thm}\label{thm_VaaWid_app}
Let $Z_n^\Iscr(\Q,G)$ be fixed to be one of the families of fields considered in Theorems \ref{thm_cheb_main_quant_pt1}, \ref{thm_cheb_main_quant_pt2}, \ref{thm_cheb_main_quant_pt3}, \ref{thm_cheb_main_quant_pt2'} and \ref{thm_cheb_main_quant_pt4}, and correspondingly assume the hypotheses (if any) of the relevant theorem.  Let the parameters $\tau_* < \be \leq d$ be those proved to exist for that family in (\ref{be_d_exist}). For every $\tau> \tau_*$ sufficiently close and every $\ep_0>0$ sufficiently small, there exists a constant
$D_9$
such that for every $X \geq 1$, aside from at most $D_9 X^{\tau +\ep_0}$ exceptions, every field $K \in Z_n^\Iscr(\Q,G;X)$ contains an element $\al$ with $K = \Q(\al)$ such that 
$ H(\al) \leq 2 D_K^{\frac{1}{2n}}.$
\end{thm}

The proof  is a simple adaptation of an observation of Vaaler and Widmer in \cite[Thm. 1.3]{VaaWid13}, which relies on finding primes that split completely in $K$ that are of size around $D_K^{1/2}$.  They showed that the bound in Question \ref{ques_Ruppert} case (\ref{Rup_weak}) holds whenever $\zeta_{\tilde{K}}$ satisfies GRH, via an application of Theorem \ref{thm_LO}.
Now, independent of GRH, for every field that is not $\del$-exceptional, with $\del$ determined by (\ref{del_final_thm}),  we  apply part (2) of Corollary \ref{cor_primes_families} with the choices $\Cscr = \{1\}$ and $\sig=1/2$, in place of  \cite[Lemma 5.1]{VaaWid13}. Then  \cite[Thm. 4.1]{VaaWid13} shows that for each field $K$ to which the conclusion (\ref{pi_diff}) applies, there exists an element $\al \in K$ with $K=\Q(\al)$ and 
$ H(\al) \leq p^{1/n} \leq 2D_K^{1/2n}.$
We use Theorem \ref{thm_cheb_main_quant_pt1}  to bound the number of $\del$-exceptional fields, with $\del$ determined by (\ref{del_final_thm}). 
We use the trivial upper bound $|Z_n^\Iscr(\Q,G;D_7)| \ll_{n,G,\Iscr} D_7^d$ for the number of fields in the family with discriminant smaller than the threshold $D_7$ required to apply part (2) of Corollary \ref{cor_primes_families}. Then upon setting $D_9= D_9(n, G,\Iscr,d,\tau,\ep_0)$ to be an appropriate maximum of $D_3$ from Theorem \ref{thm_cheb_main_quant_pt1}  and the above multiple of $D_7^d$, we may then conclude Theorem \ref{thm_VaaWid_app} holds.

\section*{Acknowledgements}
The authors thank P. Sarnak for his prescient advice and consistent encouragement over many years. We thank M. Bhargava for a number of suggestions for quantitative results for families of number fields, and A. Venkatesh and M. Bhargava for indicating a method for counting degree $n$ $A_n$-fields. We thank M. Abel, D. R. Heath-Brown, J. Getz, N. Katz, E. Kowalski, M. Milinovich, M. Nastasescu, D. Ramakrishnan, Z. Rudnick, A. Sutherland, for helpful comments and F. Thorne, A. Zaman for remarks on an earlier version of the manuscript. We also thank the referees for their close reading and helpful comments.
Pierce has been partially supported by NSF DMS-1402121, CAREER grant DMS-1652173, a Sloan Research Fellowship, a Joan and Joseph Birman Fellowship, and as a von Neumann Fellow at the Institute for Advanced Study, by the Charles Simonyi Endowment and NSF Grant No. 1128155. Pierce thanks the Max Planck Institute for Mathematics, the Hausdorff Center for Mathematics, and MSRI (under NSF Grant No. 1440140) for providing focused research environments.
Turnage-Butterbaugh is partially supported by NSF DMS-1901293 and was supported by  NSF DMS-1440140 while in residence at MSRI during the Spring 2017 semester. 
Wood has been partially supported by an American Institute of Mathematics Five-Year Fellowship, a Packard Fellowship for Science and Engineering, a Sloan Research Fellowship, National Science Foundation grant DMS-1301690 and CAREER grant DMS-1652116, and a Vilas Early Career Investigator Award.

%***************************************
\bibliographystyle{alpha}
\bibliography{NoThBibliography}
%***************************************
\label{endofproposal}

\end{document}